\documentclass{amsart}

\usepackage{amssymb}
\usepackage[dvips]{epsfig}
\usepackage{graphicx}
\usepackage{color}

\newcommand{\R}{{\mathbb R}}

\newtheorem{theorem}{Theorem}[section]

\theoremstyle{definition}

\newtheorem*{merci}{Acknowledgements}

\theoremstyle{remark}
\newtheorem{remark}{Remark}[section]
\def\R{{\mathbb R}}

\numberwithin{equation}{section}

\begin{document}
\title[On the Amick-Schonbek system]{Numerical study of the Amick-Schonbek system }

\author[C. Klein]{Christian Klein}
\address{Institut de Math\'ematiques de Bourgogne,  UMR 5584;\\
Institut Universitaire de France \\
Universit\'e de Bourgogne-Franche-Comt\'e, 9 avenue Alain Savary, 21078 Dijon
                Cedex, France} %\\
\email{Christian.Klein@u-bourgogne.fr}

\author[J.-C. Saut]{Jean-Claude Saut}
\address{Laboratoire de Math\' ematiques, UMR 8628,\\
Universit\' e Paris-Saclay et CNRS\\ 91405 Orsay, France}
\email{jean-claude.saut@universite-paris-saclay.fr}

\date{\today}

\maketitle
\begin{center}
In memoriam Vassilios Dougalis (1949-2022)

\end{center}

\begin{abstract}
%\textit{Abstract}
The aim of this paper is to present a survey and a detailed numerical 
study  on a remarkable Boussinesq system  describing weakly 
nonlinear, long surface water waves. In the one-dimensional case, this 
system can be viewed as a dispersive perturbation of the hyperbolic 
Saint-Venant (shallow water) system.   The asymptotic stability of 
the solitary waves is numerically established. Blow-up of solutions 
for initial data not satisfying the non-cavitation condition as well 
as the appearence of dispersive shock waves are studied.

\end{abstract}

\section{Introduction}

%The Euler-Korteweg equation \eqref{EK} contains  also a particular case of 
Our friend Vassilios Dougalis wrote more than ten papers on the 
theory and numerical simulations of a class of long  water wave systems, and we will focus in this paper on one  specific member of this class.

More precisely we  will consider a particular case of the so-called 
abcd Boussinesq systems for surface water waves, see \cite{BCL, 
BCS1,BCS2} \footnote{Boussinesq \cite{Bou2} was the first to derive  
particular Boussinesq systems, not in the class of those studied here 
though. We refer to \cite{Da,Da2,Da3} for details and for an excellent  history of hydrodynamics in the nineteenth century.}
\begin{equation}
    \label{abcd2}
    \left\lbrace
    \begin{array}{l}
    \eta_t+\nabla \cdot {\bf v}+\epsilon \nabla\cdot(\eta {\bf v})+\mu\lbrack a \nabla\cdot \Delta{\bf v}-b\Delta \eta_t\rbrack=0 \\
    {\bf v}_t+\nabla \eta+\epsilon \frac{1}{2}\nabla |{\bf v}|^2+\mu\lbrack c\nabla \Delta \eta-d\Delta {\bf v}_t\rbrack=0.
\end{array}\right.
    \end{equation}
Here $\eta= \eta(x,t), x\in \R^d, d=1,2, t\in \R$ is the elevation of the wave, ${\bf v}={\bf v}(x,t)$ is a measure of the horizontal velocity, $\mu$ and $\epsilon$ are the small parameters (shallowness and nonlinearity parameters respectively) defined as
$$\mu=\frac{h^2}{\lambda^2}, \quad \epsilon= \frac{\alpha}{h}$$
where  $\alpha$ is a typical amplitude of the wave, $h$ a typical depth and $\lambda$ a typical horizontal wavelength.

In the Boussinesq regime, $\epsilon$ and $\mu$ are supposed to be 
small and of same order, $\epsilon\sim\mu\ll1,$ and we will take for simplicity $\epsilon=\mu,$  writing \eqref{abcd2} as \footnote{Particular cases are formally derived in \cite{Broer, Ding, Per}.}
\begin{equation}
    \label{abcd}
    \left\lbrace
    \begin{array}{l}
    \eta_t+\nabla \cdot {\bf v}+\epsilon \lbrack\nabla\cdot(\eta {\bf v})+a \nabla\cdot \Delta{\bf v}-b\Delta \eta_t\rbrack=0 \\
    {\bf v}_t+\nabla \eta+\epsilon\lbrack \frac{1}{2}\nabla |{\bf 
	v}|^2+c\nabla \Delta \eta-d\Delta {\bf v }_t\rbrack=0.
\end{array}\right.
    \end{equation}

The coefficients (a, b, c, d) are restricted by the condition 
$$a+b+c+d=\frac{1}{3}-\tau,$$
where $\tau\geq 0$ is the surface tension coefficient.

When restricted to one-dimensional, unidirectional motions, 
\eqref{abcd} leads to the Korteweg-de Vries (KdV) equation, see 
\cite{La1} (section 9.3)

$$u_t+u_x+\epsilon\left(\frac{1}{3}-\tau\right)u_{xxx}+\epsilon uu_x=0.$$

The class of systems \eqref{abcd2}, \eqref{abcd} models  water waves on a flat bottom propagating in both directions in the aforementioned  regime (see \cite{BCL, BCS1, BCS2}). 

All (well-posed) members of the abcd class provide the same 
approximation of the water wave system in the Boussinesq regime, with an error of order $O(\epsilon^2t)$  (see \cite{BCL}), and  their dispersion relations are similar in the long wave regime. Nevertheless they possess quite different mathematical properties as nonlinear dispersive systems due to the very different behavior of the dispersion relation at high frequencies. Actually the order of their dispersive part  can vary from $+3$ to $-1$, \cite{BCS1, BCS2}. This makes their mathematical study quite fascinating.

It turns out that two particular one-dimensional cases of the abcd 
systems have remarkable properties. We focus here on one of them 
\footnote{The other one, refered to as the Kaup-Broer-Kupperschmidt 
system will be considered in a subsequent paper.} that we will refer 
to as the {\it Amick-Schonbek system} \footnote{Actually this system 
is a particular case of a system derived by Peregrine in \cite{Per}, 
and it is often referred to as the Peregrine or classical Boussinesq 
system, but Schonbek and Amick were the first to recognize its remarkable mathematical properties. A variant with slowly varying bottom is derived in \cite{TW,W}.}  corresponding to $a=b=c=0, d>0,$ say $d=1$  writing thus

\begin{equation}
    \label{AS}
    \left\lbrace
    \begin{array}{l}
    \eta_t+v_x+\epsilon(\eta v)_x=0 \\
    v_t+\eta_x+\epsilon(vv_x-v_{xxt})=0.
\end{array}\right.
    \end{equation}
    
    It appears to be the BBM regularization of the (linearly ill-posed) system
    
    \begin{equation}
    \label{IP}
   \left\lbrace
    \begin{array}{l}
    \eta_t+v_x+\epsilon(\eta v)_x=0 \\
    v_t+\eta_x+\epsilon(vv_x+\eta_{xxx})=0,
\end{array}\right.
    \end{equation}
    which is obtained from the Zakharov-Craig-Sulem formulation of 
	the water wave system after expanding the Dirichlet-to-Neumann operator at first order in $\epsilon.$ From now on we will take $\epsilon =1.$

 There have been many papers dealing with the numerical study of the 
 Amick-Schonbek system, see for instance \cite{An-Dou, An-Dou2}. The 
 case of a non-flat bottom is considered in \cite{GAD} and the 
 periodic case in \cite{ADM3}. Solitary waves were constructed in 
 \cite{chen}, where their interaction was also studied.

 As aforementioned the Amick-Schonbek is the only member of the 
 (one-dimensional) abcd systems for which the global well-posedness 
 of the Cauchy problem with arbitrary large initial data is known. 
 The asymptotics in large time of these solutions is not established, and  
 the simulations in this paper will give some insight into this issue.  
 We will consider in particular the possible soliton decomposition. 
 To this end we study the 1D case with a Fourier spectral method 
 similar to the one for the Serre-Green-Naghdi equation in \cite{DK} and 
 get, see also \cite{chen,An-Dou}, the following\\
 \textbf{Conjecture I:}\\
 The solitary waves of the Amick-Schonbeck system are asymptotically 
 stable. The long time behavior of solutions to this system for 
 Schwartz class initial data is given by solitary waves plus 
 radiation.

In addition we study the behavior of solutions when the 
 initial data $\eta_0$ do not satisfy the non-cavitation 
 condition.   We get\\  
 \textbf{Conjecture II:}\\
 Solutions to initial data not satisfying the non-cavitation 
 condition $1+\eta(x,0)>0$ can blow up in finite time. \\
 We also address
 the zero dispersion limit of the AS system in the form (\ref{ASe}) including the formation of dispersive 
 shock waves (DSWs), zones of rapid modulated oscillations. We 
 conjecture that such DSWs can be observed in the vicinity of shocks 
 of the corresponding dispersionless equation, the Saint Venant 
 system (\ref{SV}), 
 for the same initial data. There will not be a strong limit of 
 these DSWs in the limit $\epsilon\to 0$ in (\ref{ASe}) as in the 
 case of dissipative equations as the Burgers' equation.

    This paper is organized as follows. In  section 2 we review the known theoretical results and open problems on the Cauchy problem for the Amick-Schonbek system.
In section 3 we numerically construct solitary waves. Perturbations 
of these solitary waves are studied in section 4. In section 5 we 
address the long time behavior of hump-like initial data from the 
Schwartz class. In section 6 we consider initial data 
not satisfying the non-cavitation condition. Dispersive shock waves 
are studied in section 7. We add some concluding remarks in section 
8.

   % The second one, referred to as the {\it Kaup-Broer-Kupperschmidt system} ( KBK) corresponds to $a=\pm 1, b=c=d=0$ and writes

    \section{The Cauchy problem}

%\begin{equation}
  %  \label{KBK}
    %\left\lbrace
   % \begin{array}{l}
   % \eta_t+v_x+(\eta v)_x+a v_{xxx}=0 \\
   % v_t+\eta_x+vv_x=0.
%\end{array}\right.
  %  \end{equation}

   % $a=1$ corresponds to the {\it bad} KBK system and $a =-1$ to the {\it good} KBK system.
%\vspace{0.3cm}

%\section{The Amick-Schonbek system}
In this section, we collect some known facts on the Amick-Schonbek 
systems. 

As first noticed  in  a pioneering work by Schonbek \cite{Sch}, the Amick-Schonbek system is among the family of abcd systems the only one that can be viewed as a dispersive perturbation of the Saint-Venant (shallow water) system in the sense that it captures some features of the hyperbolic character of the Saint-Venant system. 

Actually Schonbek \cite{Sch} used  a parabolic regularization of the first equation in \eqref{AS}, to derive from the entropy of the Saint-Venant system an estimate leading to the global well-posedness of the Cauchy problem for arbitrary large  smooth, compactly supported initial data. Amick \cite{Am} weakened the required regularity of the initial data and proved uniqueness. 
Roughly speaking, under a non-cavitation condition that ensures that 
the underlying Saint-Venant system is hyperbolic,  the {\it positive 
convex entropy} of the Saint-Venant system implies an a  priori global bound on the solutions of \eqref{AS} leading to the global existence and uniqueness of solutions without size restrictions. 
%The best known result is obtained in \cite{MTZ} :

%\begin {theorem}\label{MoTaZa}-\cite{MTZ}
%Let $s>1/2.$ For $(\eta_0,v_0)\in H^s(\R)\times H^{s+1}(\R)$ such that $1+\epsilon\eta_0>0,$  \eqref{BBM-EK} has a unique solution $(\eta,v)\in C(\R^+;H^s(\R))\times H^{s+1}(\R))\cap C^1(\R^+,H^{s-1}(\R)\times H^s(\R)).$This solution satisfies $1+\zeta(t,x)>0$ for any $t\geq 0,$ and is the unique solution of \eqref{AS} that belongs to $L^\infty_{\text{loc}}(\R_+:H^s(\R)\times H^{s+1}(\R).$

%\end{theorem}

%One also finds in \cite {MTZ} the existence of low regularity entropic solutions of \eqref{BB       M-EK} emanating from $(\eta_0\in \Lambda_{\sigma_0} , v_0\in H^1(\R))$ where $\Lambda_{\sigma_0}$ is the Orlicz class
%$$
%\Lambda_{\sigma_0}=\lbrace \eta\;\text{measurable}\;|\;\int_\R(\eta\ln \eta-\eta+1)dx<+\infty\rbrace.
%$$

%Those results exploit the hyperbolic nature of the underlying Saint-Venant system. We do not know of similar results in the two-dimensional case.

The Cauchy problem was recently revisited in \cite{MTZ} where the following result is obtained:
\begin{theorem}\label{GWP-AS}
Let $s>1/2.$ For $(\zeta_,u_0)\in H^s(\R)\times H^{s+1}(\R)$ such that $1+\zeta_0>0,$ the Boussinesq system \eqref{AS} has  a  solution $(\zeta,u)\in C(\R_+,H^s(\R)\times H^{s+1}(\R))\cap C^1(\R_+,H^{s-1}(\R)\times H^s(\R)).$  
This solution satisfies $1+\zeta(t,x)>0$ for any $t\geq 0,$ and is the unique solution of \eqref{AS} that belongs to $L^\infty_{\text{loc}}(\R_+:H^s(\R)\times H^{s+1}(\R)).$

For any $T>0,$ the flow-map  $S:(\zeta_0,u_0)\to(\zeta,u)$ is continuous from $H^s(\R)\times H^{s+1}(\R)$ into $C([0,T];H^s(\R)\times H^{s+1}(\R)).$ 

\end{theorem}

As a consequence, Molinet {\it et al} derived the existence of a  global weak solution  $(\zeta,u)\in L^\infty(\R^+;\Lambda_{\sigma_0}\times H^1)$ of the Amick-Schonbek system where $ \Lambda_{\sigma_0}$ is the Orlicz class associated to the entropy:

$$
\Lambda_{\sigma_0}=\lbrace \eta\;\text{measurable}\;|\;\int_\R(\eta\ln \eta-\eta+1)dx<+\infty\rbrace.
$$

These results exploit the hyperbolic nature of the underlying Saint-Venant system. We do not know of similar results in the two-dimensional case (the "long time" existence of solutions to the Cauchy problem was established in \cite{SWX, Bu}).
  
  \begin{remark}
  It is very unlikely that the previous global results persist when 
  the non-cavitation condition is not satisfied by the initial data $\zeta_0.$ Then the Cauchy problem remains of course locally well-posed, but a finite time blow-up  (of which kind?) might occur. 
  %Note that the long time existence result in \cite{SWX, Bu} is obtained under a non-cavitation assumption.
  \end{remark}

  \begin{remark}
The asymptotic behavior of global solutions 
  provided by the previous theorem is unknown. In particular the 
  possible growth in time of the  $H^s\times H^{s+1}$-norm of 
  $(\zeta,u)$ and a lower bound on $(1+\zeta)$ are interesting open 
  issues. The numerical simulations below suggest a soliton 
  resolution property that would imply the boundedness of the 
  $L^\infty$  norms of solutions. No positive lower bound of the 
  infimum of $1+\zeta$ was observed.
  \end{remark}
  
  \begin{remark}
  A derivation of $L^1-L^\infty$ estimates for the linearization at the trivial solution of a large subclass of (abcd) systems, involving the Amick-Schonbek system in spatial dimensions one and two, is provided in \cite{Mel} together with the corresponding Strichartz and Morawetz estimates.
  \end{remark}
  
  \begin{remark}
  The Amick-Schonbek 
  system possesses for any $C>1$ unique solitary wave solutions 
  $V(x-Ct), Q(x-Ct)$ such that $V=V(\xi)$, $\xi=x-C t $, is even, 
  monotonically decreasing when $\xi >0,$ see \cite{chen} and section 
  3 below for details.  As far as we know the stability properties of 
  those solitary waves are unknown, and the numerics in section 4 
  will give some insight into this issue and into their possible role in the long time dynamics of solutions.
  \end{remark}
\begin{remark}
A viscous regularization of the Amick-Schonbek system different from the one used in \cite{Am, Sch, MTZ} is proposed in \cite{BMT}, namely 

\begin{equation}
    \label{AS-visco}
    \left\lbrace
    \begin{array}{l}
    \zeta_t+u_x+(\zeta u)_x=0 \\
    u_t+\zeta_x+uu_x-u_{xxt}-\nu u_{xx}=0,
\end{array}\right.
    \end{equation}

    for which the existence and uniqueness (up to translations) of 
	traveling wave solutions with non zero limit at $-\infty $ is established. They are categorized into dispersive or regularized shock waves  depending on the balance between the dispersive and dissipative effects. 
\end{remark}

\begin{remark}

In \cite{Ada} Adamy extends the Amick-Schonbek results on the Cauchy problem to  the initial value boundary problem  posed on $\R^+$ or a finite interval. She obtains in particular the existence of a weak entropy solution for $x\in \R^+$ with homogeneous or non-homogeneous Dirichlet condition at $x=0$ and similar results on a bounded interval. She also proves the uniqueness of smooth solutions.
\end{remark}

\begin{remark}
Fokas and Pelloni in \cite{FP} solved the linearized Amick-Schonbek system on either the half-line or finite interval with suitable boundary conditions,  following the unified method of Fokas \cite{Fok}. The case of the half-line has been revisited in \cite{JGM}.

\end{remark}

The long time behavior of the global solutions of the Amick-Schonbek 
system is unknown, and we present in the 
section 5  numerical simulations suggesting relevant conjectures.

%\section{Large time behavior of solutions of the Amick-Schonbek system}

%\subsection{Scattering of small solutions, soliton resolution, finite time blow-up,....}
%\subsection{The zero dispersion limit}

\subsection{The two-dimensional case}

The two-dimensional Amick-Schonbek system writes (we have kept the small parameter $\epsilon$):
\begin{equation}
    \label{2D}
    \left\lbrace
    \begin{array}{l}
    \eta_t+\nabla \cdot {\bf v}+\epsilon \lbrack\nabla\cdot(\eta {\bf v})=0 \\
    {\bf v}_t+\nabla \eta+\epsilon\lbrack \frac{1}{2}\nabla |{\bf v}|^2-d\Delta {\bf v }_t\rbrack=0.
\end{array}\right.
    \end{equation}

We are not aware of a connection between \eqref{2D} and the two-dimensional Saint-Venant system. The  existence  on long time $O(1/\epsilon)$ of solutions to the Cauchy problem for \eqref{2D} was established in \cite{SWX}, Theorem 4.4, under the non-cavitation condition $1+\eta_0\geq H_0>0$ and in \cite{Bu} without this condition.
\begin{remark}
The long time existence for a two-dimensional Amick-Schonbek system in presence of a non-trivial bathymetry is proven in \cite{Meso}.
\end{remark}

\section{Solitary waves of the Amick-Schonbek system}
 In this section, we will numerically construct localised traveling 
 wave solutions of the system (\ref{AS}), i.e., solutions of the form 
 $\eta=Q(x-Ct)$, $v = V(x-Ct)$ where $C$ is a constant real velocity. 
 We follow the presentation in \cite{chen}, where the solitary waves 
 were also studied numerically, see also \cite{An-Dou}. 
 This ansatz implies for (\ref{AS}) after one integration taking care of the 
vanishing conditions at infinity,
\begin{align}
	-CQ+V+QV & =0,
	\nonumber\\
	-CV+Q+\frac{1}{2}V^{2}+cV'' & =0,
	\label{ASsol}
\end{align}
where the prime denotes derivative with respect to the argument. Thus 
we get 
\begin{align}
	Q & =\frac{V}{C-V},
	\nonumber\\
	-CV+\frac{V}{C-V}+\frac{1}{2}V^{2}+CV'' & =0,
	\label{ASV}
\end{align}
Thus solitary waves of the Amick-Schonbek system are given by a 
single ODE for $V$, the second equation in (\ref{ASV}). In order to 
get solutions that are exponentially decreasing as $e^{-\alpha|x|}$ 
with $\alpha>0$
for $|x|\to\infty$, one must have $\alpha = \sqrt{1-1/C^{2}}$, i.e., 
$|C|>1$. Furthermore if $V_{C}$ is 
a solution to (\ref{ASV}), so is $-V_{-C}$, i.e., a change of sign 
of $C$ implies a change of sign of $V$. Therefore we will concentrate 
in the  following on the case $C>1$.
Assuming that the function and its derivatives vanish 
at infinity, the second  equation in (\ref{ASV})
can be integrated once more to give
\begin{equation}
	\frac{C}{2} ((V')^{2}-V^{2})+\frac{1}{6}V^{3}-V-C\ln(1-V/c)=0
	\label{ASVint},
\end{equation}
which implies that $V$ can be given in terms of quadratures. Due to 
the appearence of the logarithm, it does not seem possible to give 
$V$ in terms of elementary functions. Note that $Q$ gets very large 
once the maximum of $V$ is close to $C$. 
 
It is possible to numerically compute the integral with respect to $V$, but 
since we need the solution in the form $V(x)$, this has to be done 
for all interesting values of $x$ in order to be numerically able to 
invert the function $x(V)$. Therefore we use here the same approach 
as in \cite{KS15}. Since the solitary waves are rapidly decreasing, 
we can treat  them as periodic functions that are smooth on a 
sufficiently large torus within the finite numerical precision (we 
work with double precision where the smallest difference between 
rational numbers is of the order of $10^{-16}$). Thus we consider 
$x\in L[-\pi,\pi]$ with $L\gg1$ and use the standard discretisation of 
the discrete Fourier transform (DFT) $x_{n} = -\pi L +nh$, $n=1,2,\ldots,N$ 
with $h = 2\pi L/N$. The DFT can be conveniently treated with a fast 
Fourier transform (FFT). We denote the DFT of a discretised function 
$V$  by 
$\hat{V}$. Thus we approximate the second ODE in (\ref{ASV})
\begin{equation}
	F:=-\hat{V}+\frac{1}{C(1+k^{2})}\left(\frac{1}{2}\widehat{V^{2}}+\widehat{\frac{V}{C-V}}\right)=0
	\label{ASF},
\end{equation}
which means we have to find the zero of an $N$-dimensional function, 
where $N$ is the number of DFT modes. This zero is identified 
iteratively with a Newton method where the action of the inverse of 
the Jacobian is computed as in \cite{KS15} with the Krylov subspace 
technique GMRES \cite{GMRES}. 

\begin{figure}[htb!]
 \includegraphics[width=0.49\textwidth]{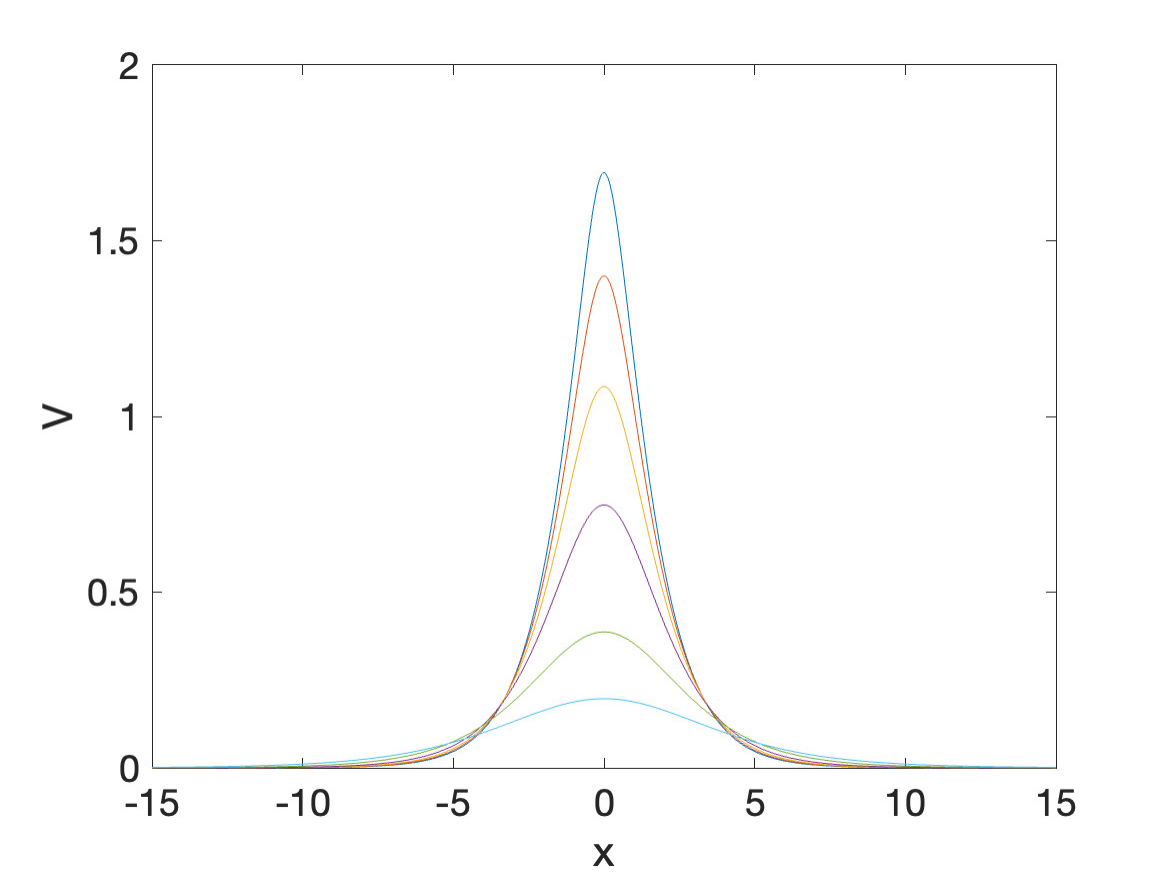}
 \includegraphics[width=0.49\textwidth]{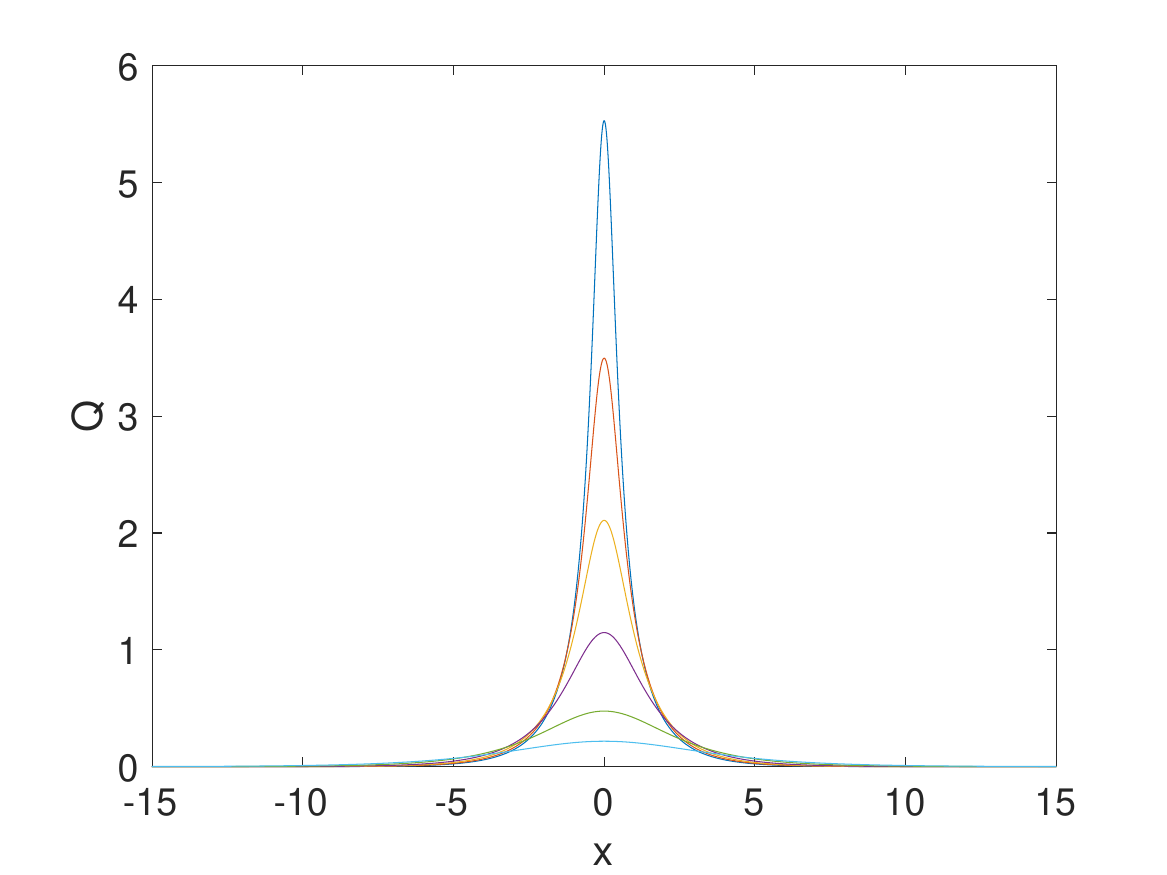}
 \caption{Solitary waves of the Amick-Schonbeck system for different 
 values of the velocity $C=2,1.8,1.6,1.4,1.2,1.1$ from top to bottom, 
 on the left $V$, on the right $Q$.}
 \label{figASsol}
\end{figure}

For concrete computations, we use for $C=2$, $N=2^{12}$ DFT modes and 
$x\in 10[-\pi,\pi]$ with the initial iterate 
$V^{(0)}=1.5\mbox{sech}^2(\alpha x/2)$ where 
$\alpha=\sqrt{1-1/C^{2}}$. For smaller or larger values of $C$ we apply a 
continuation technique, i.e., we use the result for $C=2$ as the 
initial iterate for slightly larger or smaller values of $C$. We show 
the solitary waves for different velocities in Fig.~\ref{figASsol}. 
It can be seen that the amplitude of the waves decreases with $C$, 
and that the fall-off to infinity becomes slower as expected. It 
appears that the amplitude of the waves tends to zero for $C\to1$. 

\begin{figure}[htb!]
 \includegraphics[width=0.49\textwidth]{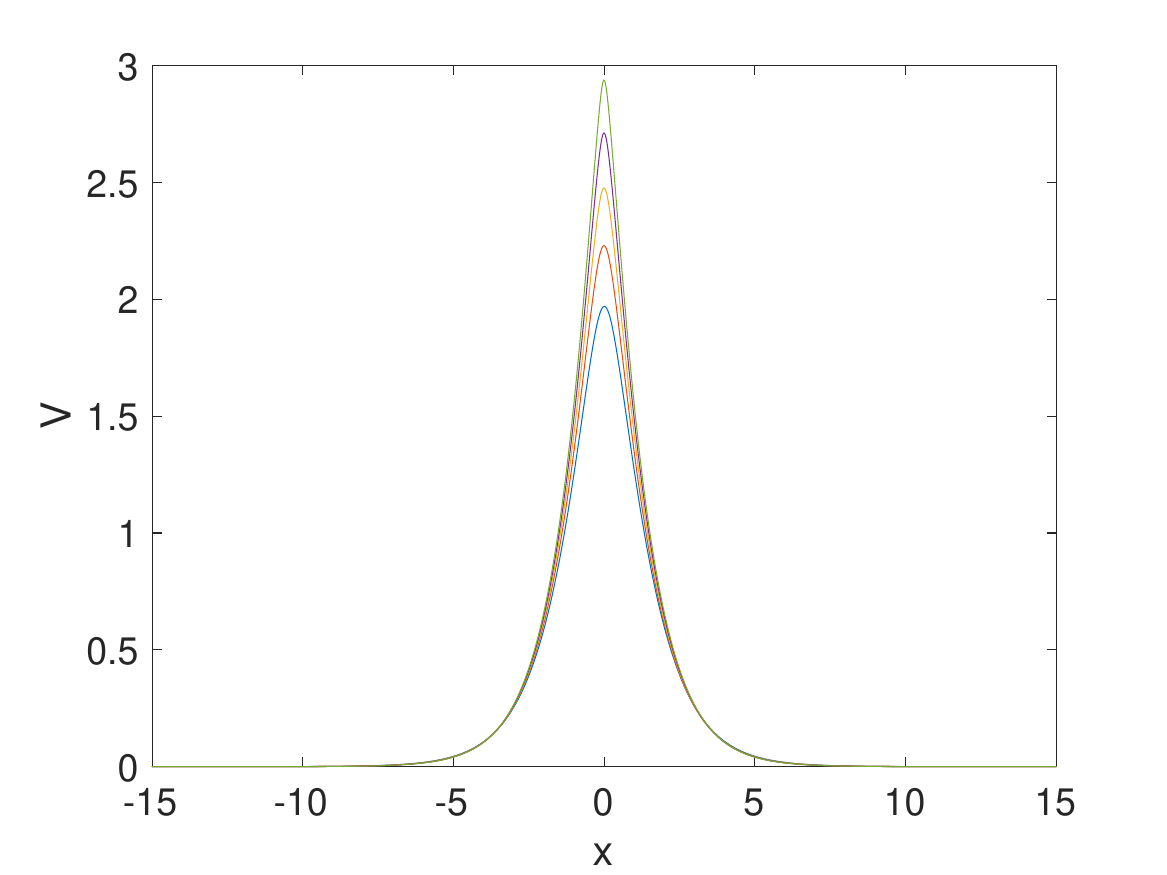}
 \includegraphics[width=0.49\textwidth]{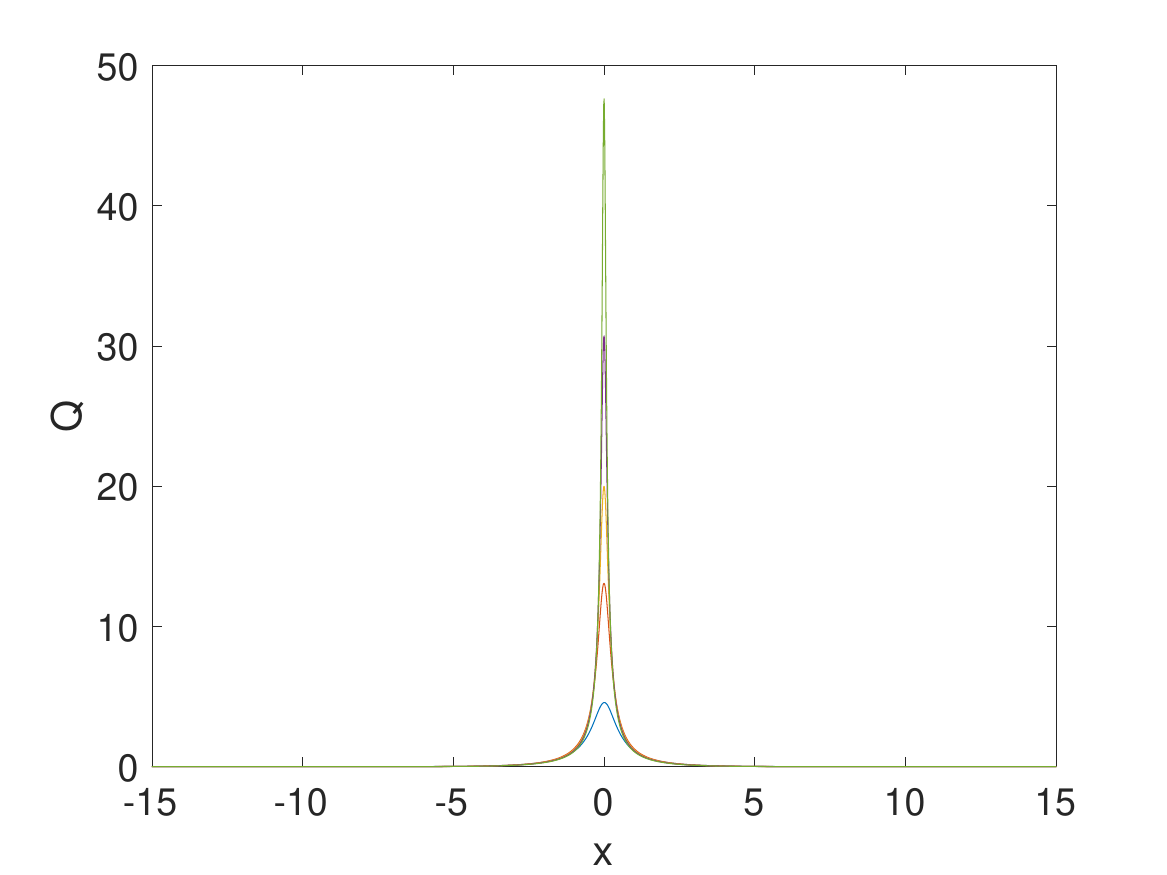}
 \caption{Solitary waves of the Amick-Schonbeck system for different 
 values of the velocity $C=3, 2.8, 2.6,  2.4, 2.2$ from top to bottom, 
 on the left $V$, on the right $Q$.}
 \label{figASsolc23}
\end{figure}

For values of $C$ larger than 2, the maximum of $V$ tends closer to 
$C$ as can be seen in Fig.~\ref{figASsolc23}. This implies that the 
function $Q$ grows strongly. But there is no indication that there is 
an upper limit for $C$ beyond which there are no solitary waves. 

\begin{figure}[htb!]
 \includegraphics[width=0.49\textwidth]{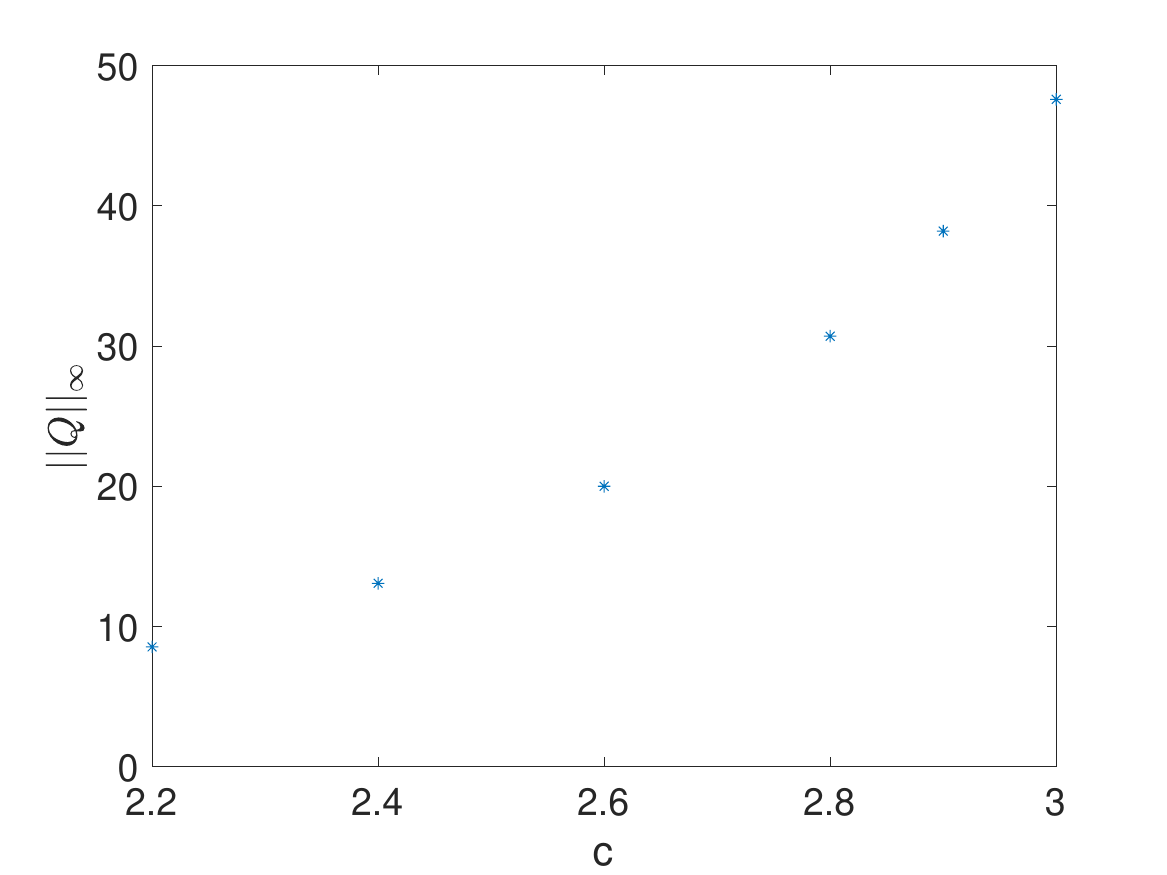}
 \includegraphics[width=0.49\textwidth]{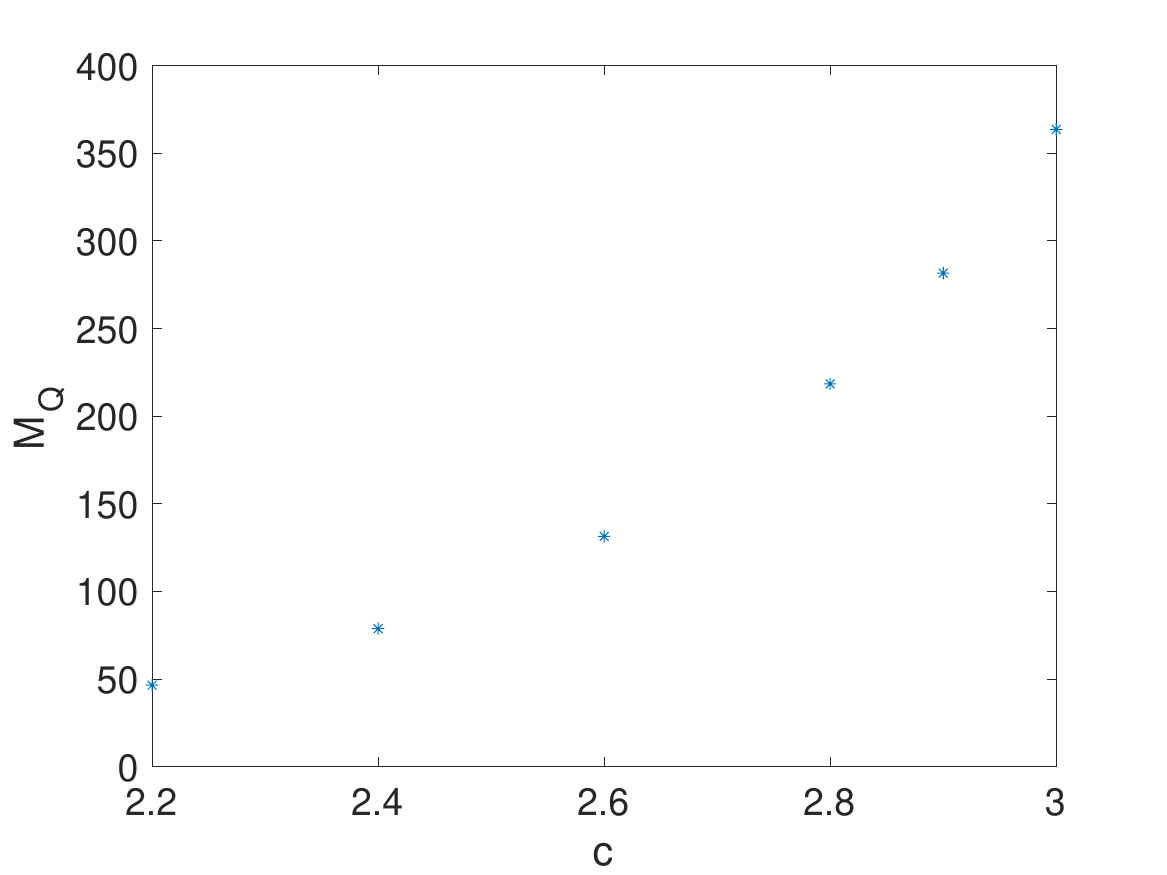}
 \caption{The $L^{\infty}$ norm of the functions $Q$ of 
 Fig.~\ref{figASsolc23}.}
 \label{figASsolc23max}
\end{figure}

We show in Fig.~\ref{figASsolc23max} the $L^{\infty}$ norm of the functions $Q$ and the related 
mass (the square of the $L^{2}$ norm) for the considered values of 
$C$. Both appear to grow algebraically, but it is difficult to access 
much larger values of $C$ in order to make a numerical conjecture.

\section{Numerical study of perturbed solitary waves}
In this section, we numerically study the stability of the solitary 
waves constructed in the previous section. The results give strong 
evidence for the first conjecture in the introduction. For the 
interaction of solitary waves the reader is referred to 
\cite{chen,An-Dou}. 

For a numerical time evolution of solutions of the 
Amick-Schonbeck system (\ref{AS}), we apply the same DFT 
discretisation as for the construction of the solitary waves. The 
time integration is done with the standard explicit 4th order 
Runge-Kutta method. Since there are apparently no conserved 
quantities for the system (\ref{AS}), which could be used to control 
the accuracy of the time integration, we test the code at the example 
of  the 
solitary waves. For the case $C=2$ in Fig.~\ref{figASsol}, we use 
$N_{t}=4000$ time steps for $t\leq 1$ and the same spatial resolution 
as in the previous section. For $t=1$, the difference between the 
numerical solution to the system (\ref{AS}) for solitary wave initial 
data and the propagation of the initial data with velocity $C=2$ is 
smaller than $10^{-13}$ on the whole computational domain. 
Thus the initial data can be propagated with 
essentially machine precision. Note that this also tests the 
numerically constructed solitary waves since a non-trivial error there would 
show up in the time evolution even for arbitrary small time steps. 

\begin{figure}[htb!]
 \includegraphics[width=0.49\textwidth]{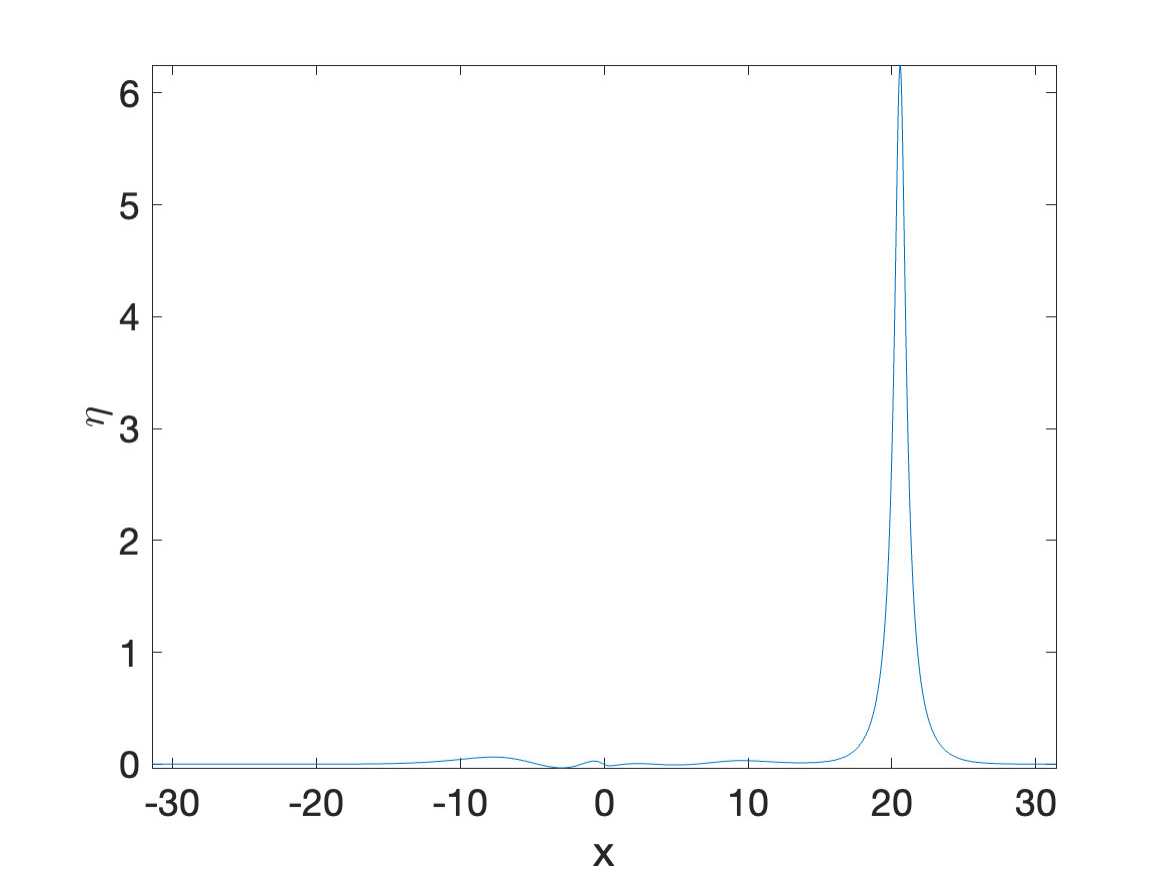}
 \includegraphics[width=0.49\textwidth]{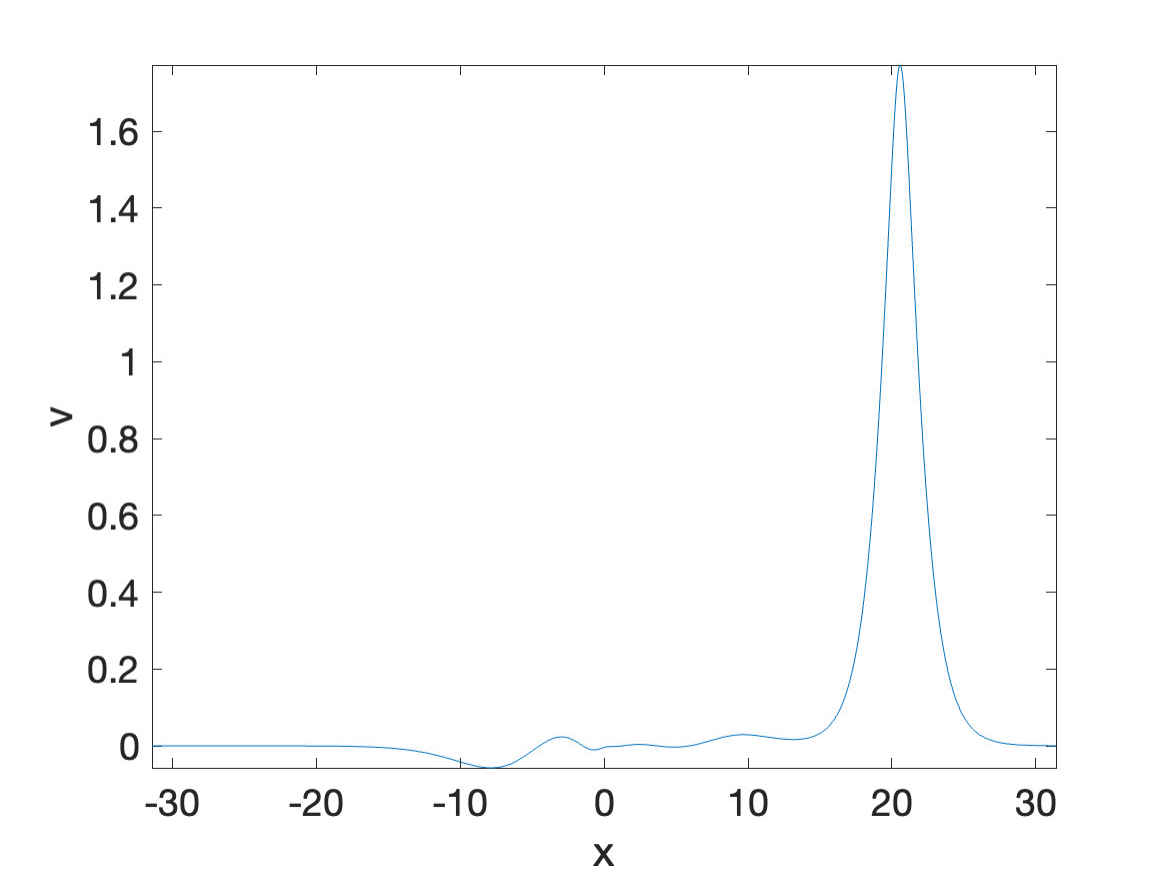}
 \caption{Solution to the Amick-Schonbeck system (\ref{AS}) for the 
 initial data $\eta(x,0) = 1.1 Q_{2}(x)$, $v(x,0) = V_{2}(x)$ for 
 $t=10$, on the left $\eta$, on the right $v$.}
 \label{figASc2la11}
\end{figure}

To consider perturbations of the solitary waves, we use initial 
data for the system (\ref{AS}) which are perturbed solitary waves. We 
start with initial data of the form 
\begin{equation}
	\eta(x,0) = \lambda Q(x),\quad v(x,0) = V(x)
	\label{lambda}
\end{equation}
where $\lambda\in \mathbb{R}$ with $\lambda\sim 1$. In practice we 
consider values of $\lambda$ between $0.9$ and $1.1$.  Note that 
these are
perturbations of the order of 10\% and thus by no means small. 
In a numerical context, one always has to consider  
perturbations of a certain finite magnitude in order to see effects of the 
perturbation on the solution in finite time. We use 
$N_{t}=10^{4}$ time steps for $t\leq10$ for $\lambda=1.1$. The solution at the final 
time is shown in Fig.~\ref{figASc2la11}. It can be seen that 
radiation is propagating to the left as a consequence of the 
perturbation, but that the initial hump is essentially unchanged.

\begin{figure}[htb!]
 \includegraphics[width=0.49\textwidth]{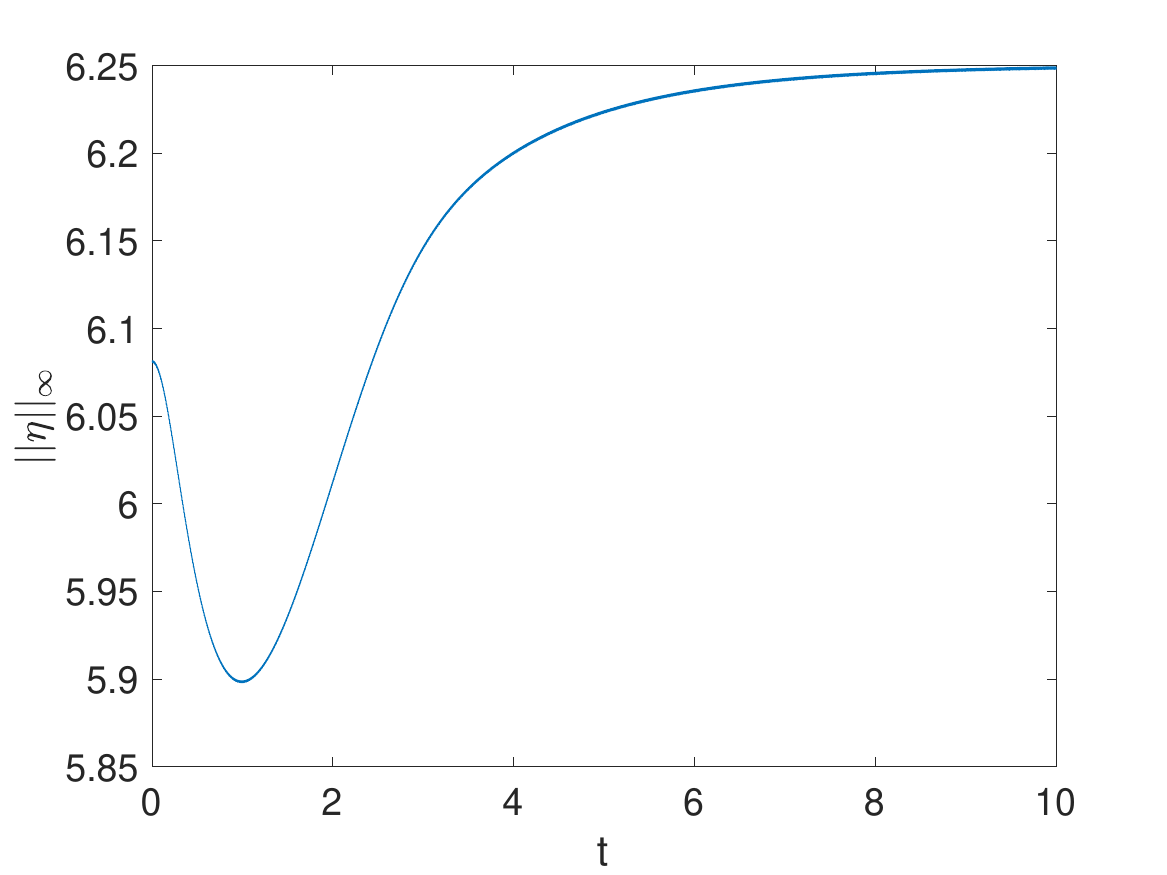}
 \includegraphics[width=0.49\textwidth]{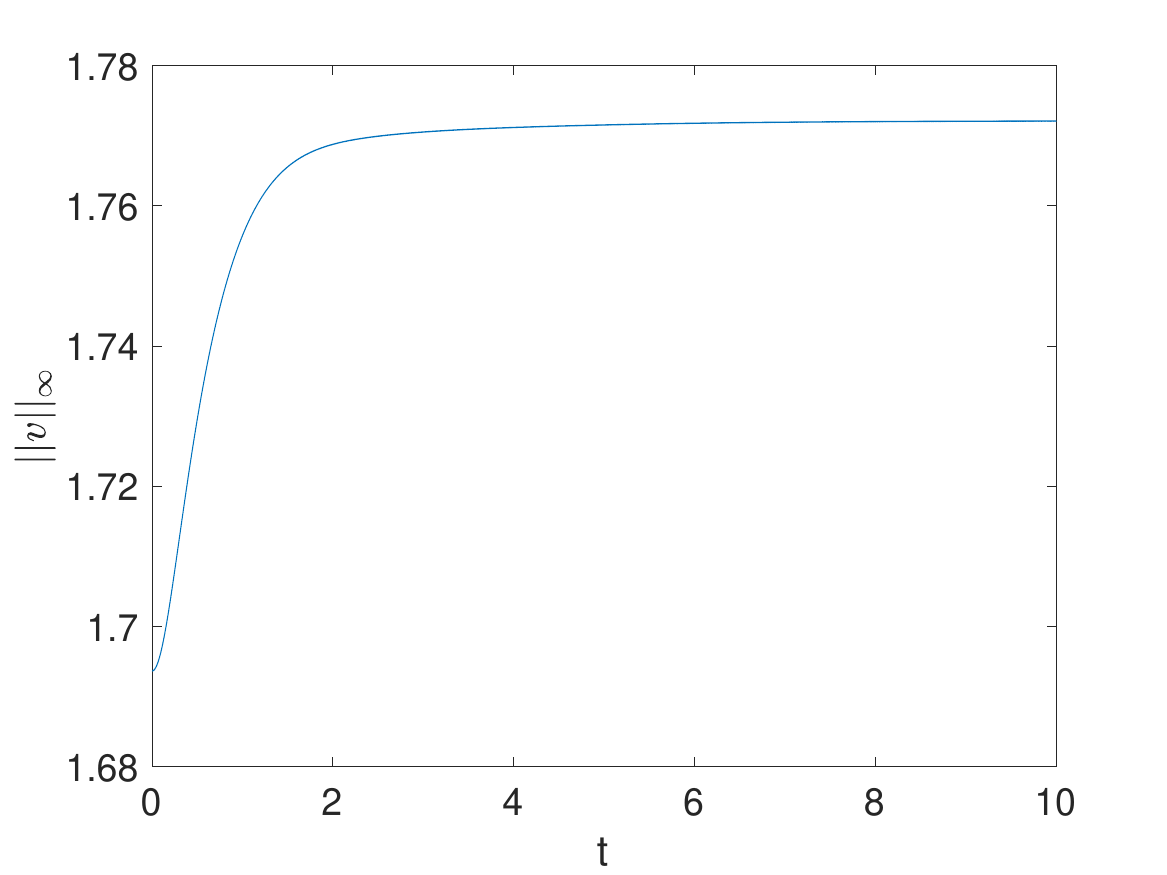}
 \caption{$L^{\infty}$ norms of the solution to the Amick-Schonbeck system (\ref{AS}) for the 
 initial data $\eta(x,0) = 1.1 Q_{2}(x)$, $v(x,0) = V_{2}(x)$ in dependence of time, 
 on the left for $\eta$, on the right for $v$.}
 \label{figASc2la11inf}
\end{figure}

The interpretation of the humps as a solitary wave of slightly larger 
velocity (due to the large perturbation, the final state will be a 
slightly different solitary wave) is confirmed by the $L^{\infty}$ 
norms of $\eta$ and $v$ in Fig.~\ref{figASc2la11inf}. Both norms 
appear to reach asymptotically a constant value at a slightly higher 
level than for the unperturbed solitary wave. Note that the solitary 
waves do not show a simple scaling in $C$. Thus it is not obvious 
with which velocity the final state in Fig.~\ref{figASc2la11inf} 
propagates. The small oscillations in the $L^{\infty}$ norm are due 
to the norm being evaluated on the numerical grid whilst the maximum 
does not have to be located on a grid point, and that radiation can 
reenter the computational domain on the other side since we work on a 
torus. 

\begin{figure}[htb!]
 \includegraphics[width=0.49\textwidth]{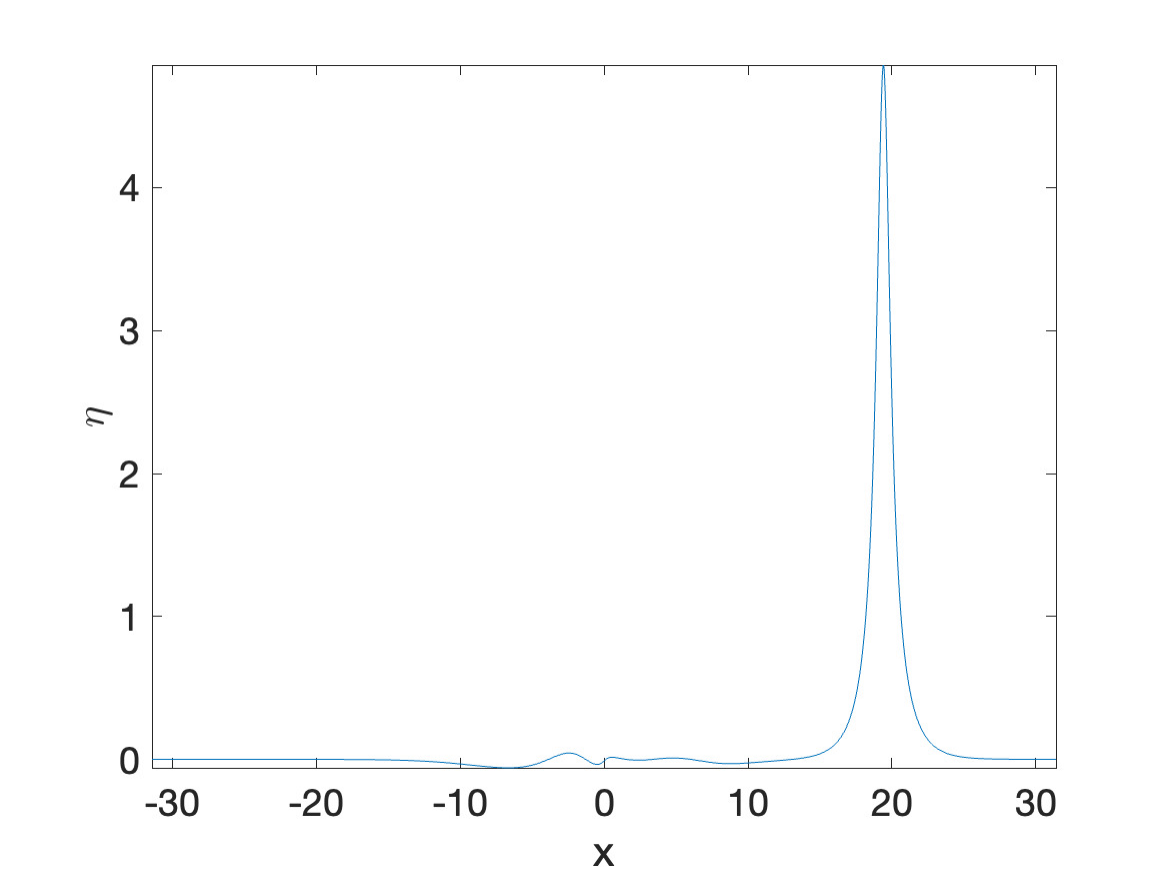}
 \includegraphics[width=0.49\textwidth]{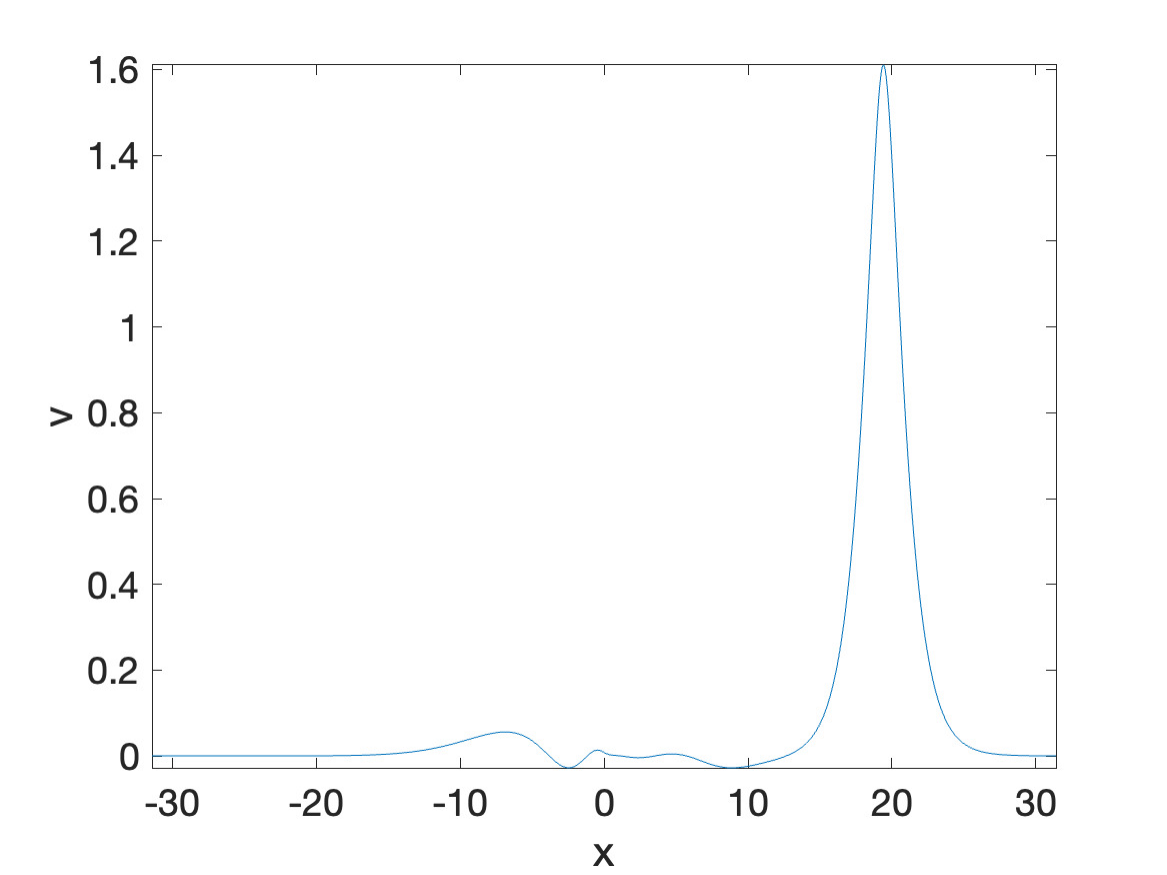}\\
  \includegraphics[width=0.49\textwidth]{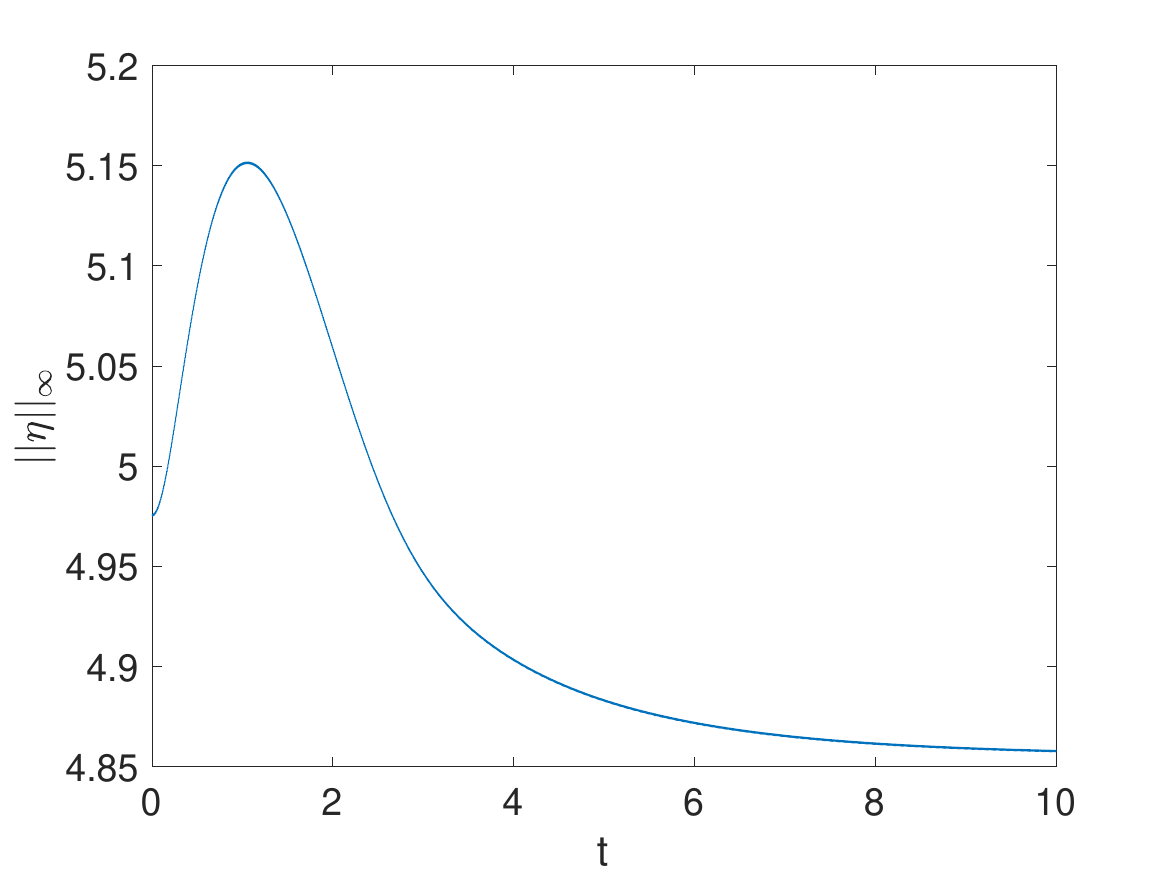}
 \includegraphics[width=0.49\textwidth]{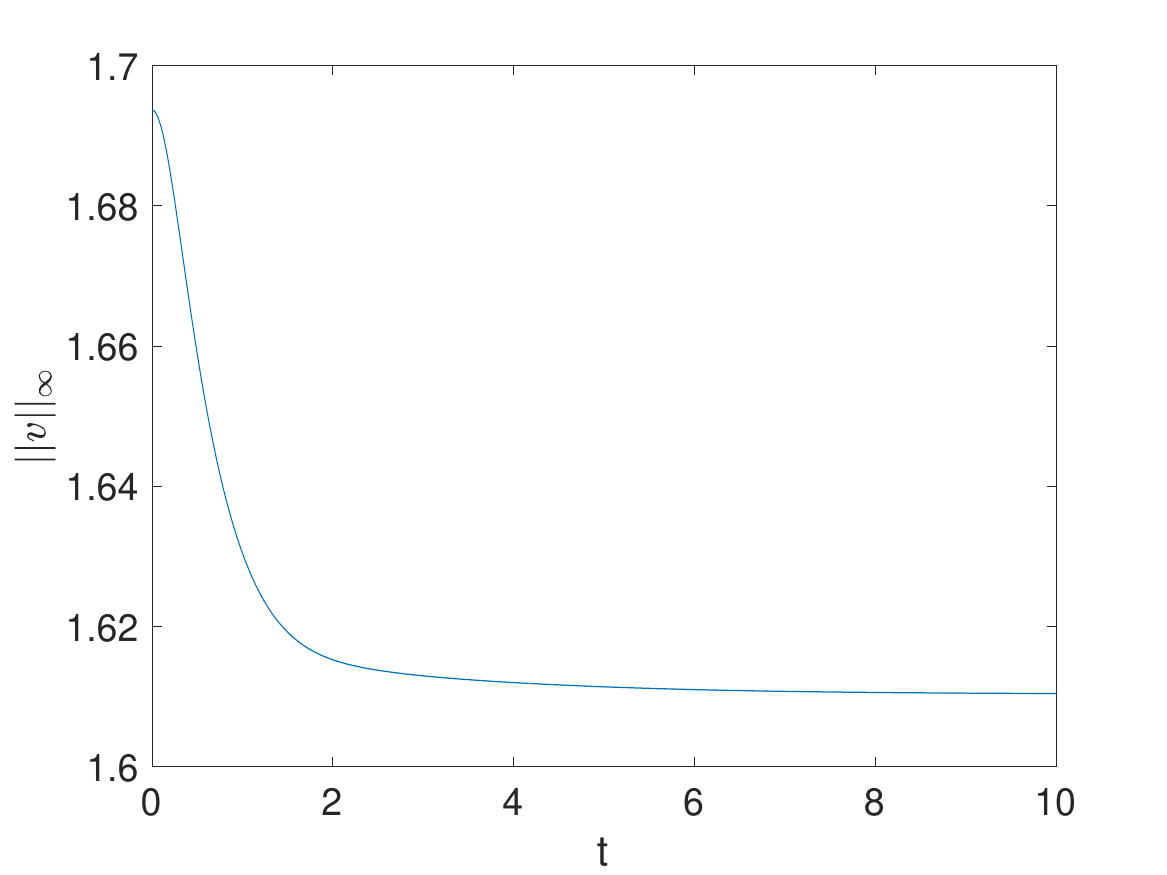}
 \caption{Solution to the Amick-Schonbeck system (\ref{AS}) for the 
 initial data $\eta(x,0) = 0.9 Q_{2}(x)$, $v(x,0) = V_{2}(x)$ for 
 $t=10$, on the left $\eta$, on the right $v$, in the upper row the 
 solutions for $t=10$, in the lower row the respective $L^{\infty}$ 
 norms in dependence of time.}
 \label{figASc2la09}
\end{figure}

The results are similar for initial data of the form (\ref{lambda}) 
with $\lambda=0.9$. The solutions for $t=10$ can be seen in the upper 
row of Fig.~\ref{figASc2la09}. Once more they appear to be solitary 
waves plus radiation, but this time with slightly smaller velocity 
than the unperturbed solitary wave. This interpretation is confirmed 
by the $L^{\infty}$ norms shown in the lower row of the same figure. 

Thus the solitary waves for $C=2$ appear to be stable for the 
considered perturbation. In order to show that this result does not 
depend on a specific perturbation, we study also perturbations of 
the form
\begin{equation}
	\eta(x,0) = Q_{c}(x)\pm \mu \exp(-x^{2}),\quad v(x,0) = V_{c}(x),
	\label{gausspert}
\end{equation}
i.e., we add or subtract a small Gaussian to the solitary wave. Once 
more the solution at the final time appears to be a solitary wave of 
slighty larger (for the $+$ sign) or smaller (for the $-$ sign) 
velocity. 
\begin{figure}[htb!]
 \includegraphics[width=0.49\textwidth]{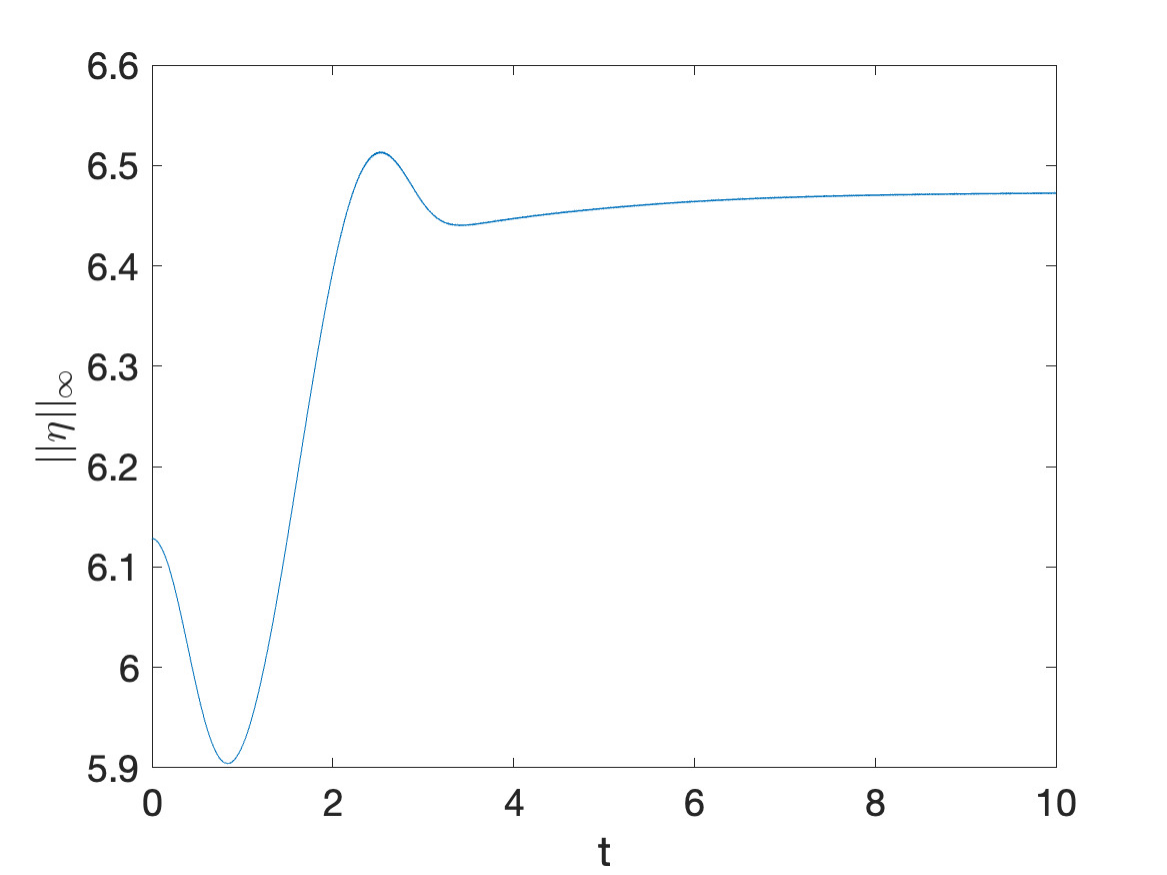}
 \includegraphics[width=0.49\textwidth]{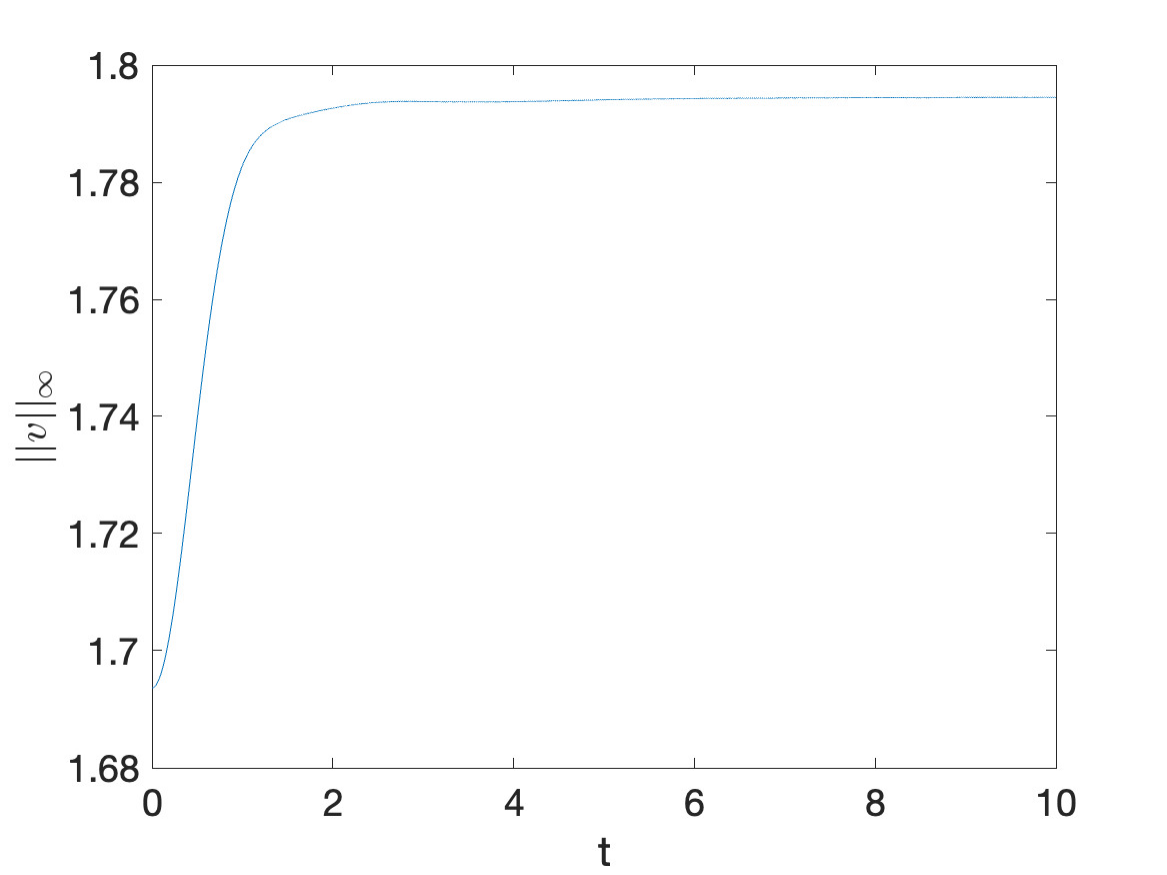}\\
  \includegraphics[width=0.49\textwidth]{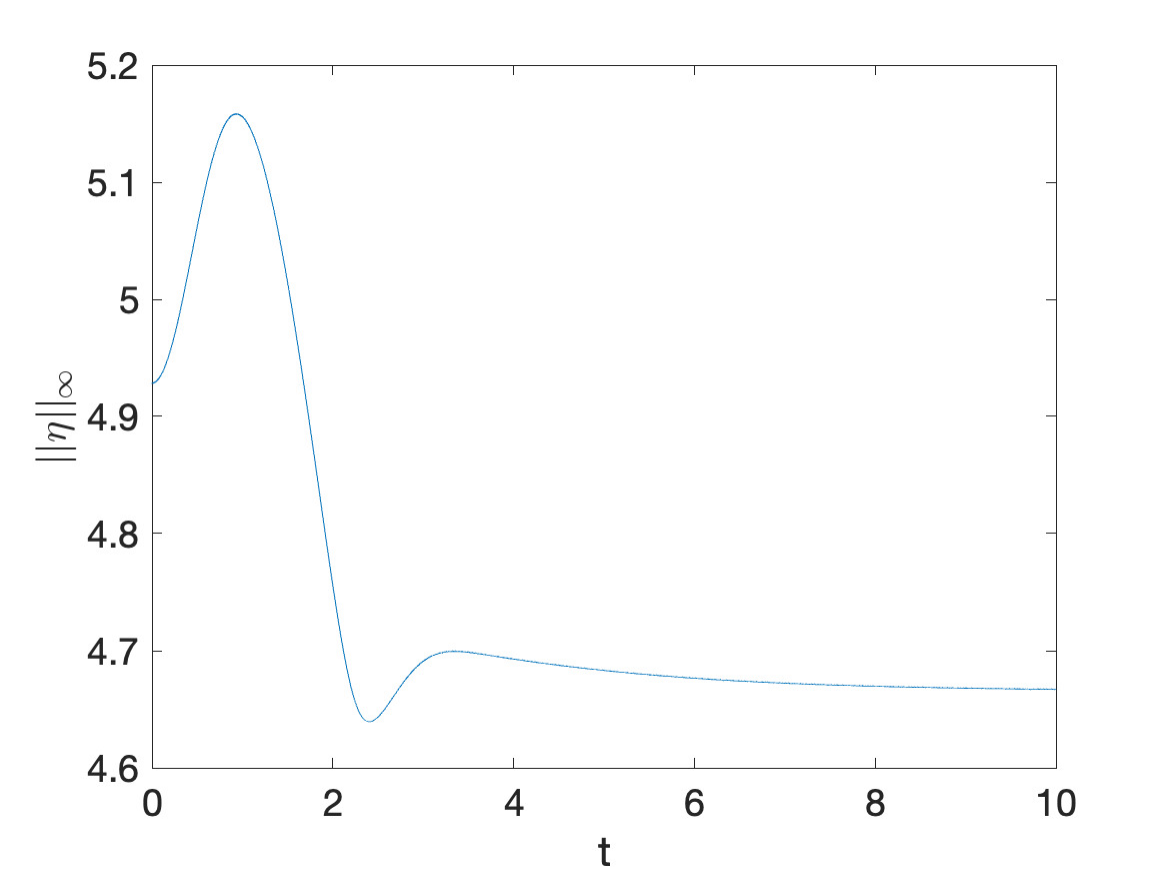}
 \includegraphics[width=0.49\textwidth]{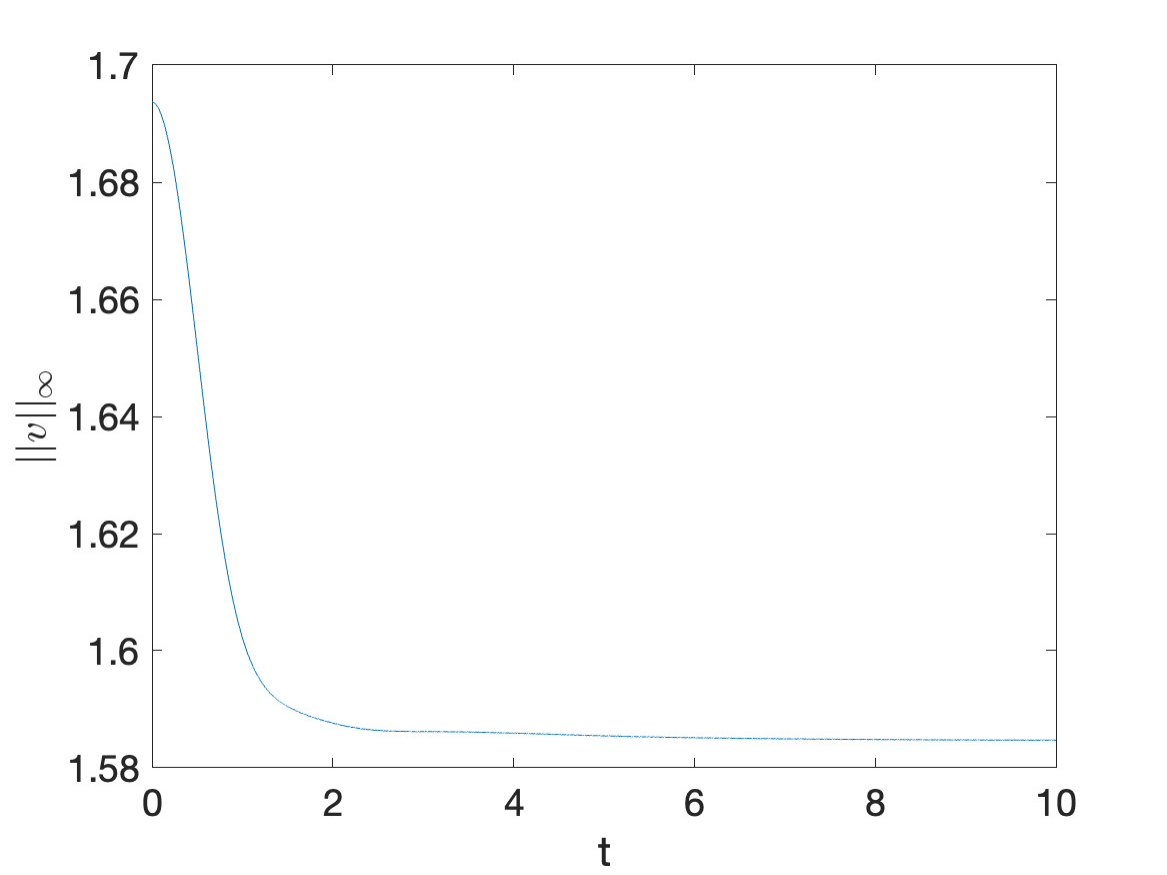}
 \caption{$L^{\infty}$ norms of the solution to the Amick-Schonbeck system (\ref{AS}) for the 
 initial data $\eta(x,0) = Q_{c}(x)\pm 0.6 \exp(-x^{2})$, $v(x,0) = 
 V_{2}(x)$, on the left $\eta$, on the right $v$, in the upper row 
 for the $+$ sign, in the lower row for the $-$ sign.}
 \label{figASc2gauss}
\end{figure}

\begin{figure}[htb!]
 \includegraphics[width=0.49\textwidth]{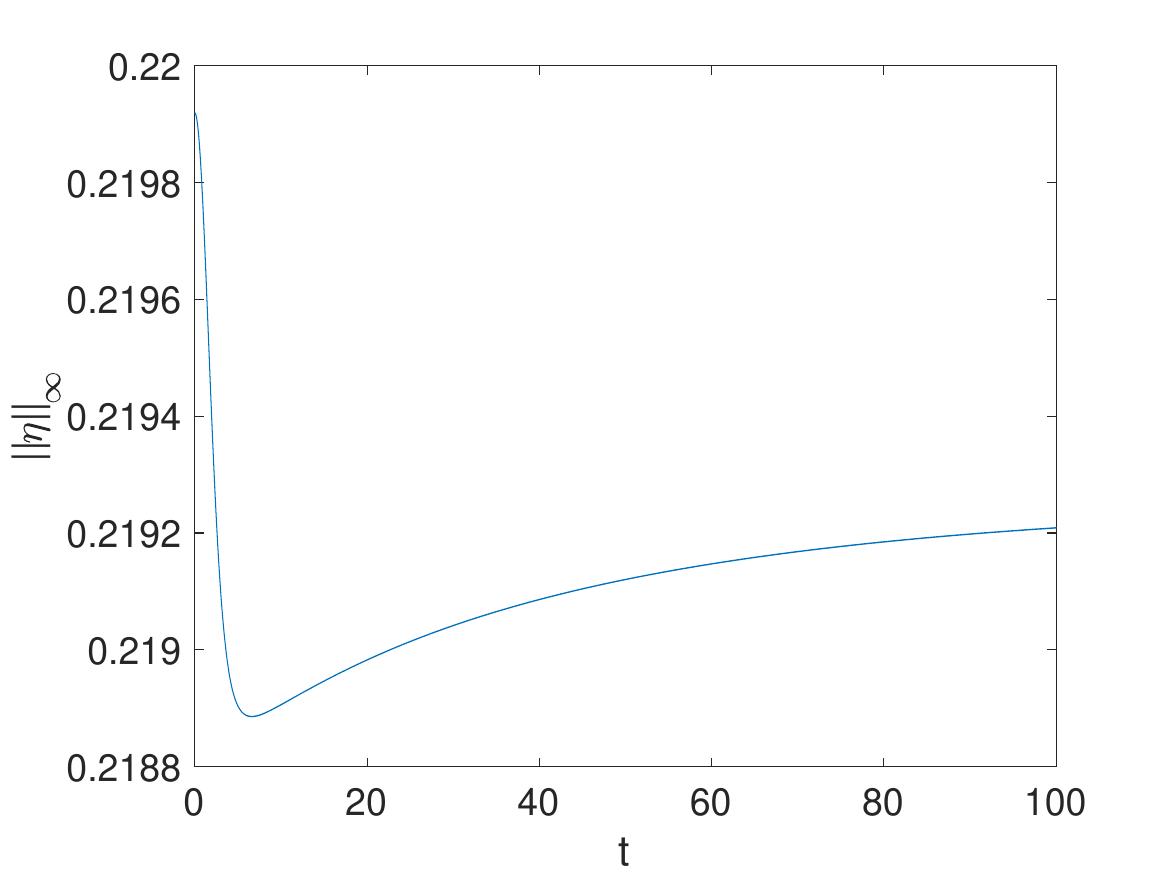}
 \includegraphics[width=0.49\textwidth]{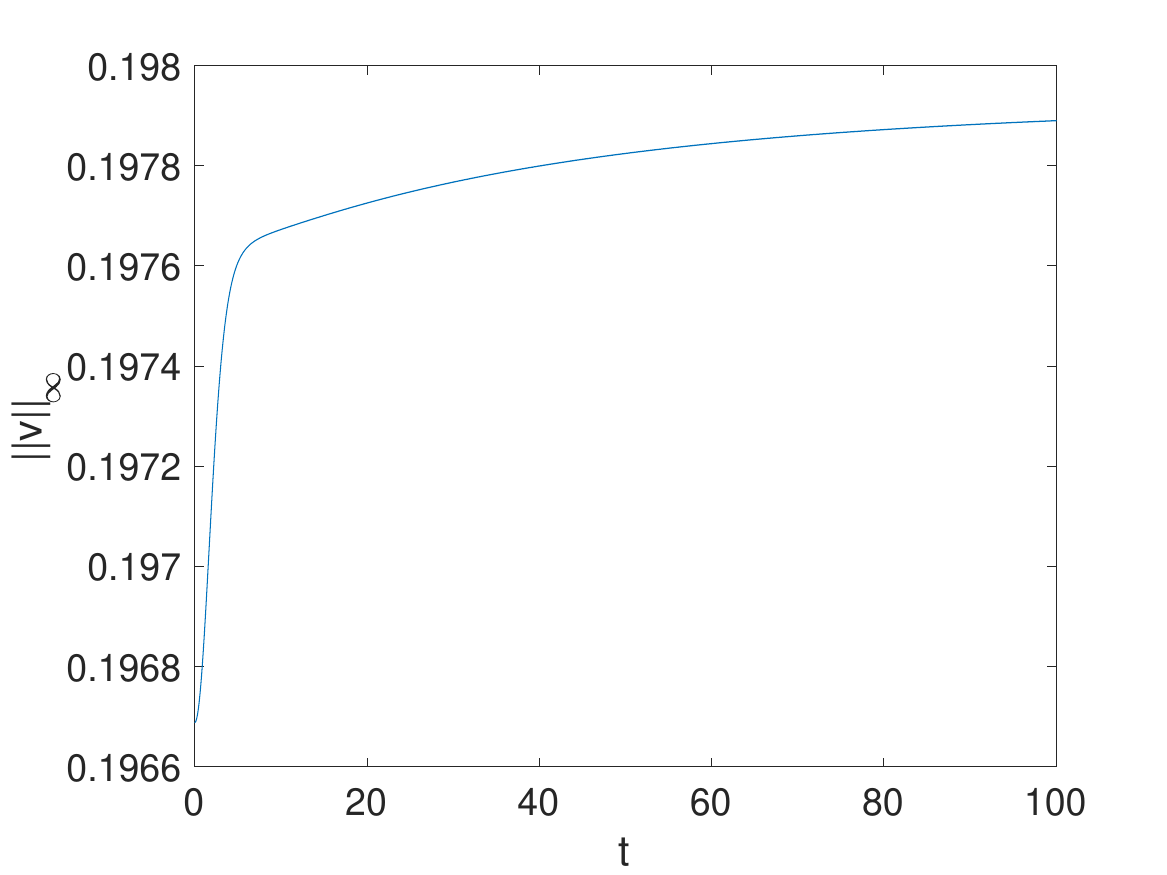}\\
  \includegraphics[width=0.49\textwidth]{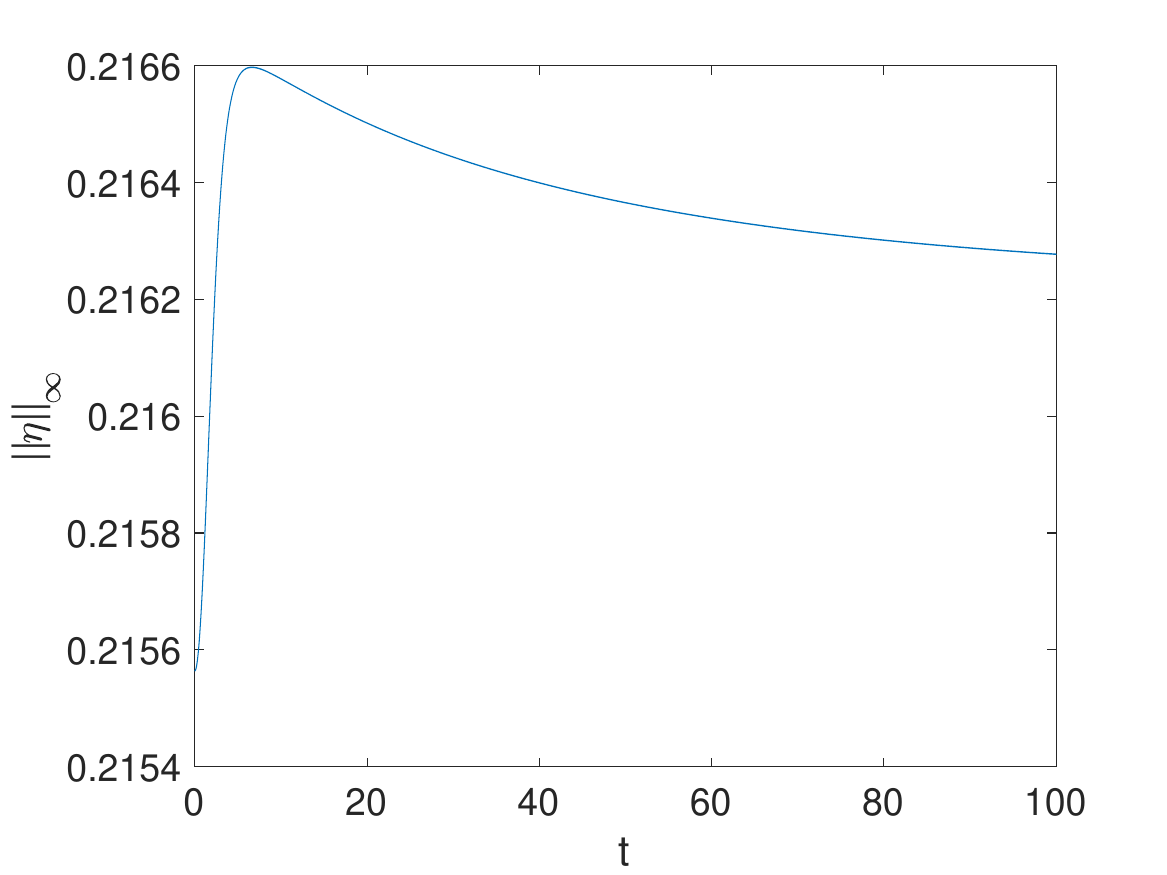}
 \includegraphics[width=0.49\textwidth]{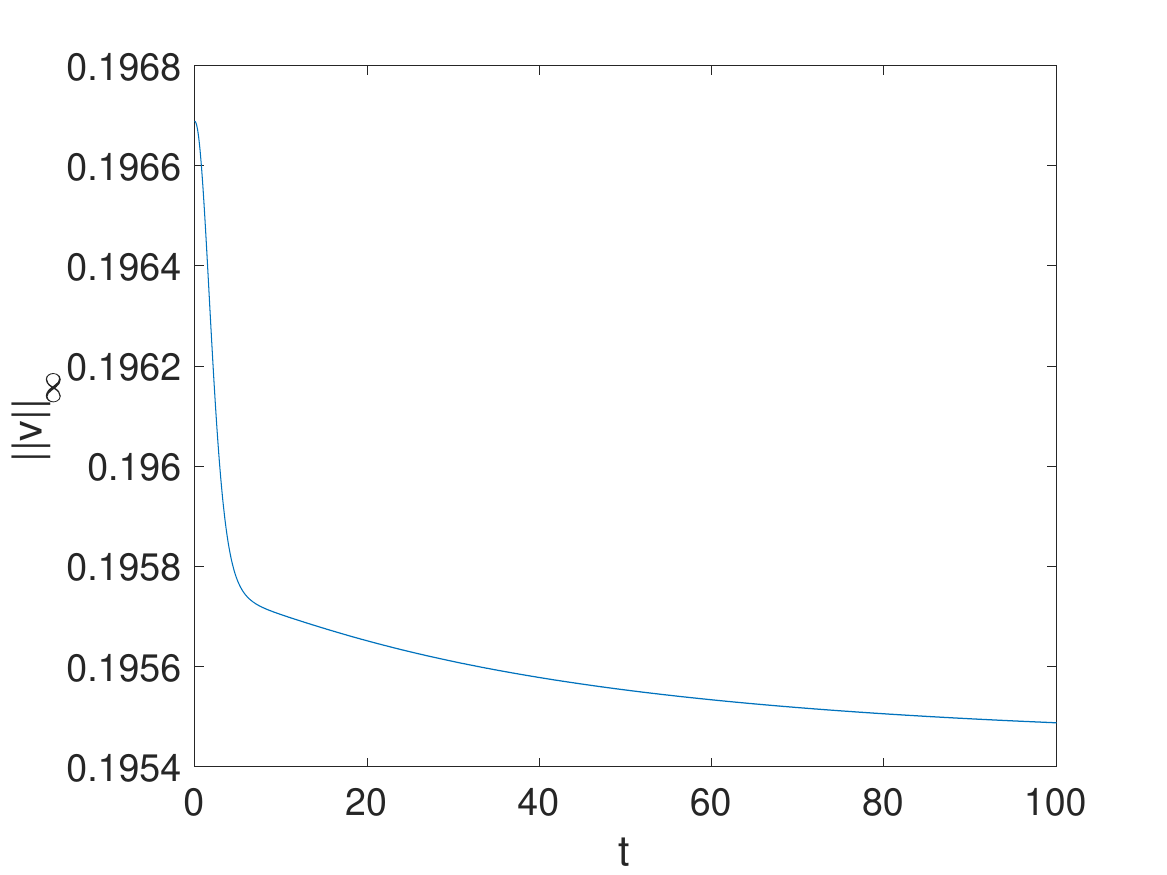}
 \caption{$L^{\infty}$-norms of the solution to the Amick-Schonbeck system (\ref{AS}) for 
 initial data of the form (\ref{lambda}) with $C=1.1$, on the left $\eta$, on the right $v$, in the upper row 
 for the $\lambda=1.01$, in the lower row for the $\lambda=0.99$.}
 \label{figASc11la101}
\end{figure}

The situation is somewhat different if a solitary wave with a 
velocity close to 1 is considered. Here  smaller 
perturbations have to be applied, and it will take much longer 
times to get close to a final state. Therefore we consider the same 
perturbations as above with smaller values for $\lambda$, $\mu$ 
for the solitary wave with $C=1.1$ on a larger 
torus $x\in40[-\pi,\pi]$. We use $N=2^{14}$ DFT modes and 
$N_{t}=2*10^{4}$ time steps for $t\leq 100$. In 
Fig.~\ref{figASc11la101} we show the $L^{\infty}$ norms of the 
solutions to the Amick-Schonbeck system (\ref{AS}) for initial data 
of the form (\ref{lambda}) for $\lambda=1.01$ and $\lambda=0.99$, 
i.e., perturbations of the order of 1\%. It can be seen that the 
$L^{\infty}$ norms only reach final states very slowly, but it 
appears that the latter are again a solitary wave plus radiation. 

\begin{figure}[htb!]
 \includegraphics[width=0.49\textwidth]{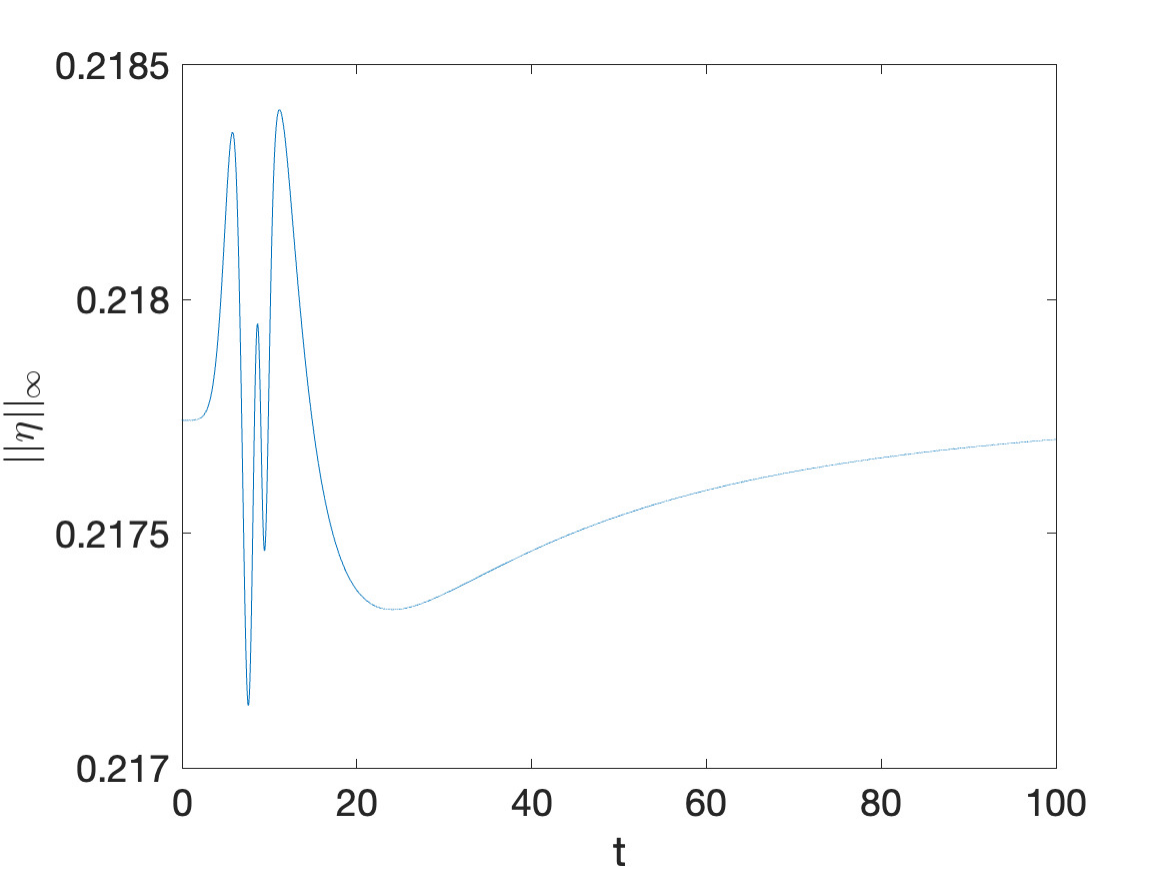}
 \includegraphics[width=0.49\textwidth]{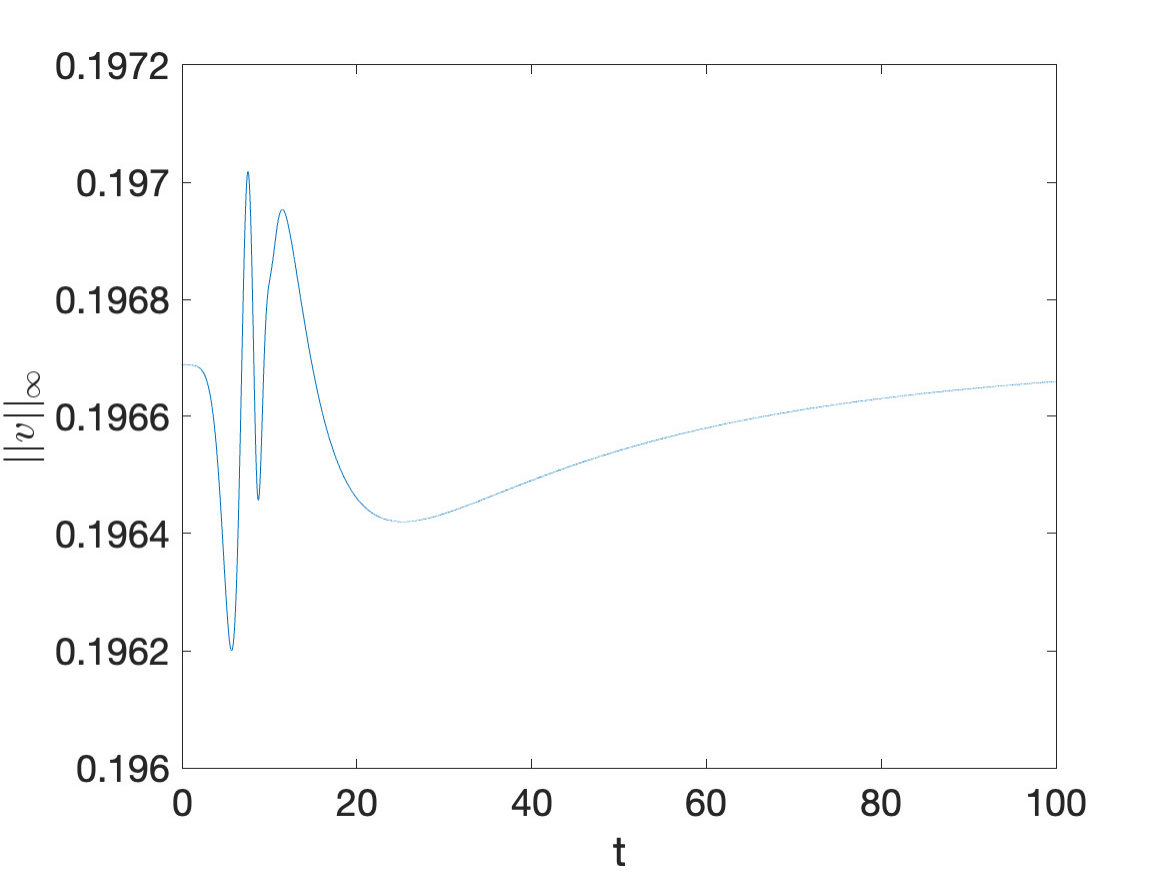}\\
  \includegraphics[width=0.49\textwidth]{Amicksolc11_0002gaussetainf.eps}
 \includegraphics[width=0.49\textwidth]{Amicksolc11_0002gaussvinf.eps}
 \caption{$L^{\infty}$ norms of the solution to the Amick-Schonbeck system (\ref{AS}) for the 
 initial data $\eta(x,0) = Q_{C}(x)\pm 0.002 \exp(-x^{2})$, $v(x,0) = 
 V_{C}(x)$ with $C=1.1$, on the left $\eta$, on the right $v$, in the upper row 
 for the $+$ sign, in the lower row for the $-$ sign.}
 \label{figASc11mu0002}
\end{figure}

We also consider initial data of the form (\ref{gausspert}) for 
solitary waves with velocity $C=1.1$, but again a perturbation of the 
order of 1\%, $\mu=0.002$. The $L^{\infty}$ norms of the solutions 
can be seen in Fig.~\ref{figASc11mu0002}. 

\begin{figure}[htb!]
 \includegraphics[width=0.49\textwidth]{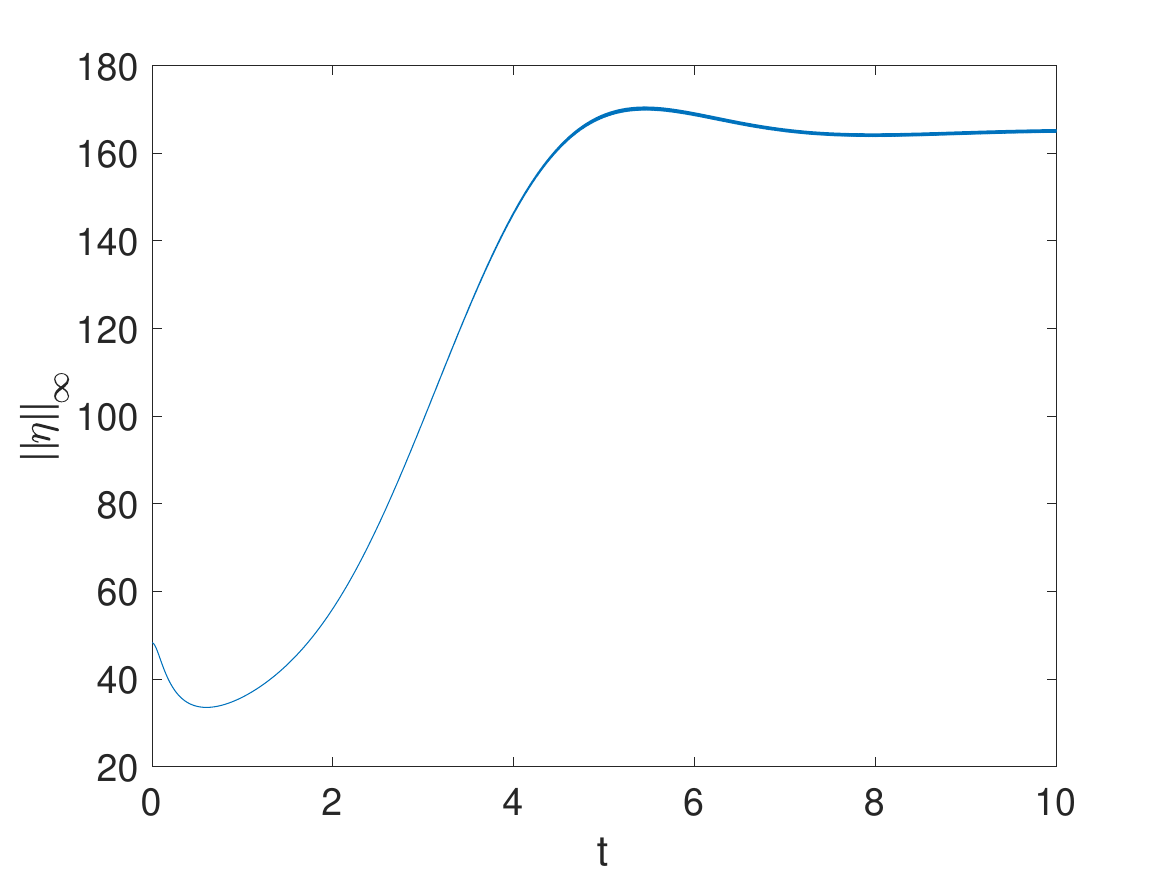}
 \includegraphics[width=0.49\textwidth]{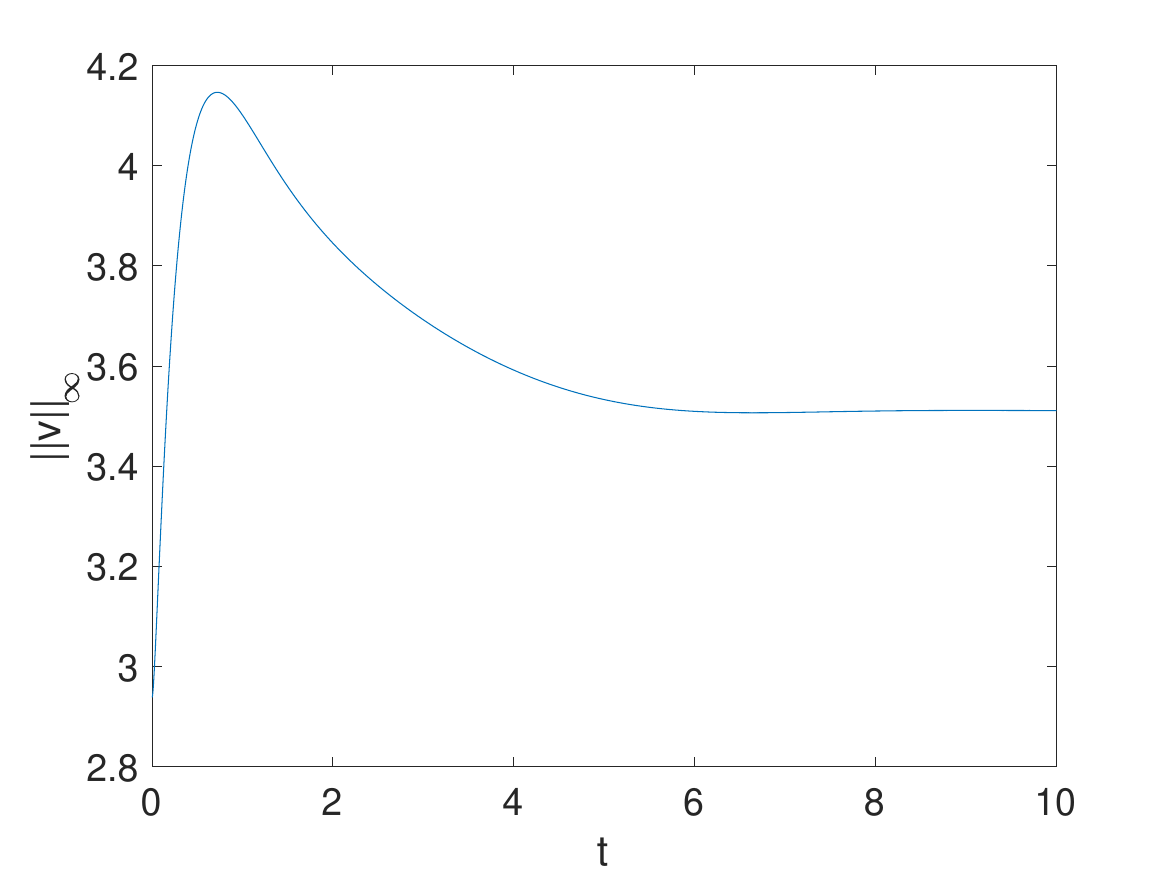}\\
  \includegraphics[width=0.49\textwidth]{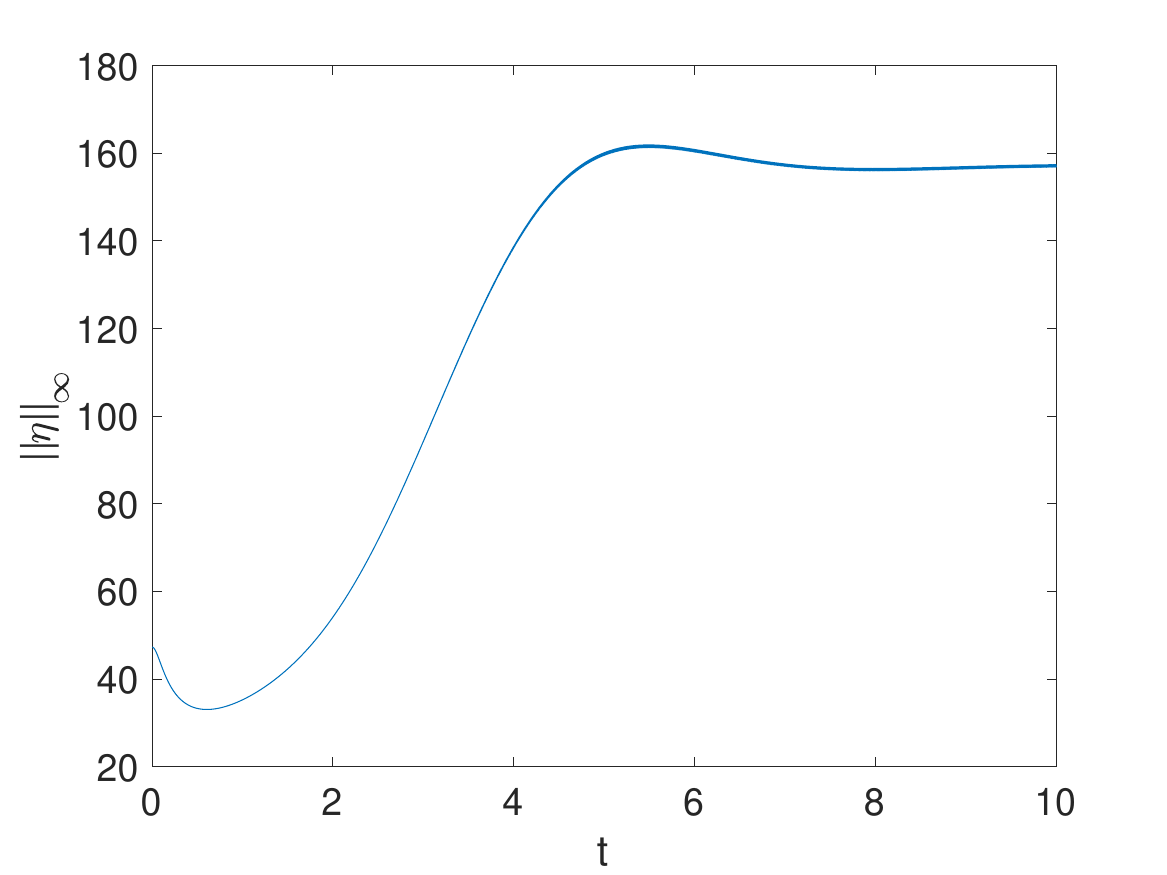}
 \includegraphics[width=0.49\textwidth]{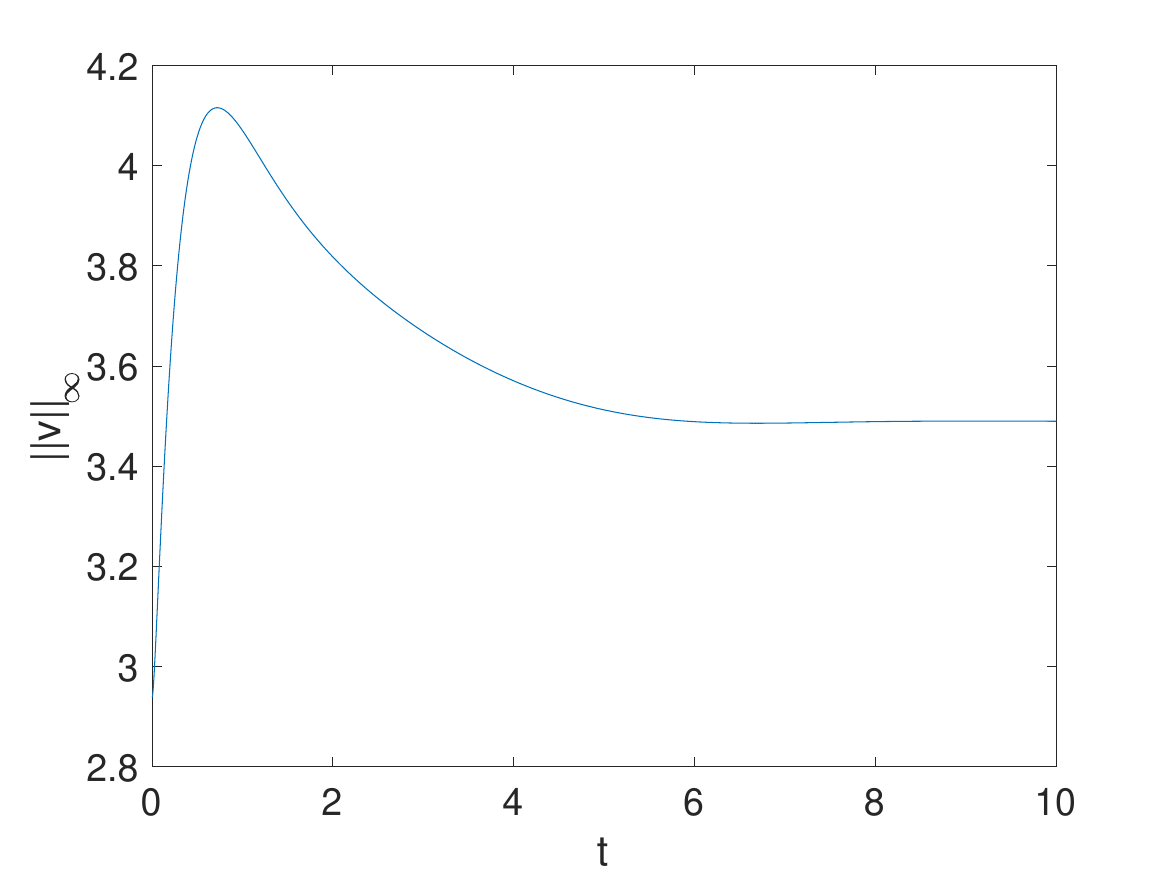}
 \caption{$L^{\infty}$-norms of the solution to the Amick-Schonbeck system (\ref{AS}) for 
 initial data of the form (\ref{lambda}) with $C=3$, on the left $\eta$, on the right $v$, in the upper row 
 for the $\lambda=1.01$, in the lower row for the $\lambda=0.99$.}
 \label{figASc3la101}
\end{figure}

The situation for values of $C$ larger than 2 is similar to the one 
for $c\sim 1$. We use $N=2^{15}$ DFT modes for $x\in 20[-\pi,\pi]$ 
and $N_{t}=10^{5}$ time steps for $t\leq 10$. A small perturbation of the form (\ref{lambda}) as in 
Fig.~\ref{figASc11la101} leads for $C=3$ to a small perturbation in 
the $L^{\infty}$ norm of $v$ indicating the final state is a solitary 
wave with a slightly higher velocity. But as can be seen in 
Fig.~\ref{figASc3la101} this implies a strong change of the 
$L^{\infty}$ norm of $\eta$ since the latter varies rapidly with $C$ 
for larger values of $C$, see Fig.~\ref{figASsolc23}.
latter is strongly growing

\section{Localised initial data}
In this section we discuss the time evolution of localised initial 
data for the Amick-Schonbeck system (\ref{AS}). We consider initial 
data of the form
\begin{equation}
	\eta(x,0) = A\exp(-x^{2}),\quad v(x,0) = 0,
	\label{gauss}
\end{equation}
where $A$ is a positive constant. The results give strong evidence 
to the soliton resolution conjecture for this system. 

\begin{figure}[htb!]
 \includegraphics[width=0.49\textwidth]{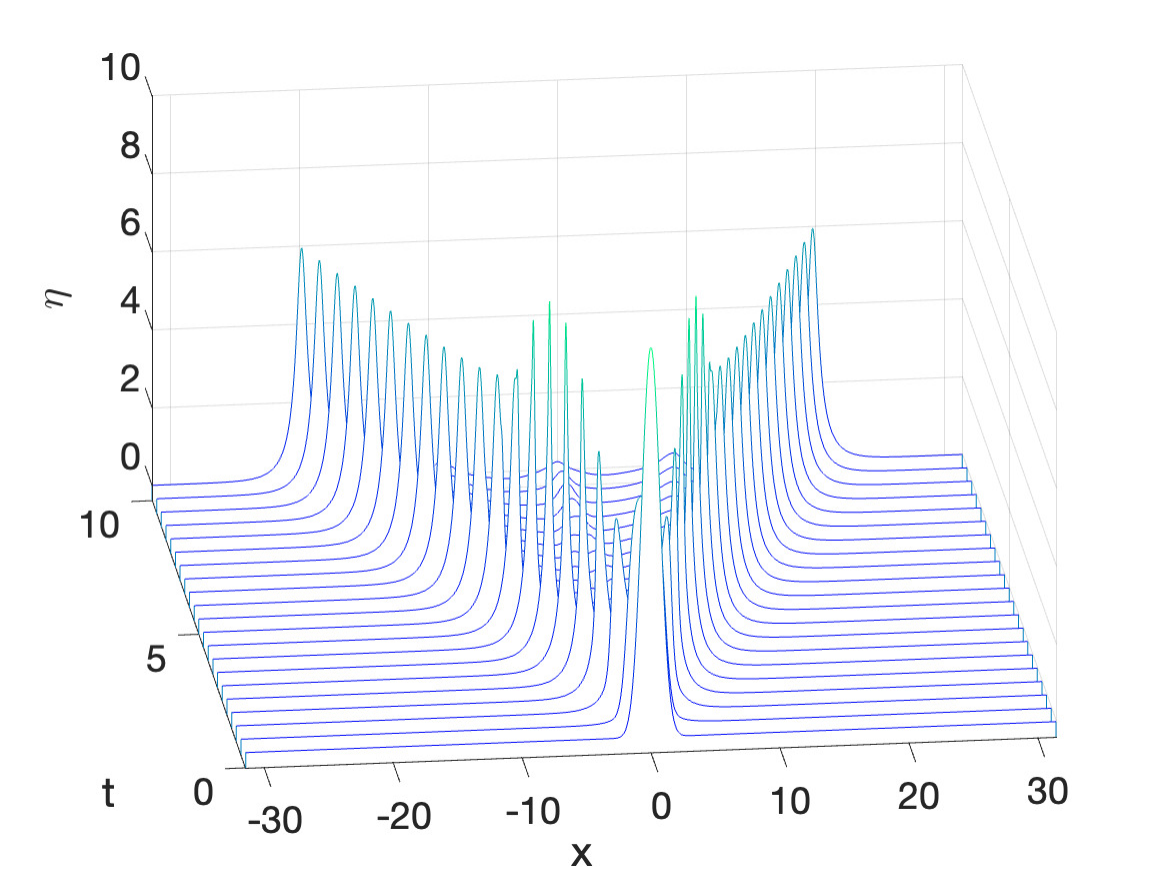}
 \includegraphics[width=0.49\textwidth]{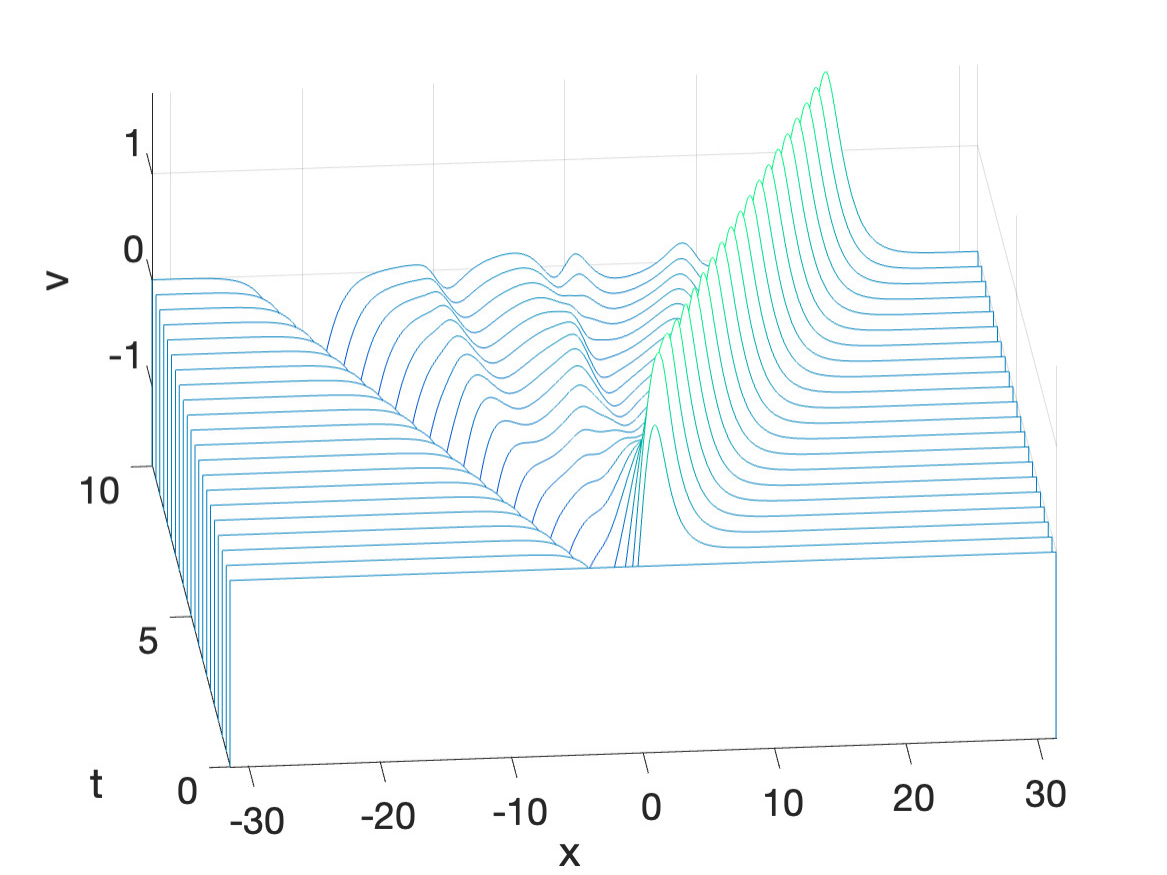}
 \caption{Solution to the Amick-Schonbeck system (\ref{AS}) for  
 initial data of the form (\ref{gauss}) with $A=10$, on the left $\eta$, on the right $v$.}
 \label{figAS10gauss}
\end{figure}

We first consider the case of a large constant, $A=10$. We use 
$N=2^{14}$ DFT modes for $x\in10[-\pi,\pi]$ and $N_{t}=10^{4}$ time steps for $t\leq 10$. 
For symmetry 
reasons the initial hump in $\eta$ splits into two humps traveling 
towards positive and negative values of $x$ respectively. This can be seen in 
Fig.~\ref{figAS10gauss} on the left. The function $v$ is positive 
for positive values of $x$ and negative for negative $x$. 

\begin{figure}[htb!]
 \includegraphics[width=0.49\textwidth]{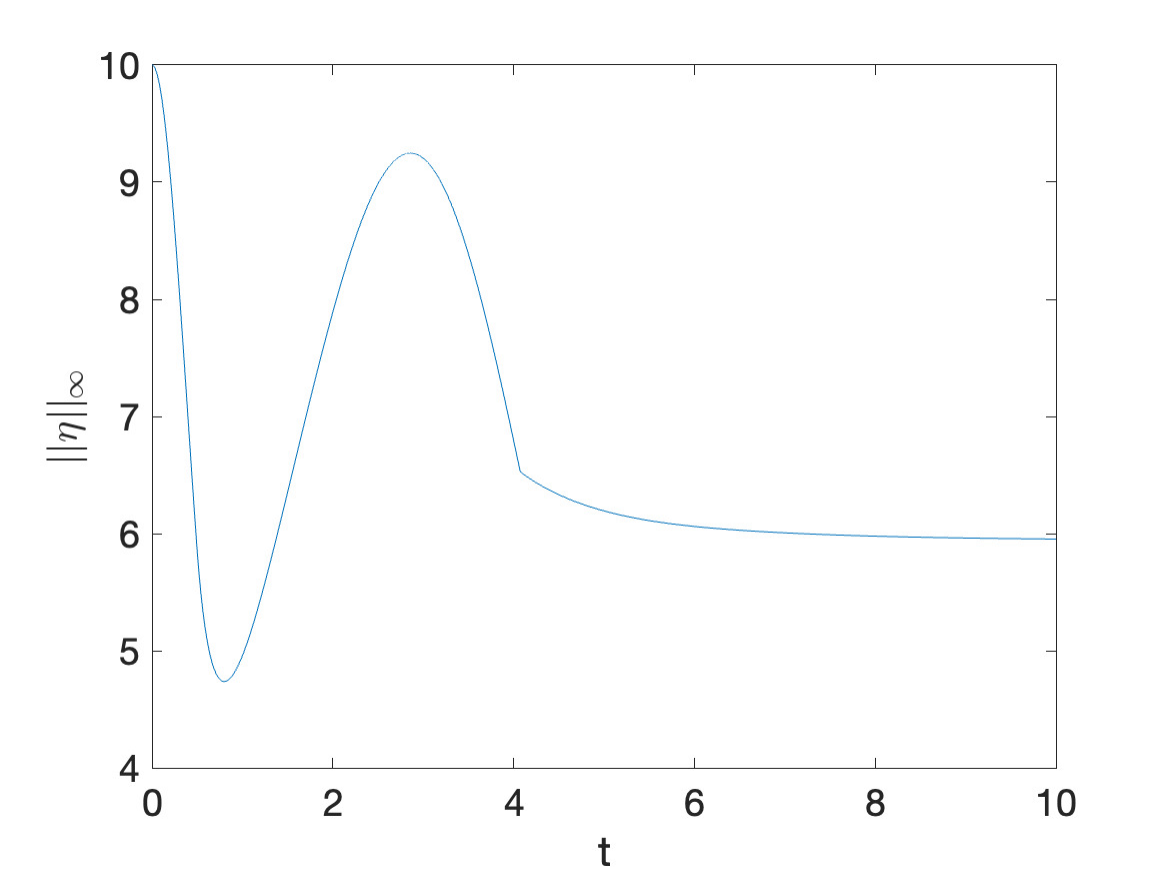}
 \includegraphics[width=0.49\textwidth]{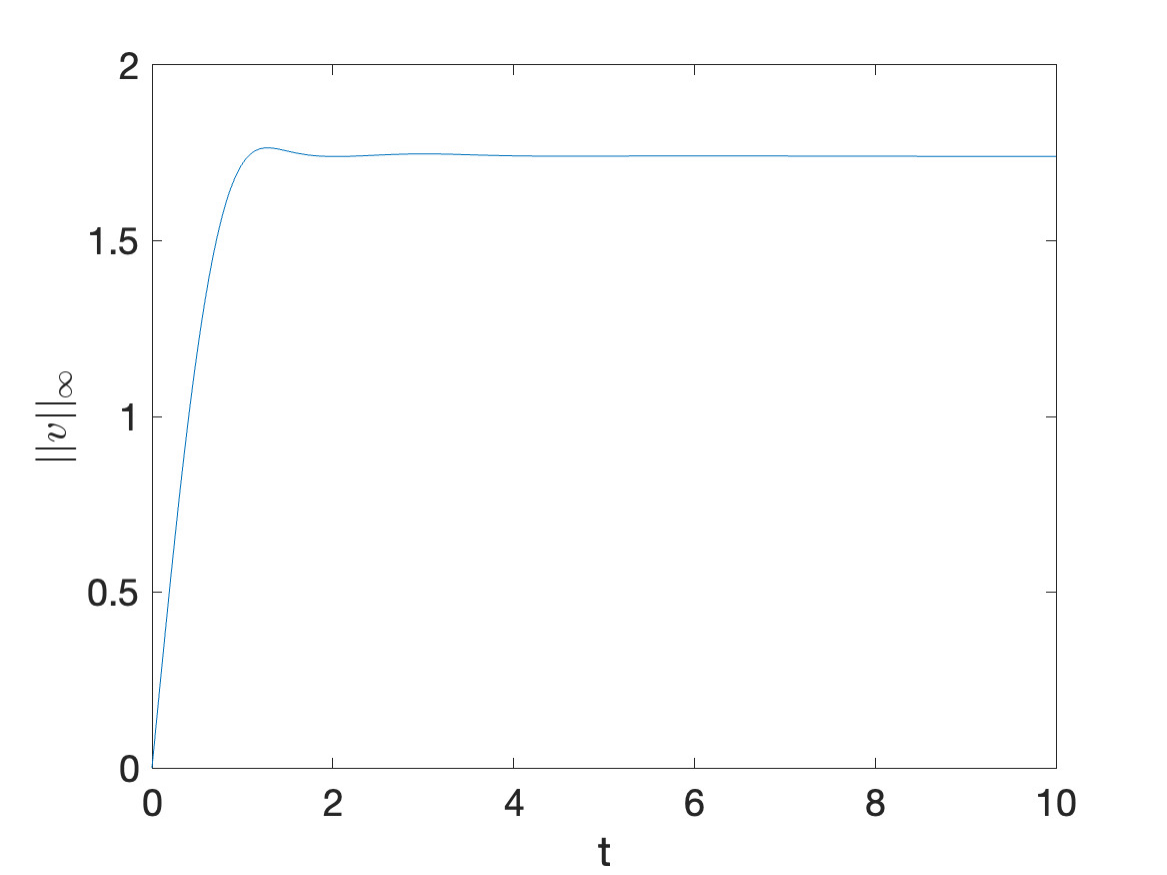}
 \caption{$L^{\infty}$ norms of the solution to the Amick-Schonbeck system (\ref{AS}) for  
 initial data of the form (\ref{gauss}) with $A=10$, on the left $\eta$, on the right $v$.}
 \label{figAS10gaussinf}
\end{figure}

The final state of the solution appears to be two solitary waves, one 
with positive velocity, and one with negative velocity. This 
interpretation is confirmed by the $L^{\infty}$ norms of the solution 
in Fig.~\ref{figAS10gaussinf}. As expected from the stability of the 
solitary waves studied in the previous section, solitary waves seem 
to appear in the long time behavior of the solutions in accordance 
with the soliton resolution conjecture.

\begin{figure}[htb!]
 \includegraphics[width=0.49\textwidth]{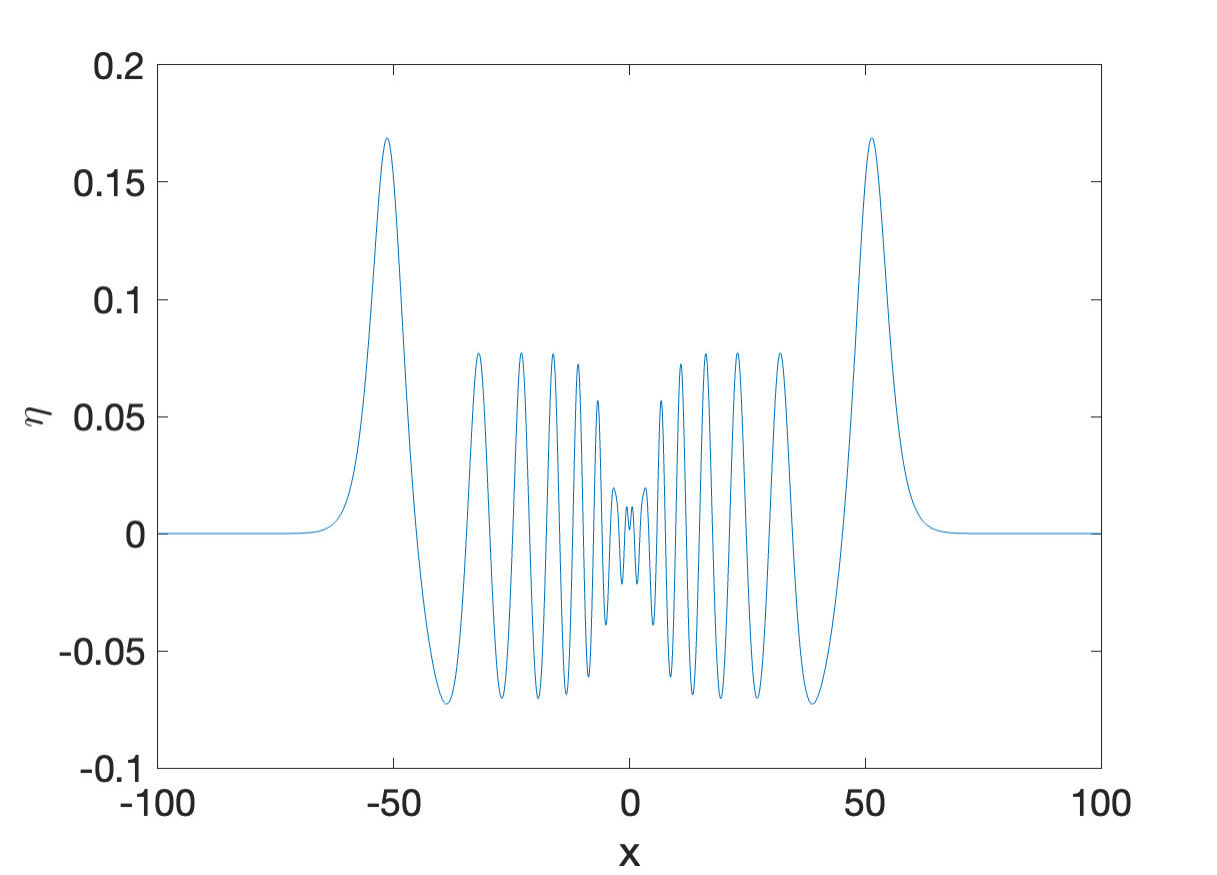}
 \includegraphics[width=0.49\textwidth]{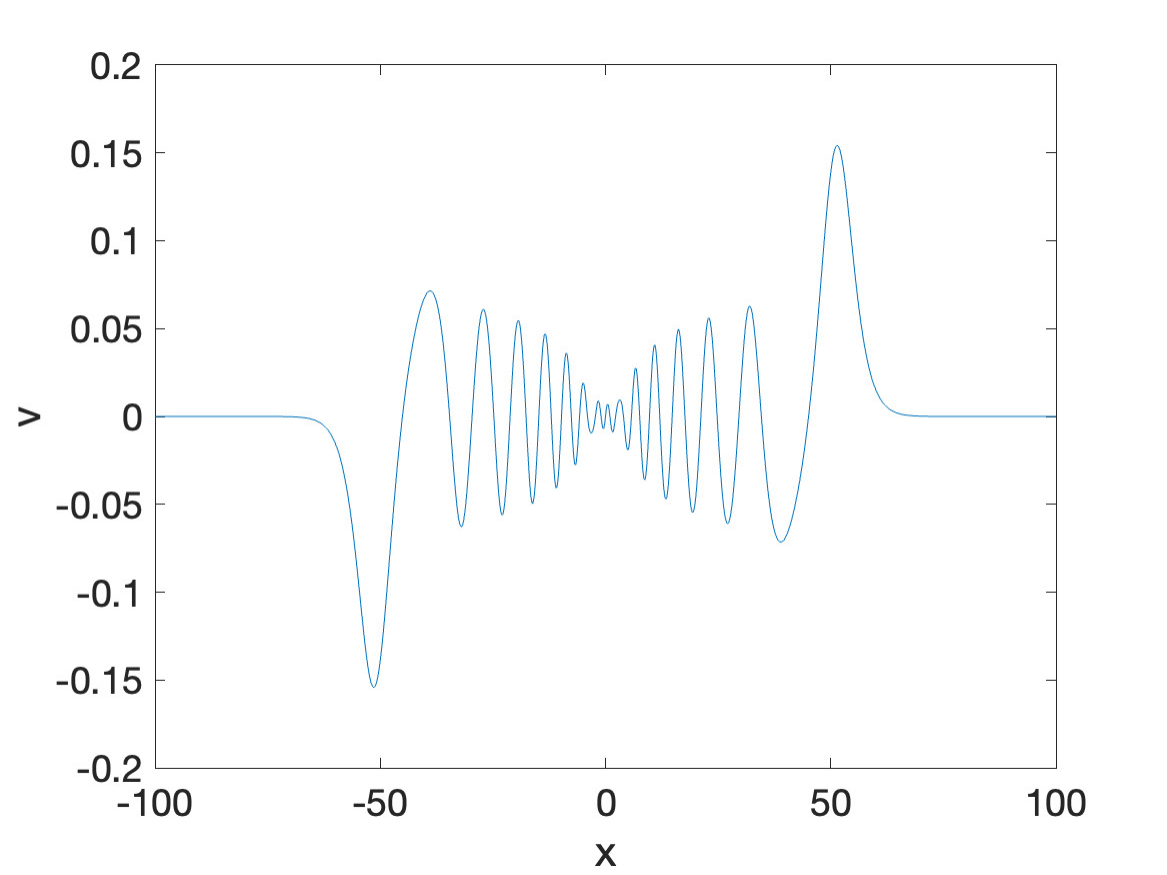}
 \caption{Solution to the Amick-Schonbeck system (\ref{AS}) for  
 initial data of the form (\ref{gauss}) with $A=1$ for $t=50$, on the left $\eta$, on the right $v$.}
 \label{figASgauss}
\end{figure}

For smaller values of $A$, say $A=1$, the same phenomenon as in the 
previous section is observed, that much larger times are needed to 
observe the formation of solitary waves with $C\sim1$. We use the 
same number of DFT modes, but on a larger torus, $x\in40[-\pi,\pi]$ 
and $N_{t}=2*10^{4}$ time steps for $t\leq 50$. The solution at the 
final time is shown in Fig.~\ref{figASgauss}. There are again two 
humps, but now with strong radiation in comparison to the amplitude 
of the presumed solitary waves. 

\begin{figure}[htb!]
 \includegraphics[width=0.49\textwidth]{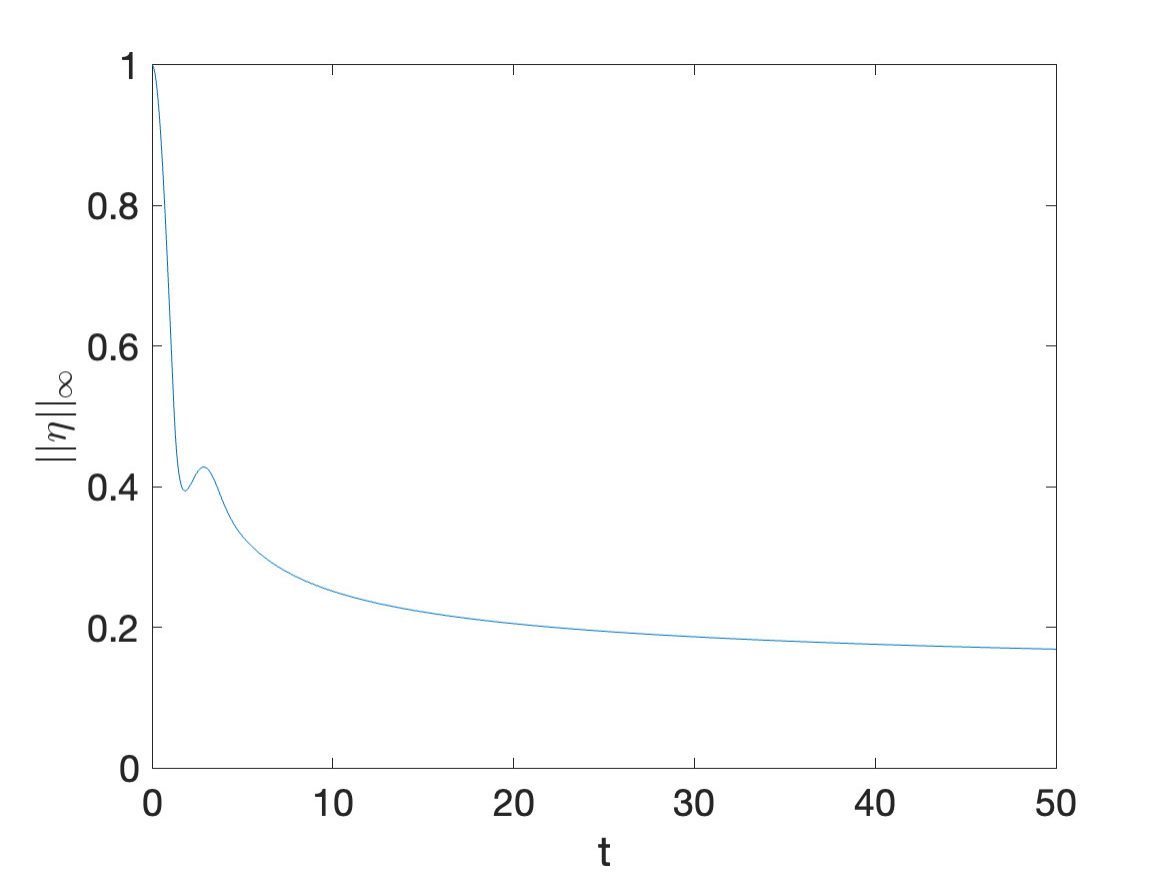}
 \includegraphics[width=0.49\textwidth]{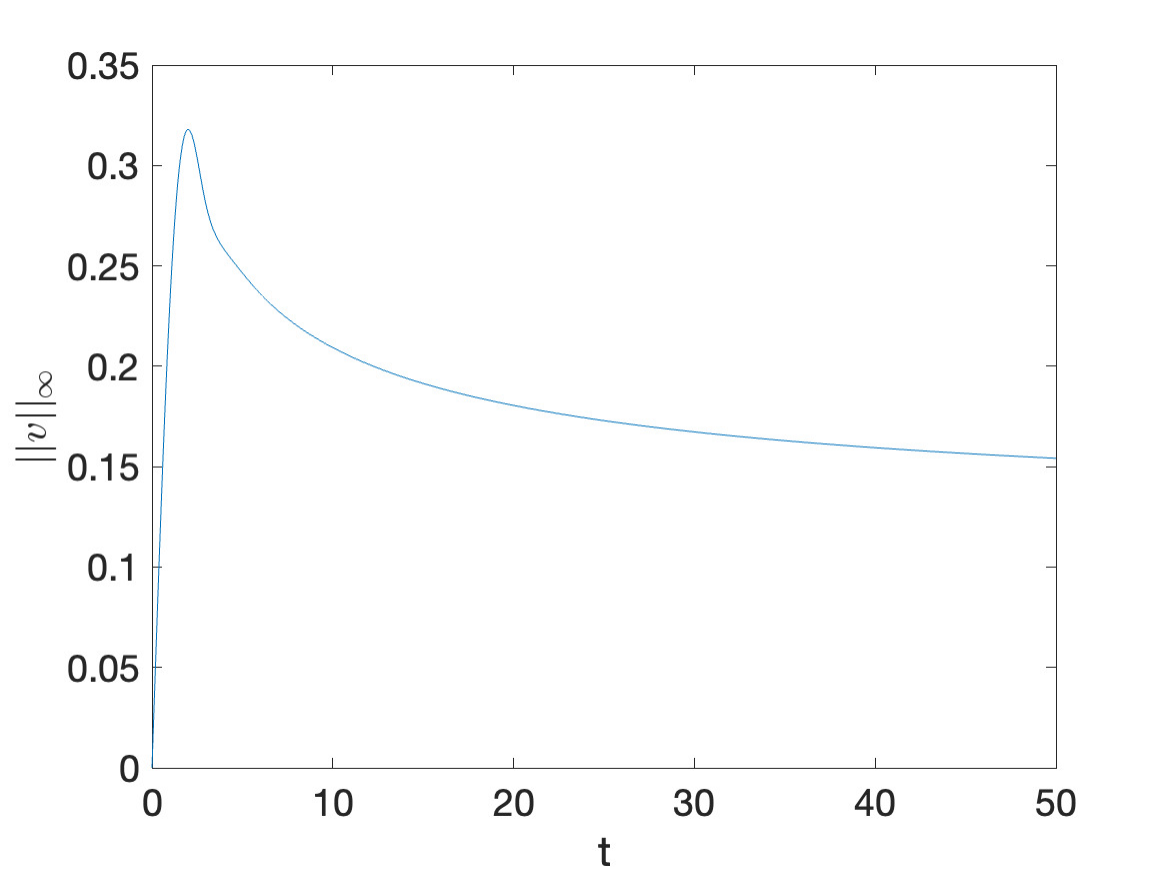}
 \caption{$L^{\infty}$ norms of the solution to the Amick-Schonbeck system (\ref{AS}) for  
 initial data of the form (\ref{gauss}) with $A=1$, on the left $\eta$, on the right $v$.}
 \label{figASgaussinf}
\end{figure}

The $L^{\infty}$ norms of these solutions can be seen in 
Fig.~\ref{figASgaussinf}. They appear to saturate, but the final 
state, which seems to be solitary waves plus radiation,  is obviously 
not yet reached. 

\section{Non-cavitation initial data}
In this section we study the Amick-Schonbeck system (\ref{AS}) for 
initial data that violate the non-cavitation  condition 
$1+\eta>0$ or that are close to violating it. The results are similar 
to what was found in \cite{DK} for the Serre-Green-Naghdi equation 
and give strong evidence for the second conjecture in the introduction. 

As we will show below, it is numerically challenging to study the 
solution of the system (\ref{AS}) for initial data that are close to 
or actually violating the non-cavitation condition $1+\eta>0$ since 
the solution appears to become singular or almost singular in finite 
time. We 
consider initial data of the form 
\begin{equation}
	\eta(x,0) = -\kappa \exp(-x^{2}),\quad v(x,0) = 0,
	\label{NC}
\end{equation}
where $\kappa$ is a positive constant. First we consider the case 
$\kappa=1$, i.e., an example violating the non-cavitation condition 
in one point, $x=0$ for $t=0$. We use $N=2^{18}$ DFT modes for 
$x\in3[-\pi,\pi]$ and $N_{t}=10^{5}$ time steps for $t\leq 5$, close 
to the accessible accuracy with double precision arithmetics. It can 
be seen that the initial depression in $\eta$ disappears and that it 
eventually
develops a strongly peaked maximum at a high elevation. The function $v$ on the other 
hand is as before odd in $x$ with a maximum to the left and a
minimum to the right of the origin. For larger times, a new small  
oscillation of this function appears near the origin.
\begin{figure}[htb!]
 \includegraphics[width=0.49\textwidth]{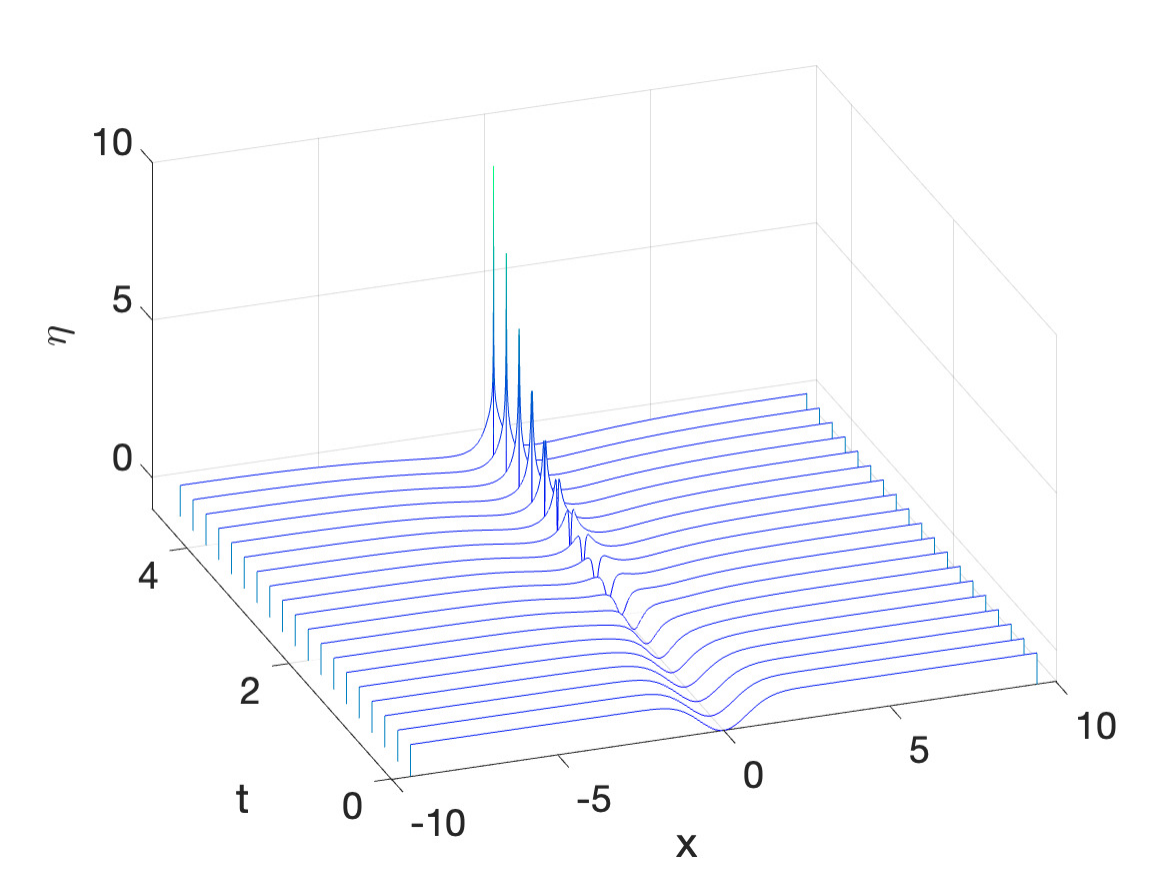}
 \includegraphics[width=0.49\textwidth]{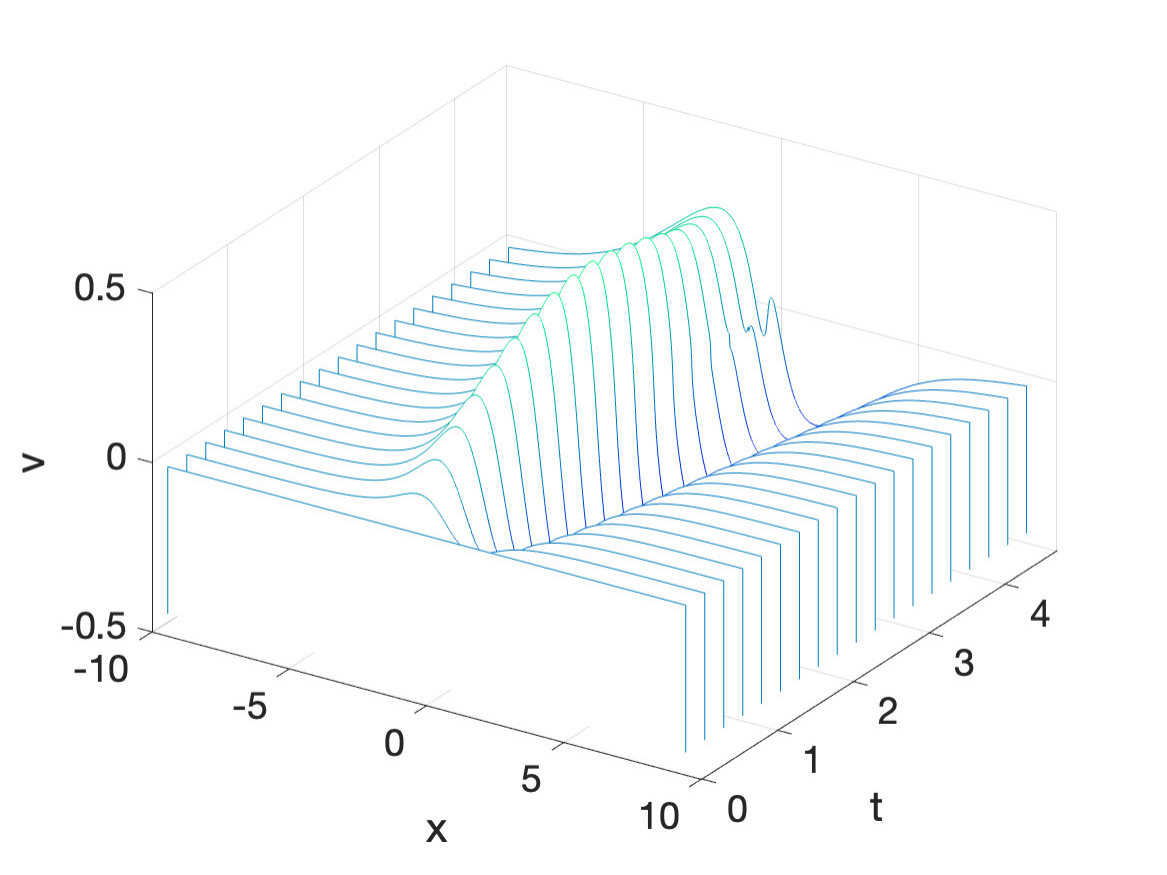}
 \caption{Solution to the Amick-Schonbeck system (\ref{AS}) for  
 initial data of the form (\ref{NC}) with $\kappa=1$, on the left $\eta$, on the right $v$.}
 \label{figASmgauss}
\end{figure}

\begin{figure}[htb!]
 \includegraphics[width=0.49\textwidth]{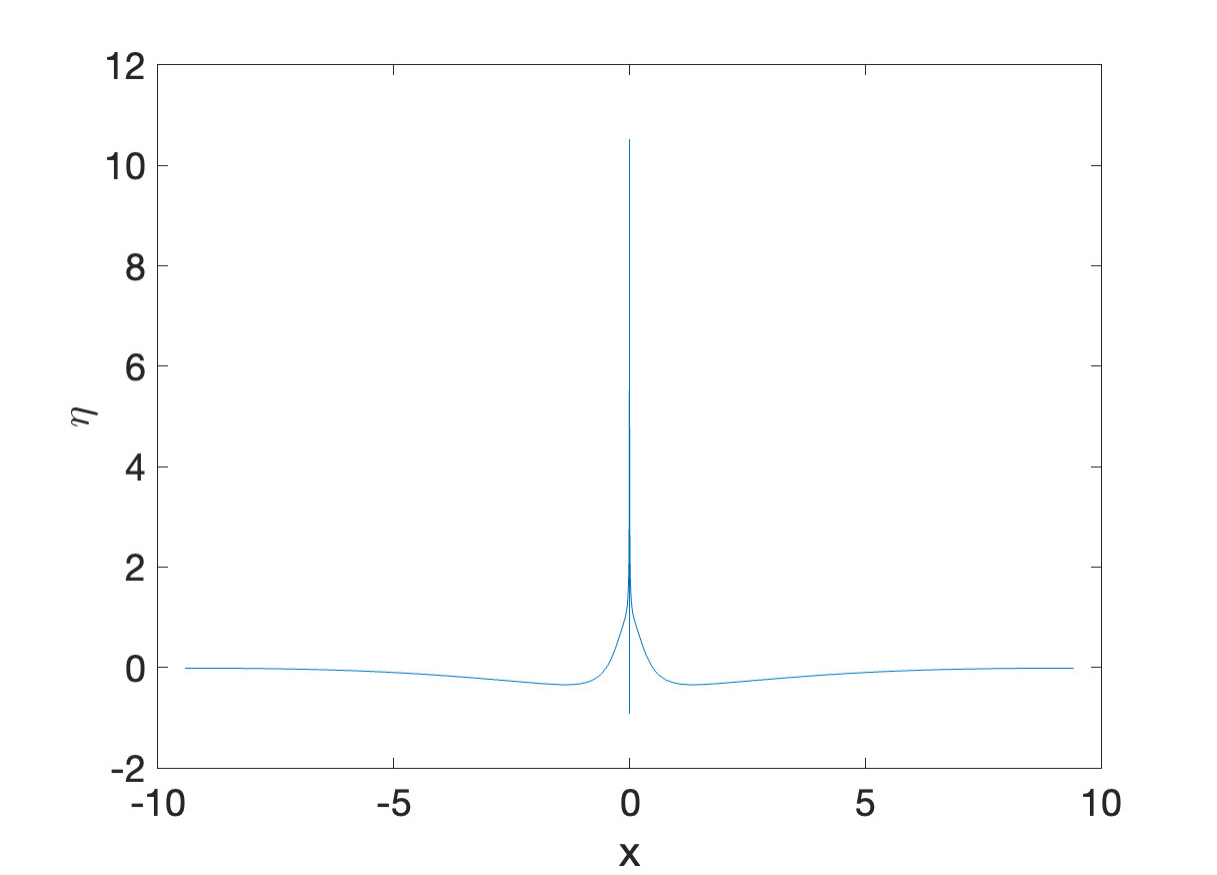}
 \includegraphics[width=0.49\textwidth]{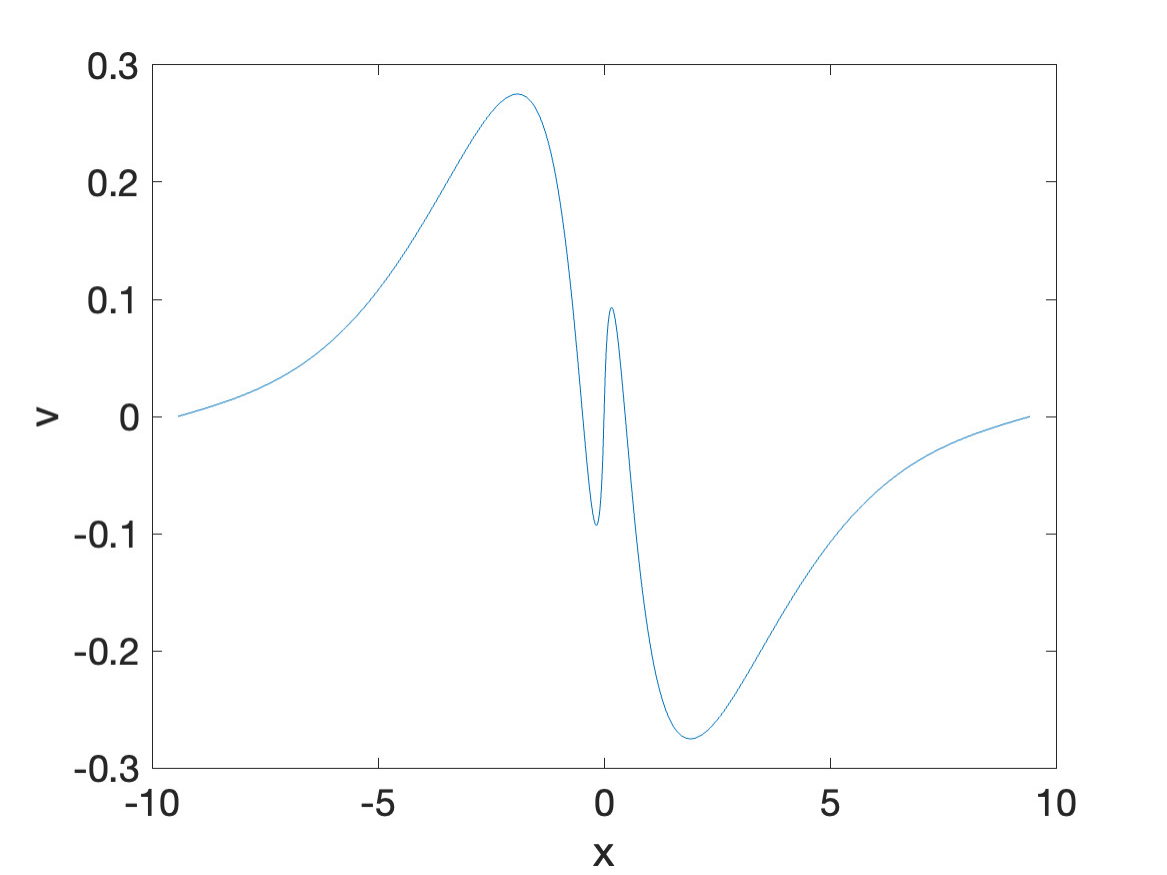}
 \caption{Solution to the Amick-Schonbeck system (\ref{AS}) for  
 initial data of the form (\ref{NC}) with $\kappa=1$ for $t\sim 4.68$, on the left $\eta$, on the right $v$.}
 \label{figASmgausst486}
\end{figure}

We show the solution at the last recorded time $t\sim4.68$ in 
Fig.~\ref{figASmgausst486}. It can be seen that the function $\eta$ at 
this time is strongly peaked near the origin, and that $v$ has a 
strong gradient there. 

\begin{figure}[htb!]
 \includegraphics[width=0.49\textwidth]{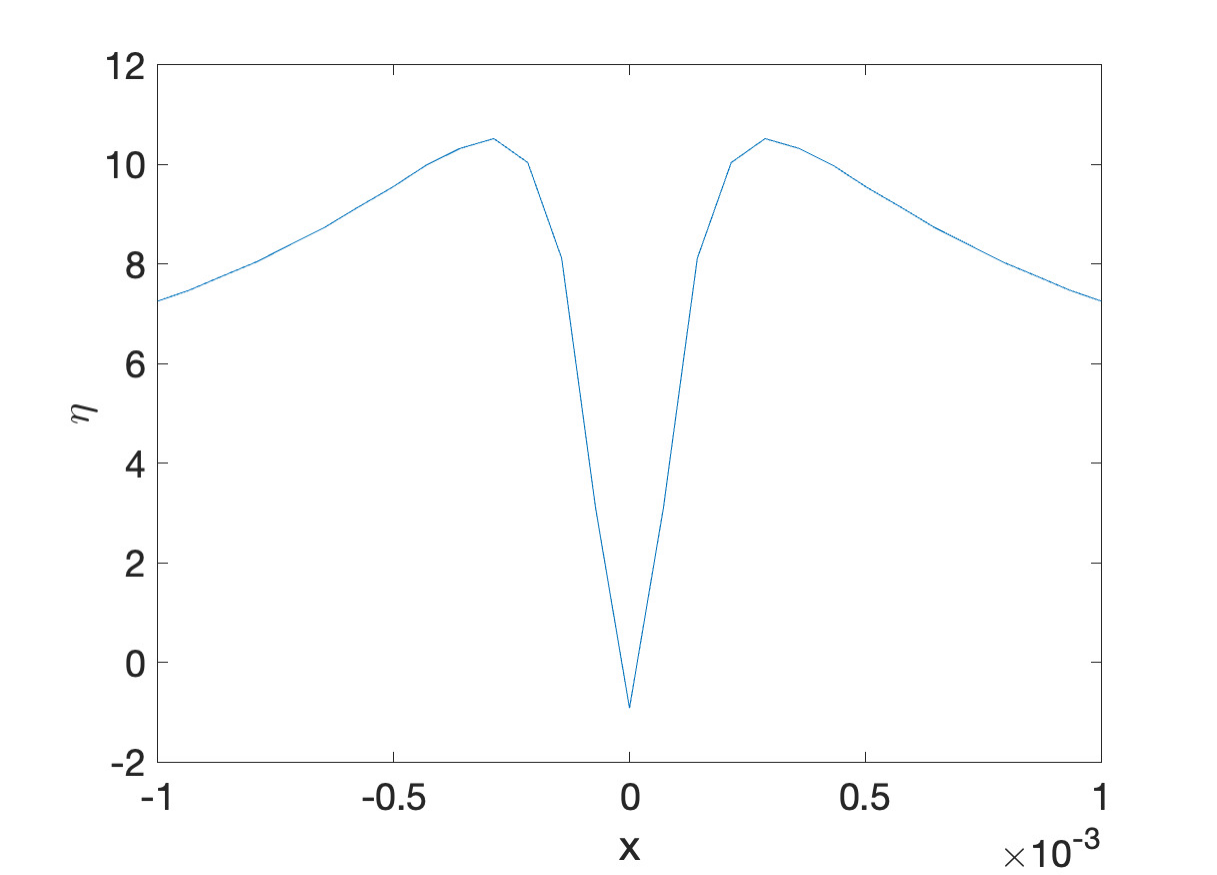}
 \includegraphics[width=0.49\textwidth]{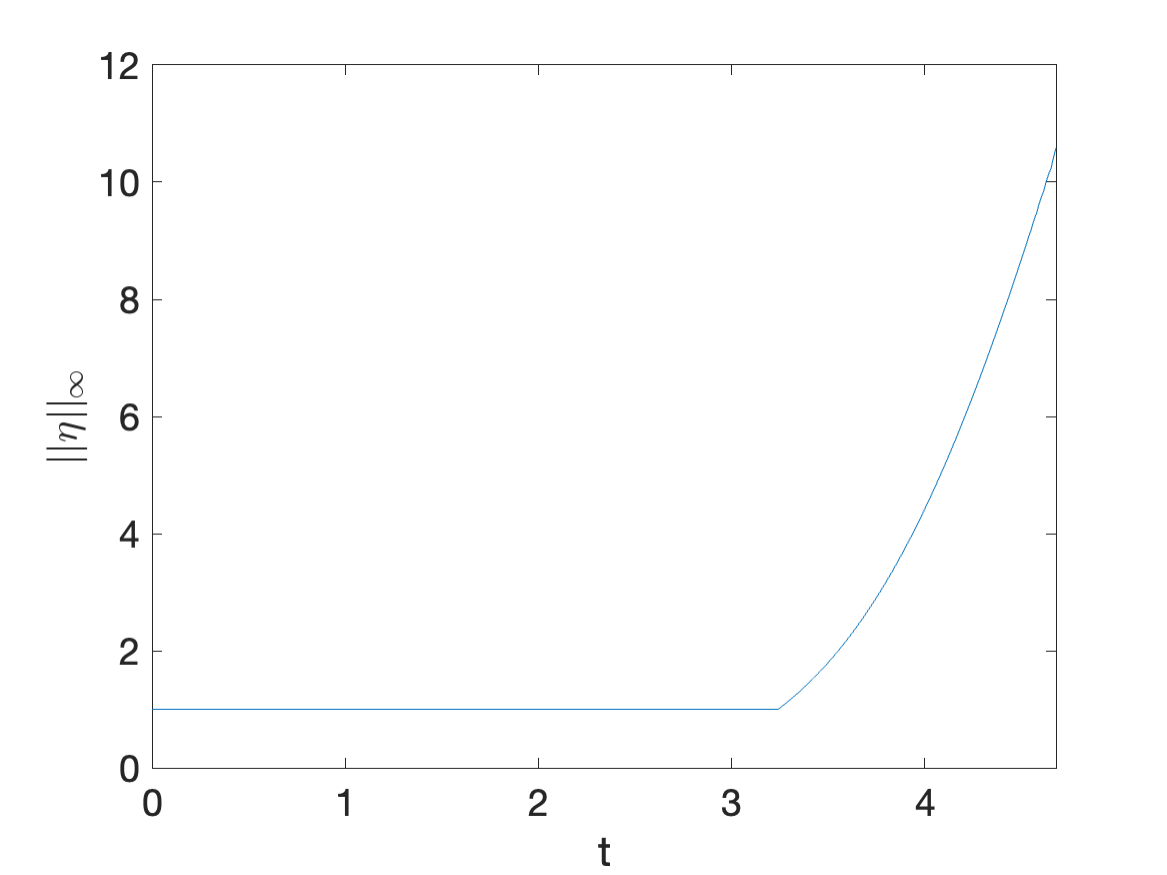}
 \caption{Close-up of the function $\eta$ on the left of 
 Fig.~\ref{figASmgausst486} on the left and the $L^{\infty}$ norm of 
 $\eta$ in dependence of time on the right.}
 \label{figASmgausst486cu}
\end{figure}

A close-up of the the function $\eta$ near the origin, see 
Fig.~\ref{figASmgausst486cu},  reveals that 
there is a strong cusp-like gap near the origin. It is this cusp-like 
structure that is numerically difficult to resolve and potentially 
singular, not the peak. Note that the function $\eta$ appears to be 
negative at the minimum in the close-up, but just at this point. 
The $L^{\infty}$ norm of the function $\eta$ on the right of the same figure 
shows that the peak near the origin is growing which could correspond 
to an $L^{\infty}$ blow-up. 

\begin{figure}[htb!]
 \includegraphics[width=0.49\textwidth]{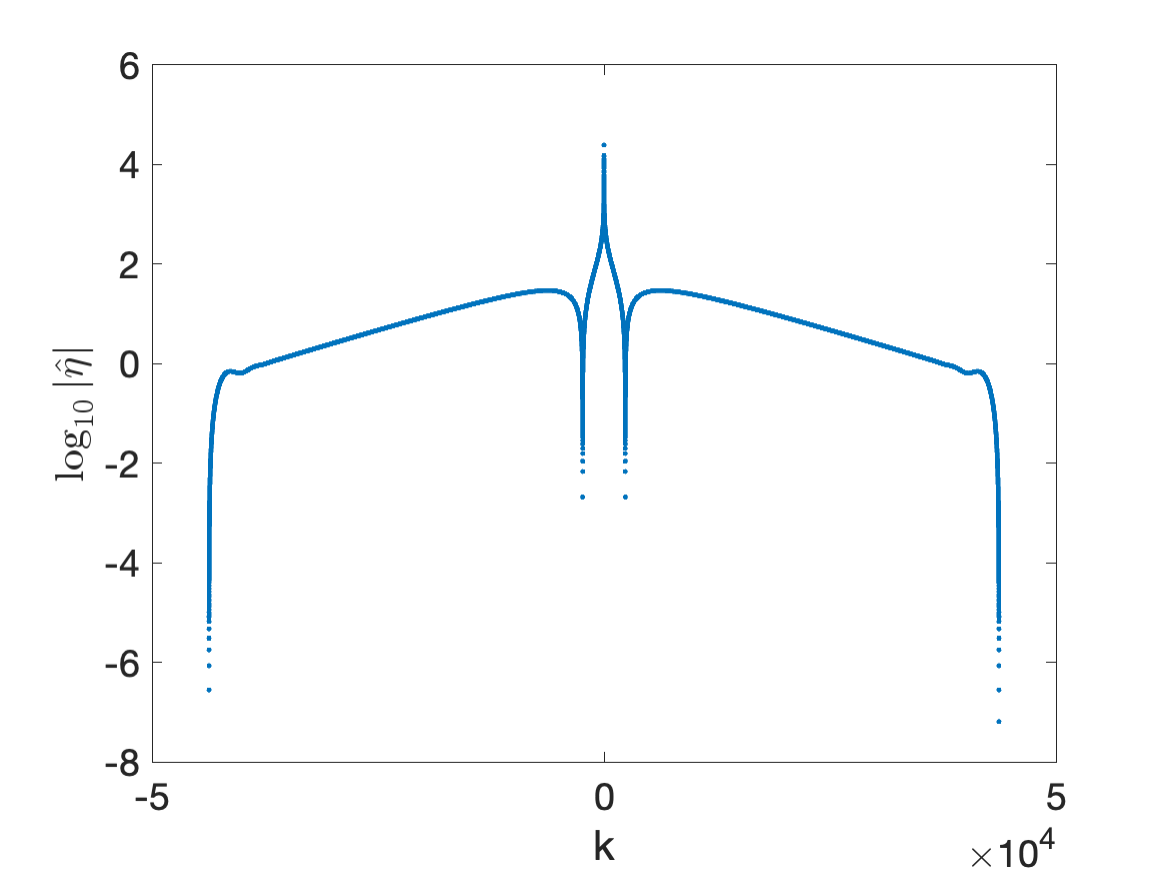}
 \includegraphics[width=0.49\textwidth]{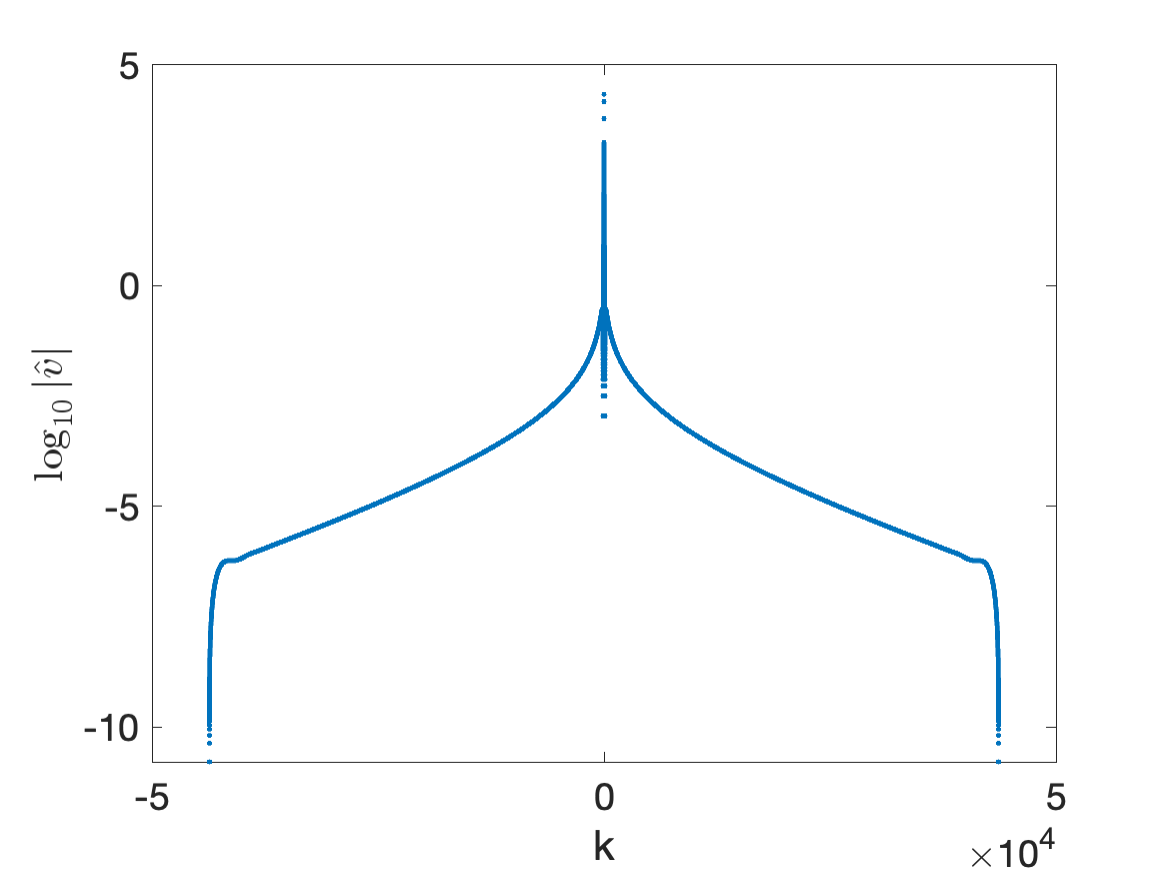}
 \caption{Fourier coefficients of  the solution to the Amick-Schonbeck system (\ref{AS}) for  
 initial data of the form (\ref{NC}) with $\kappa=1$ for $t\sim 4.68$, on the left $\eta$, on the right $v$.}
 \label{figASmgausst486fourier}
\end{figure}

An interesting approach to numerically study singularities in a 
Fourier spectral method was introduced by Sulem, Sulem and Frisch in 
\cite{SSF}.  It is well 
known that an essential singularity in the complex plane of the form $u\sim 
(z-z_{j})^{\mu_{j}}$, $\mu_{j}\notin \mathbb{Z}$, 
with $z_{j}=\alpha_{j}-i\delta_{j}$ in the lower half 
plane ($\delta_{j}\geq 0$) implies    for $k\to\infty$  
the following asymptotic behavior of the Fourier 
transform (see e.g.~\cite{asymbook}, here denoted in the same way as the DFT),
\begin{equation}
    \hat{u}(k)\sim 
    \sqrt{2\pi}\mu_{j}^{\mu_{j}+\frac{1}{2}}e^{-\mu_{j}}\frac{(-i)^{\mu_{j}+1}}{k^{\mu_{j}+1}} e^{-ik\alpha_{j}-k\delta_{j}}
    \label{fourierasym}.
\end{equation}
For a single such singularity with positive $\delta_{j}$, the modulus of the Fourier 
transform decreases exponentially for large $k$ until
$\delta_{j}=0$, when 
this modulus  has an algebraic dependence 
on $k$ for large $|k|$. This behavior can also be seen in the discrete version of the 
Fourier transform, the DFT. For the example in 
Fig.~\ref{figASmgausst486}, the DFT coefficients are shown in 
Fig.~\ref{figASmgausst486fourier}. The algebraic decay of the DFT 
coefficients for large $|k|$ indicates that a singularity of the 
approximated function in the complex plane hits the real axis in this 
case. The DFT coefficients are fitted as discussed in \cite{KR}, $\ln 
|\hat{\eta}|$ with linear regression to $\mu \ln |k| 
-\delta|k|+const$. Concretely two values of $k$ are 
chosen, $k_{min}\gg1$ and $k_{max}<N/2/L$ such that the following 
results do not change considerably if these values are slightly 
modified (in other words, an asymptotic region is identified with 
$k\gg 1$ but a cutoff before reaching the level of rounding errors). 
A standard regression analysis of the relation (\ref{fourierasym}) 
leads to the matrix 
$$
   \mathcal{M} =  \sum_{k=k_{min}}^{k_{max}}
   \begin{pmatrix}
   	(\ln k)^{2} & k\ln k & \ln k\ \\
   	k\ln k & k^{2} & k\\
	\ln k & k & 1
   \end{pmatrix}
   $$
   and $\mathbf{\alpha}=\mathcal{M}^{-1}\mathbf{v}$, where 
   $$\mathbf{v} = \sum_{k=k_{min}}^{k_{max}}
     \begin{pmatrix}
   	(\ln(k))\hat{u}(k) \\ k\hat{u}(k) \\ \hat{u}(k)    
\end{pmatrix}
  $$
  and where the vector $\mathbf{\alpha}$ has the components $\begin{pmatrix}
   	\mu+1 \\ \delta \\ const\end{pmatrix}$.
The code is stopped when $\delta$ is smaller than the smallest distance on the numerical grid since smaller values 
cannot be numerically distinguished from zero in this case.  In the considered 
example this happens for $t\sim4.6809$. We get for the critical 
exponent of $\hat{v}$ the value $0.37$ and for the one for 
$\hat{\eta}$ the value $-0.85$. Note that the fitting of these 
critical exponents comes with a non negligible uncertainty, see the 
discussion in \cite{KR},  since it 
corresponds to fitting a logarithmic correction of a linear 
dependence in $k$ where the latter is close to vanishing. The data is 
compatible with $1/3$ for $v$, a cubic cusp as in the generic break-up of 
solutions to the Hopf equation, and with $-1$ for $\eta$, a simple 
pole. However, 
the behavior of the DFT coefficients clearly indicates that there is an $L^{\infty}$ blow-up in $\eta$ 
whereas a cusp forms in $v$. 

\begin{figure}[htb!]
 \includegraphics[width=0.49\textwidth]{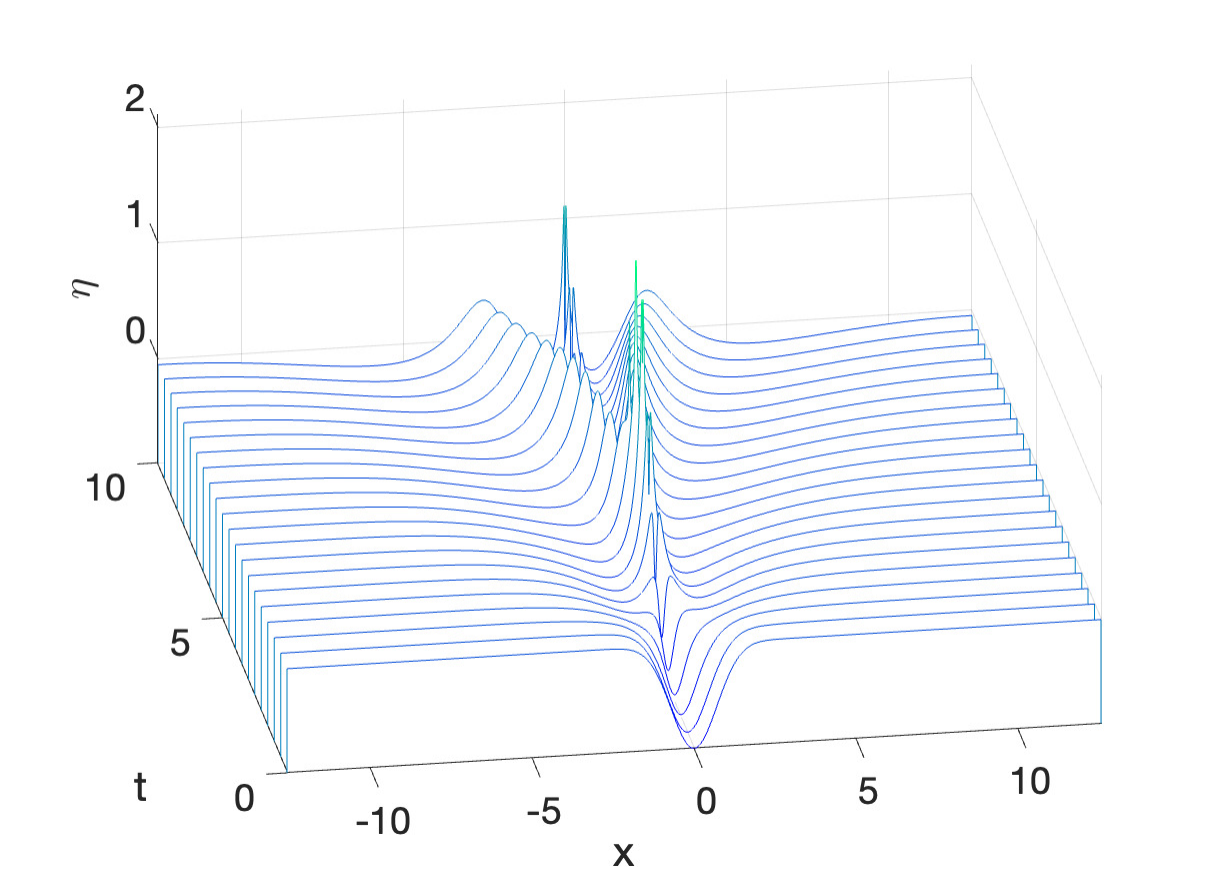}
 \includegraphics[width=0.49\textwidth]{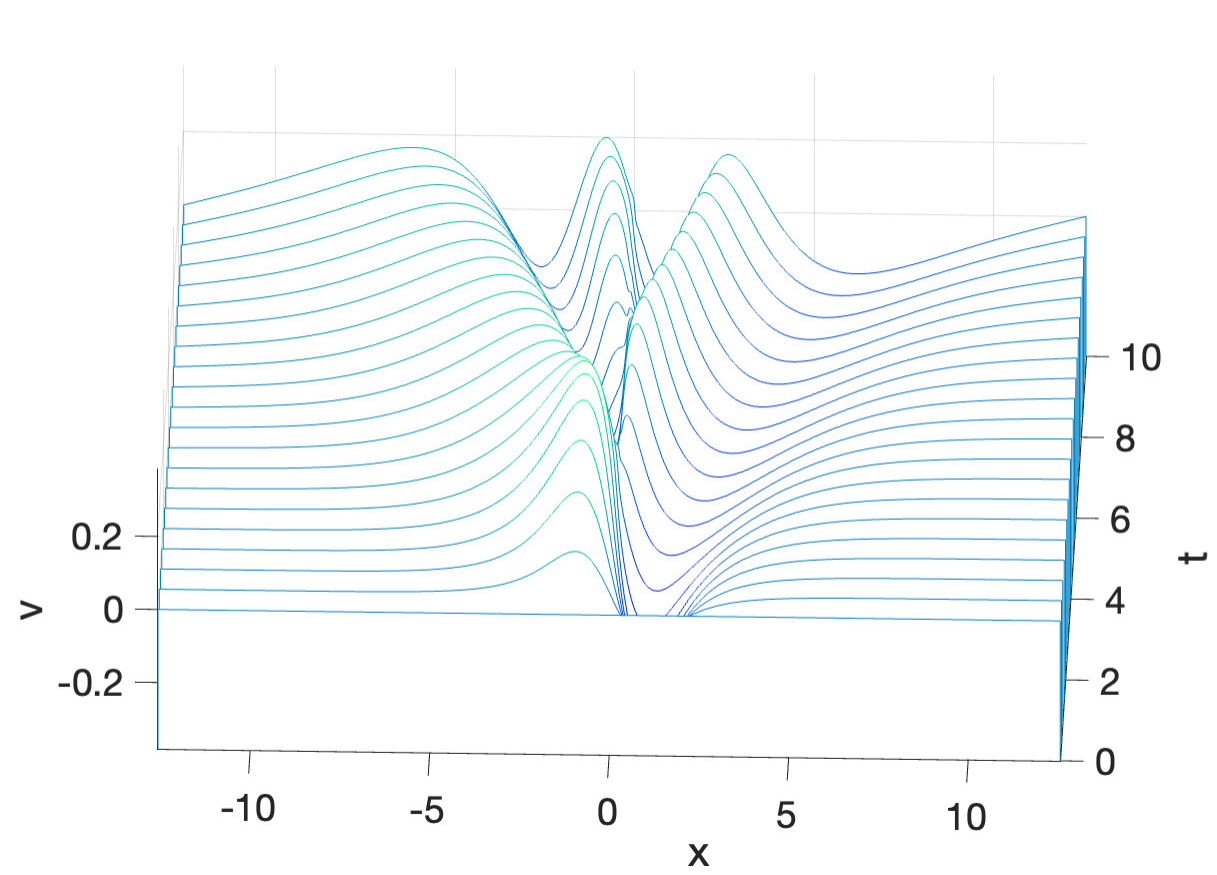}
 \caption{Solution to the Amick-Schonbeck system (\ref{AS}) for  
 initial data of the form (\ref{NC}) with $\kappa=0.9$, on the left $\eta$, on the right $v$.}
 \label{figASm09gauss}
\end{figure}

For initial data which come close to violating the non-cavitation 
condition, for instance (\ref{NC}) with $\kappa=0.9$, 
the solution shows qualitatively the same behavior as can be seen in 
Fig.~\ref{figASm09gauss}. But there is no indication of the formation 
of a singularity in this case. 

\begin{figure}[htb!]
 \includegraphics[width=0.49\textwidth]{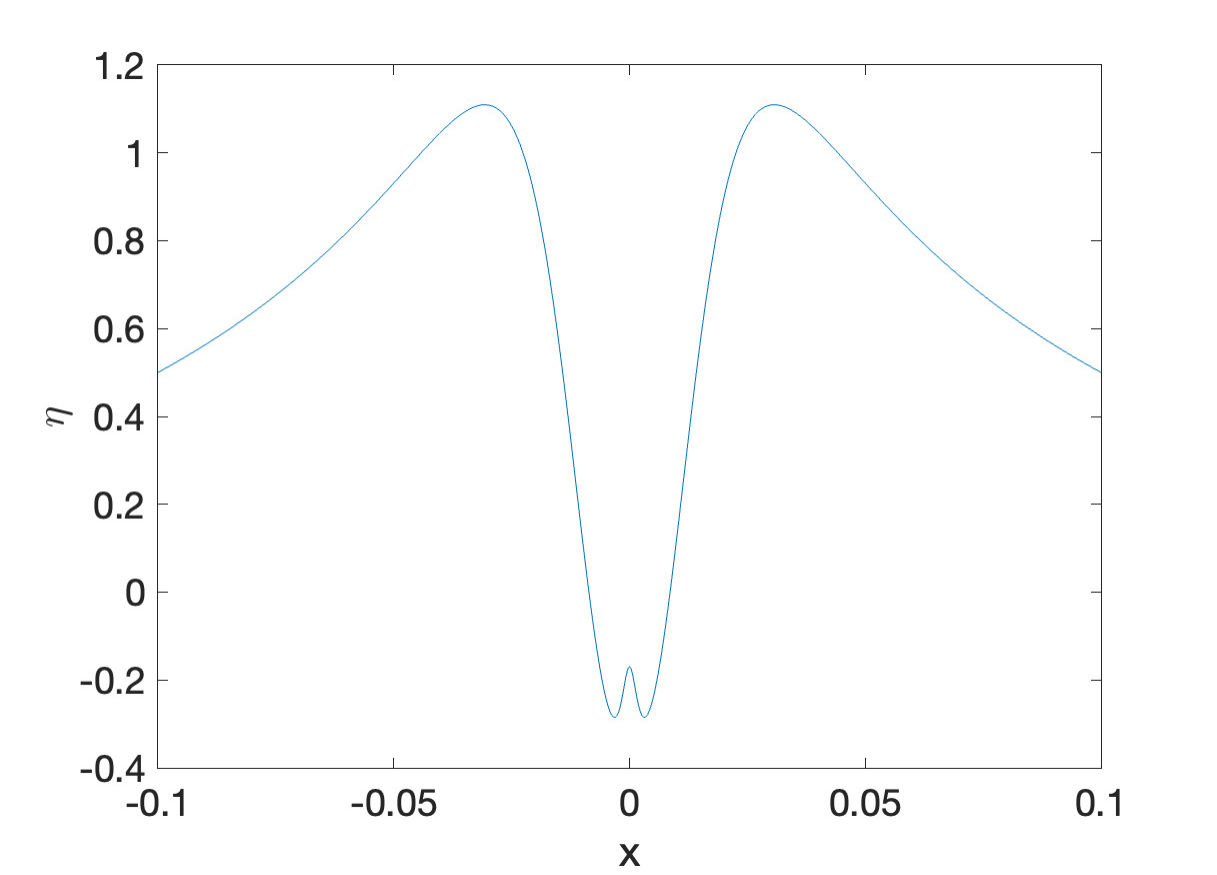}
 \includegraphics[width=0.49\textwidth]{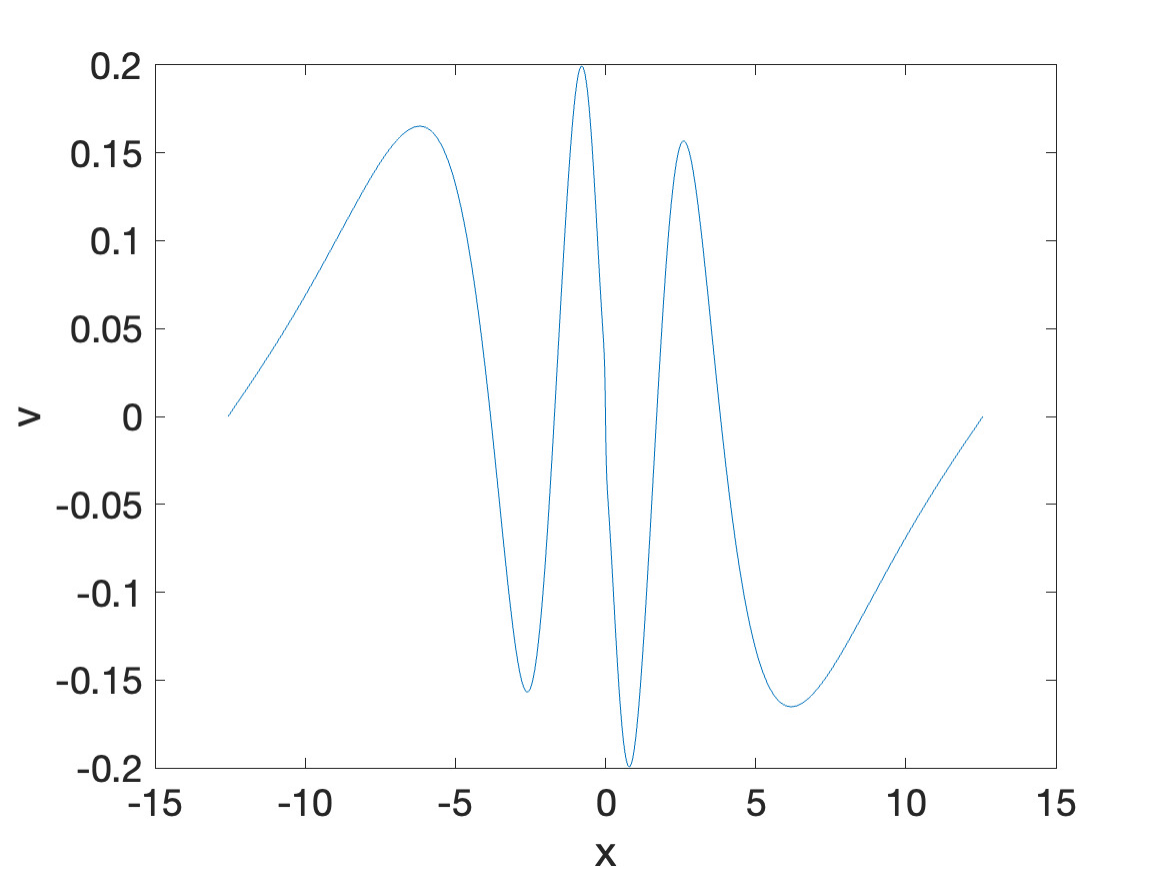}
 \caption{Solution to the Amick-Schonbeck system (\ref{AS}) for  
 initial data of the form (\ref{NC}) with $\kappa=0.9$ for $t=10$, on 
 the left a close-up of $\eta$, on the right $v$.}
 \label{figASm09gausst10}
\end{figure}

The DFT coefficients show exponential decrease during the whole 
computation. The solution for $t=10$ in Fig~\ref{figASm09gausst10} 
shows no indication of the formation of a singularity both in $v$ and 
in a close-up of $\eta$. The $L^{\infty}$ norms of both functions 
decrease with time for large $t$. 

\section{Dispersive shock waves}
In this section we study the possible formation of zones of rapid 
modulated oscillations called dispersive shock waves (DSWs). 

A convenient way to study the appearence of DSWs in nonlinear 
dispersive PDEs is to consider data with support of length scales of 
order $1/\varepsilon$, $\varepsilon\gg 1$, on time scales of order 
$1/\varepsilon$. This can be achieved by rescaling time and space by a 
factor of $\varepsilon$, which leads for the system (\ref{AS}) to a 
system containing a small parameter (we keep the same notation as in 
(\ref{AS}) for simplicity)
\begin{equation}
    \label{ASe}
    \left\lbrace
    \begin{array}{l}
    \eta_t+v_x+(\eta v)_x=0 \\
    v_t+\eta_x+vv_x-\varepsilon^{2}v_{xxt}=0.
\end{array}\right.
    \end{equation}
The formal limit $\varepsilon\to0$ leads to the Saint-Venant system
\begin{equation}
    \label{SV}
    \left\lbrace
    \begin{array}{l}
    \eta_t+v_x+(\eta v)_x=0 \\
    v_t+\eta_x+vv_x=0,
\end{array}\right.
    \end{equation}
expected to develop shocks for general initial data in finite time. 

\begin{figure}[htb!]
 \includegraphics[width=0.49\textwidth]{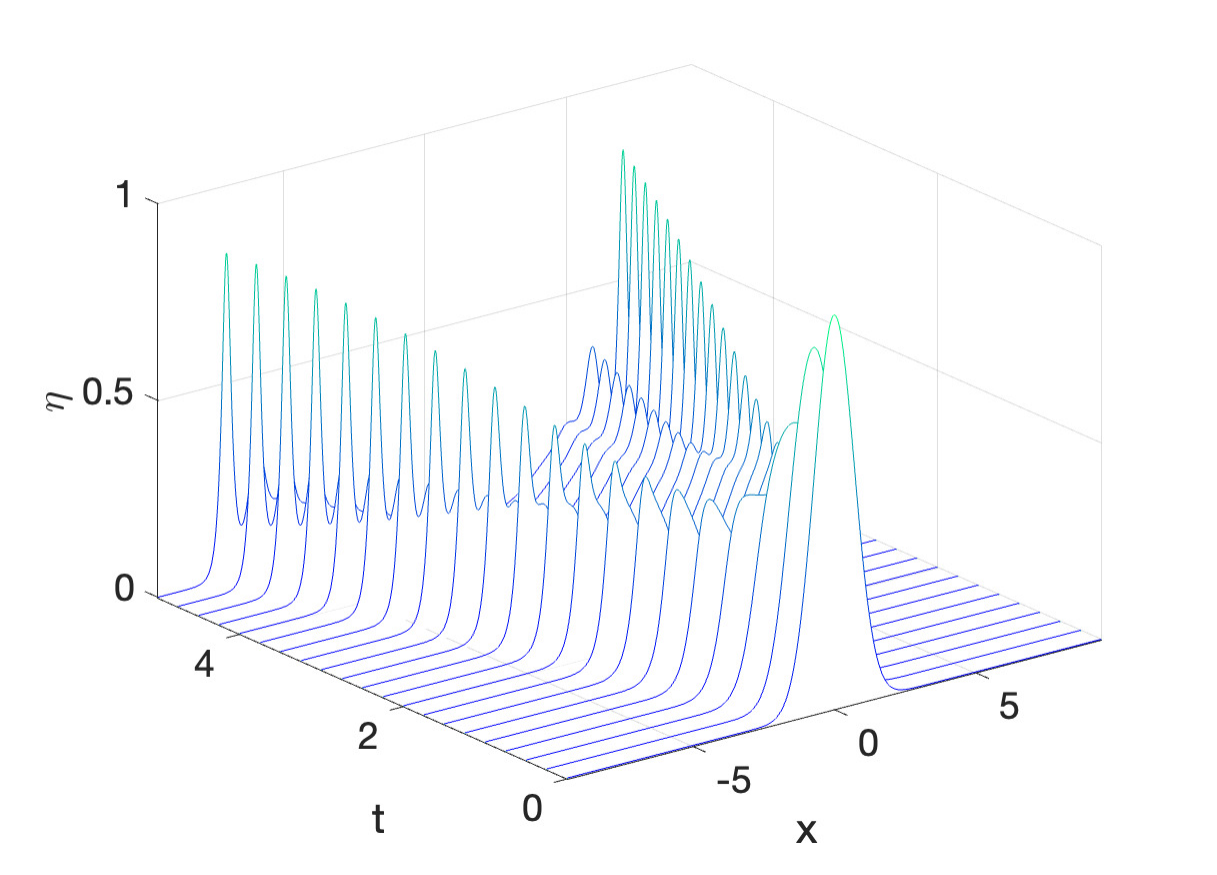}
 \includegraphics[width=0.49\textwidth]{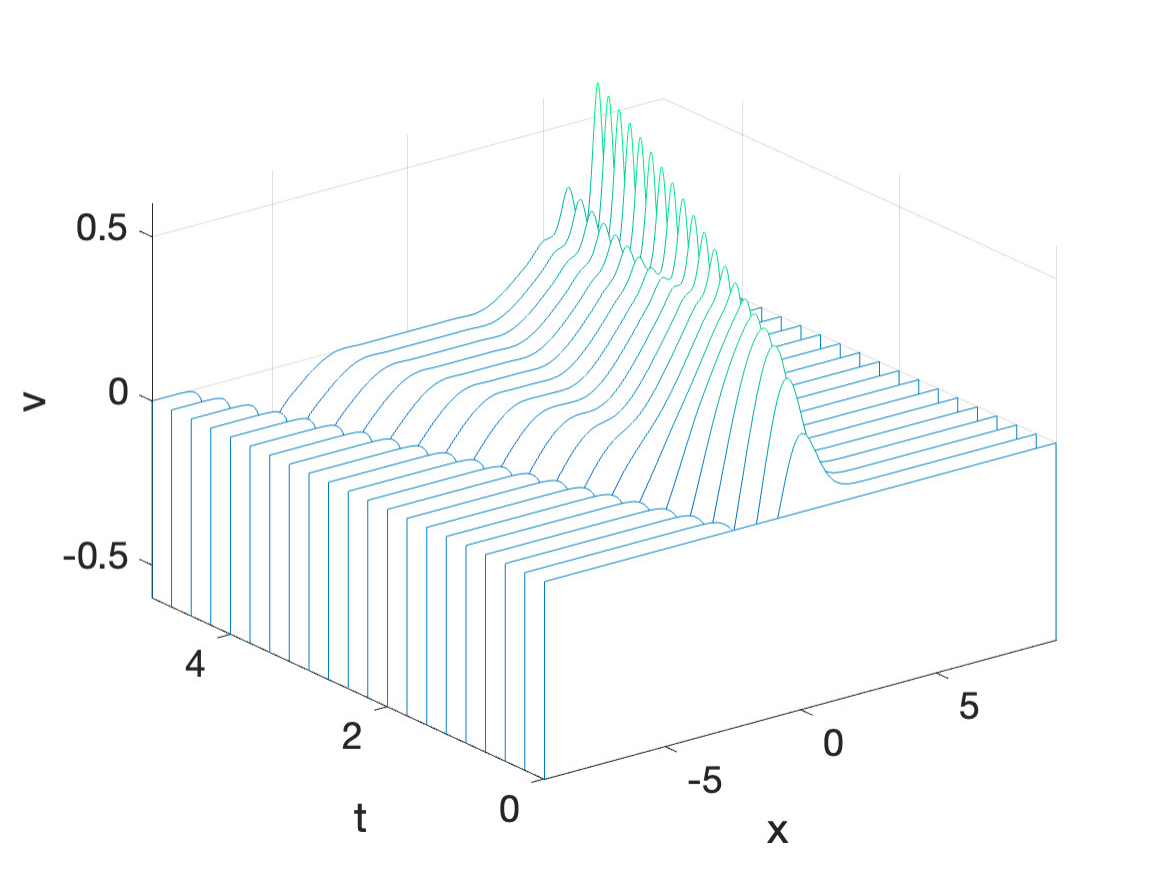}
 \caption{Solution to the Amick-Schonbeck system (\ref{ASe}) for 
 $\varepsilon=0.1$ and  
 initial data $\eta(x,0)=\exp(-x^{2})$, $v(x,0)=0$, on 
 the left  $\eta$, on the right $v$.}
 \label{Amickgauss1e2}
\end{figure}

We do not study this aspect here, but the dispersive regularisation 
of these shocks via (\ref{ASe}). We consider the initial data 
$\eta(x,0)=\exp(-x^{2})$, $v(x,0)=0$ as before for several values of 
$\varepsilon$ with $N=2^{14}$ DFT modes for $x\in3[-\pi,\pi]$ and 
$N_{t}=10^{4}$ time steps for $t\leq 5$. In Fig.~\ref{Amickgauss1e2} we show the 
solution for $\varepsilon=0.1$ 
in dependence of time. It can be seen that as before the function 
$\eta$ remains symmetric, and that two humps form. 
There is a steepening of the gradients near these humps which would 
form shocks in the Saint Venant system (\ref{SV}), but which
will lead to the formation of oscillations in the dispersive system 
(\ref{ASe}) near the location of the shocks of the Saint-Venant system. 
The function $v$ shows a 
similar behavior, but is odd. Note that the humps in the DSWs can be 
interpreted as solitary waves that look slightly different than 
before since both $x$ and $t$ have been rescaled. 

\begin{figure}[htb!]
 \includegraphics[width=0.49\textwidth]{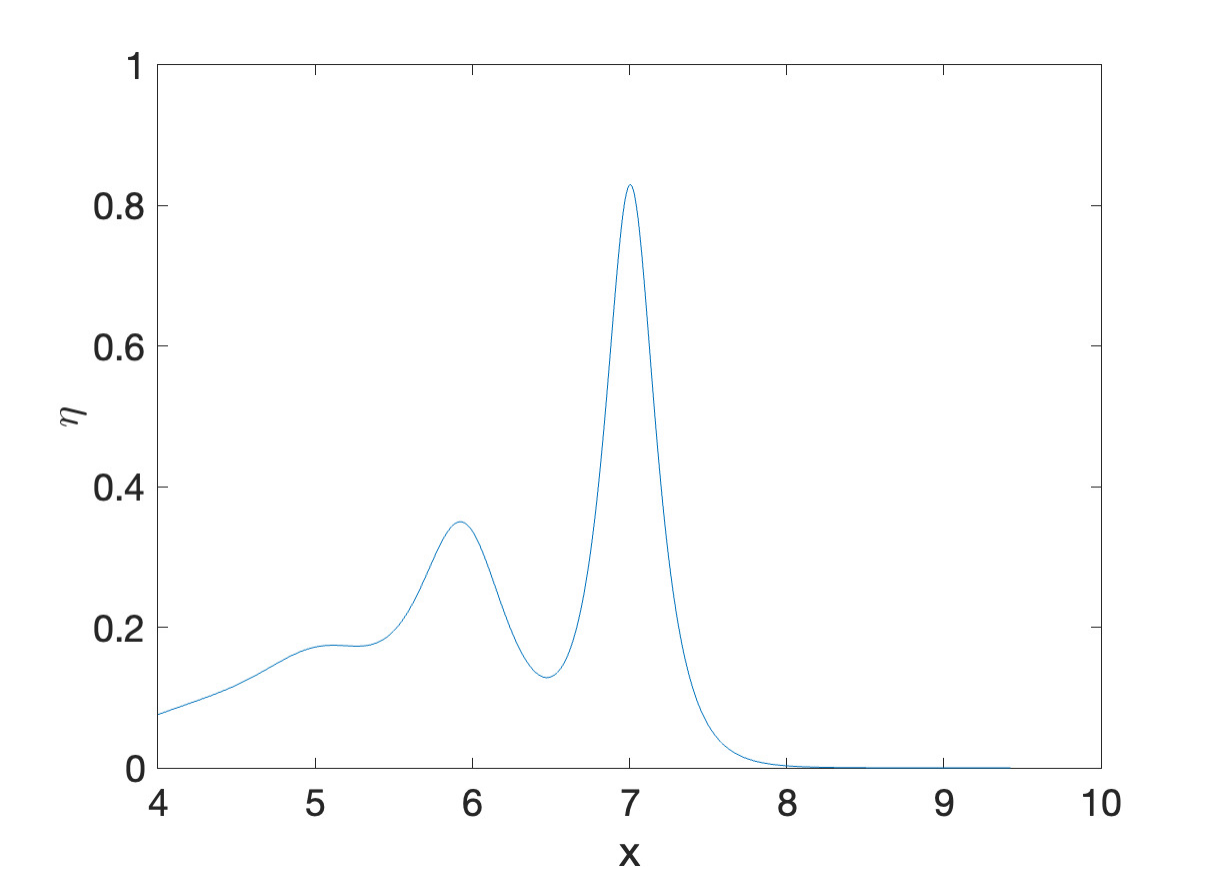}
 \includegraphics[width=0.49\textwidth]{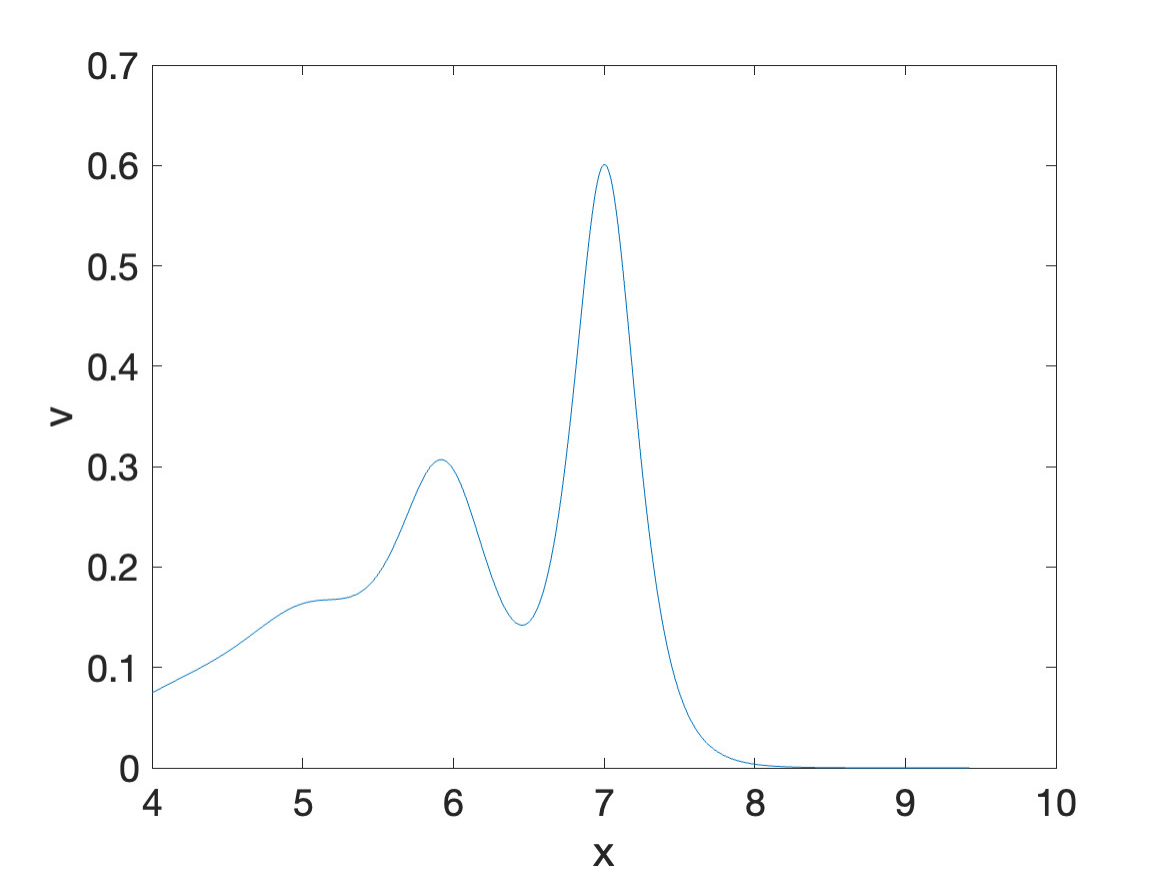}
 \caption{Close-up of the solutions to the Amick-Schonbeck system (\ref{ASe}) for 
 $\varepsilon=0.1$ and  
 initial data $\eta(x,0)=\exp(-x^{2})$, $v(x,0)=0$ 
for $t=5$, on 
 the left  $\eta$, on the right $v$.}
 \label{Amickgauss1e2cu}
\end{figure}

A close-up of the oscillatory zone for the solutions at the final 
time in Fig.~\ref{Amickgauss1e2} is shown in 
Fig.~\ref{Amickgauss1e2cu} supporting the soliton resolution 
conjecture. 

\begin{figure}[htb!]
 \includegraphics[width=0.32\textwidth]{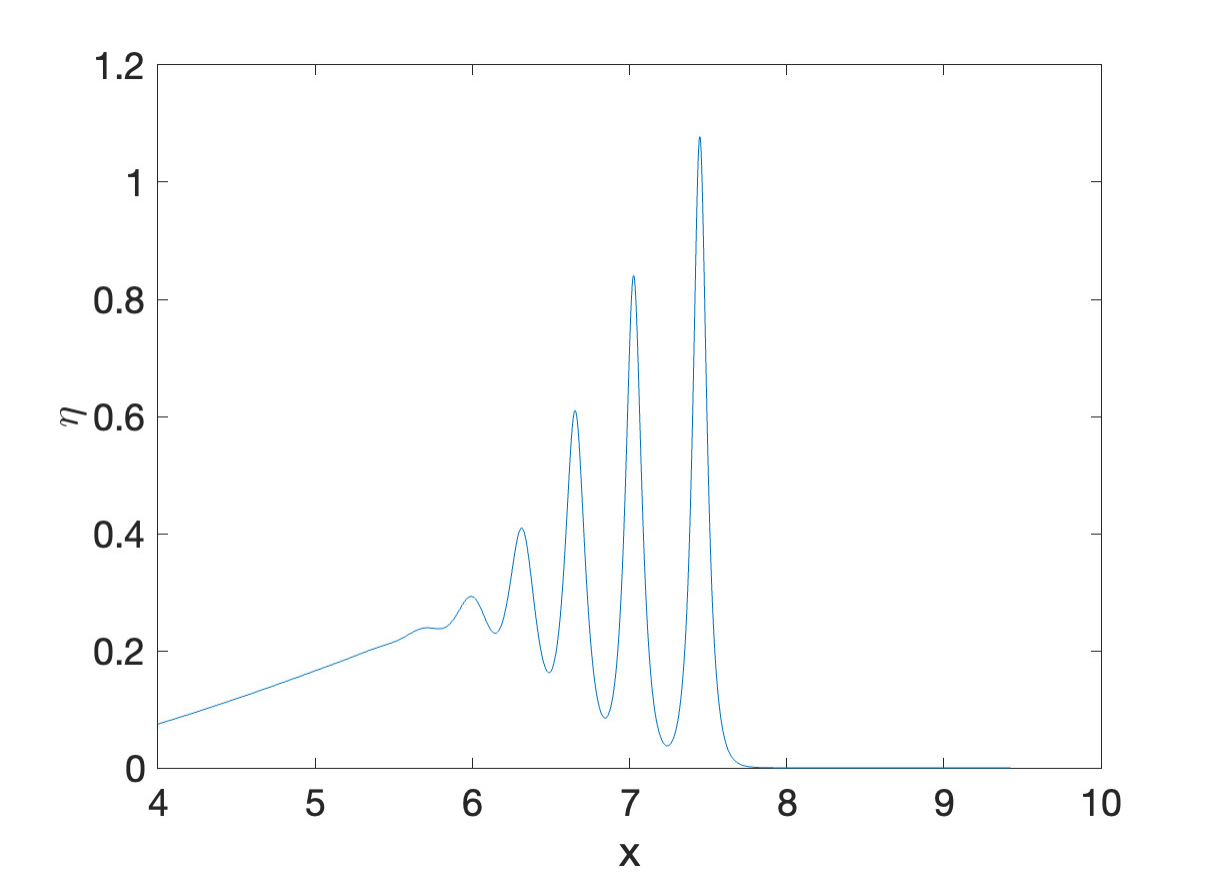}
 \includegraphics[width=0.32\textwidth]{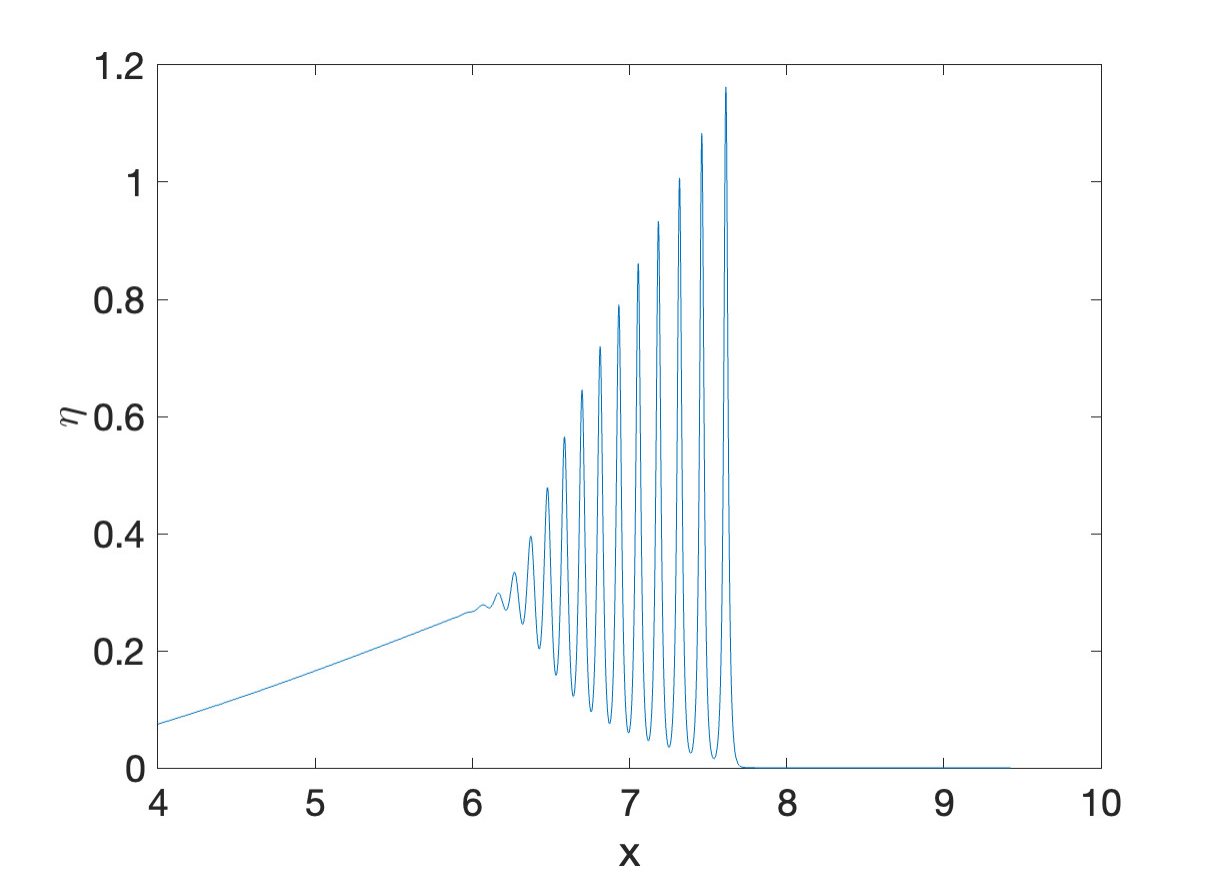}
 \includegraphics[width=0.32\textwidth]{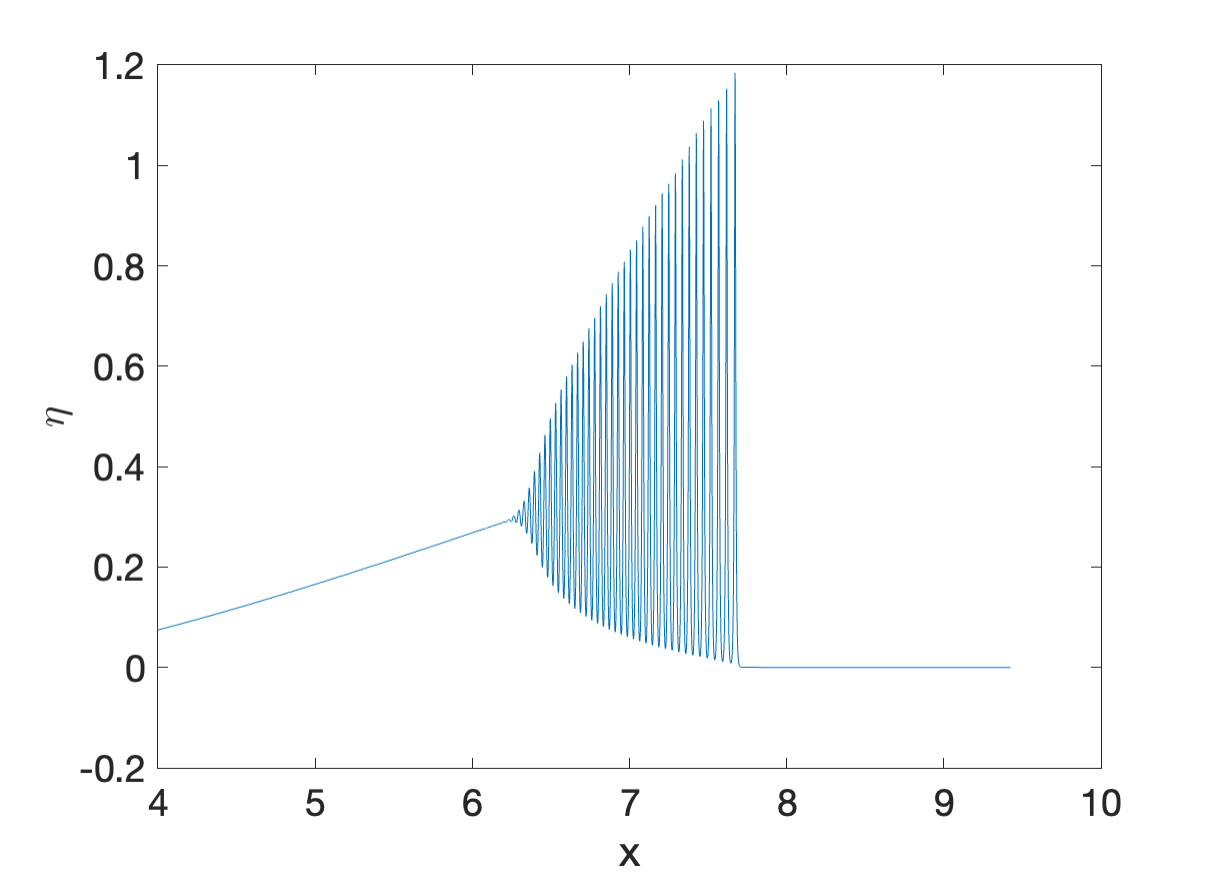}\\
 \includegraphics[width=0.32\textwidth]{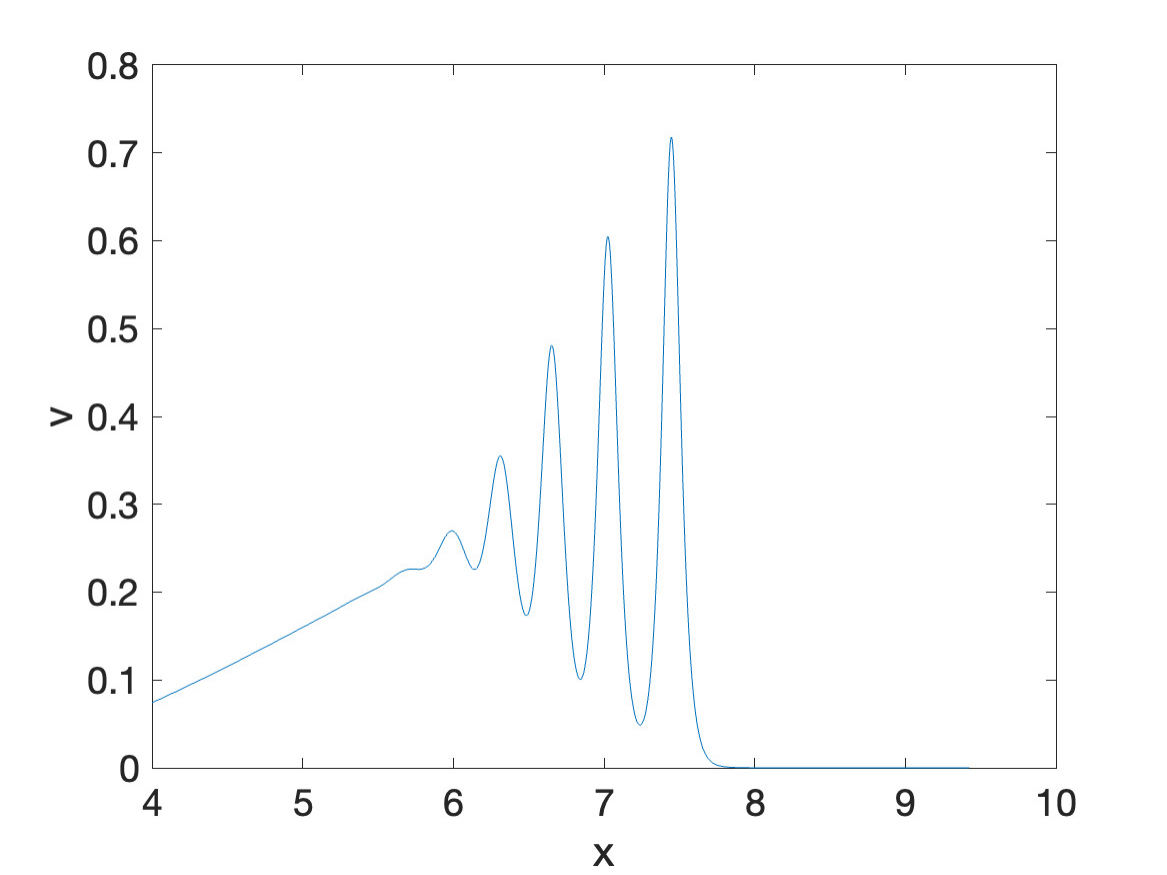}
 \includegraphics[width=0.32\textwidth]{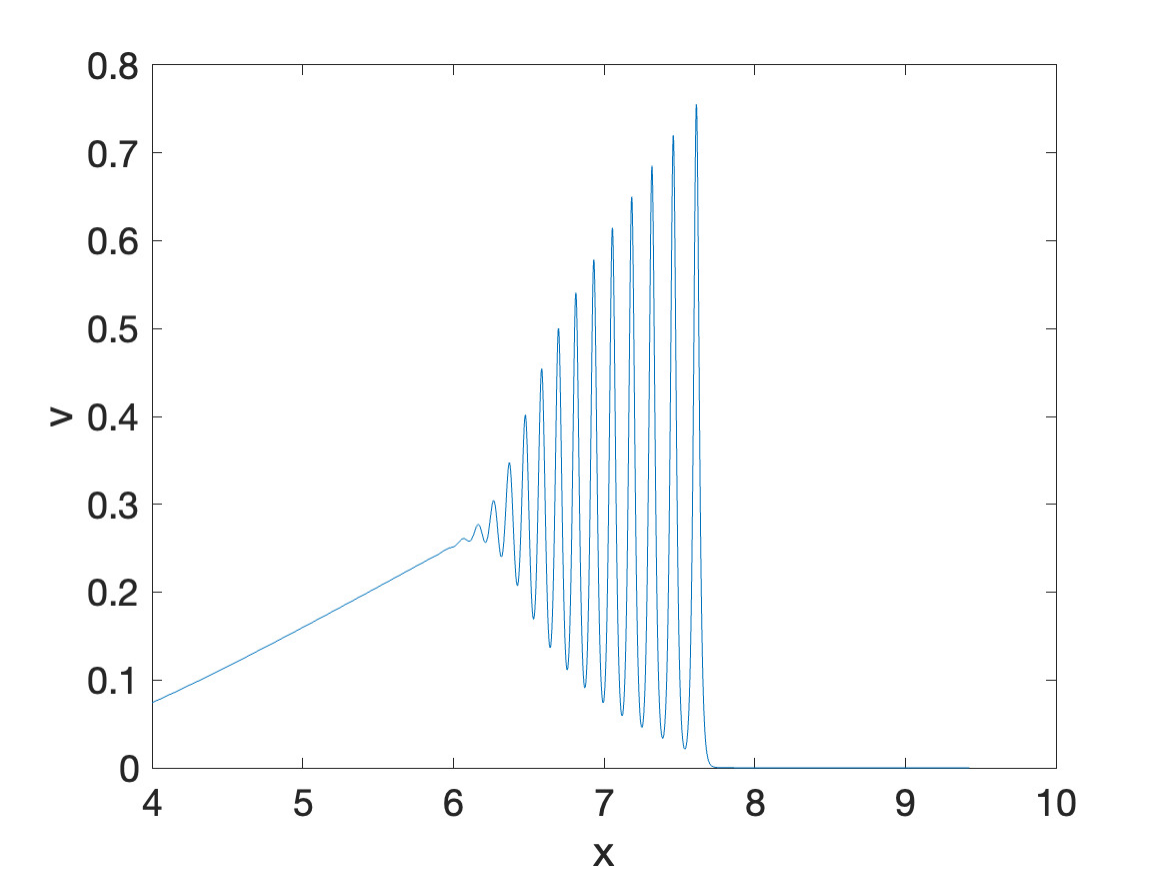}
 \includegraphics[width=0.32\textwidth]{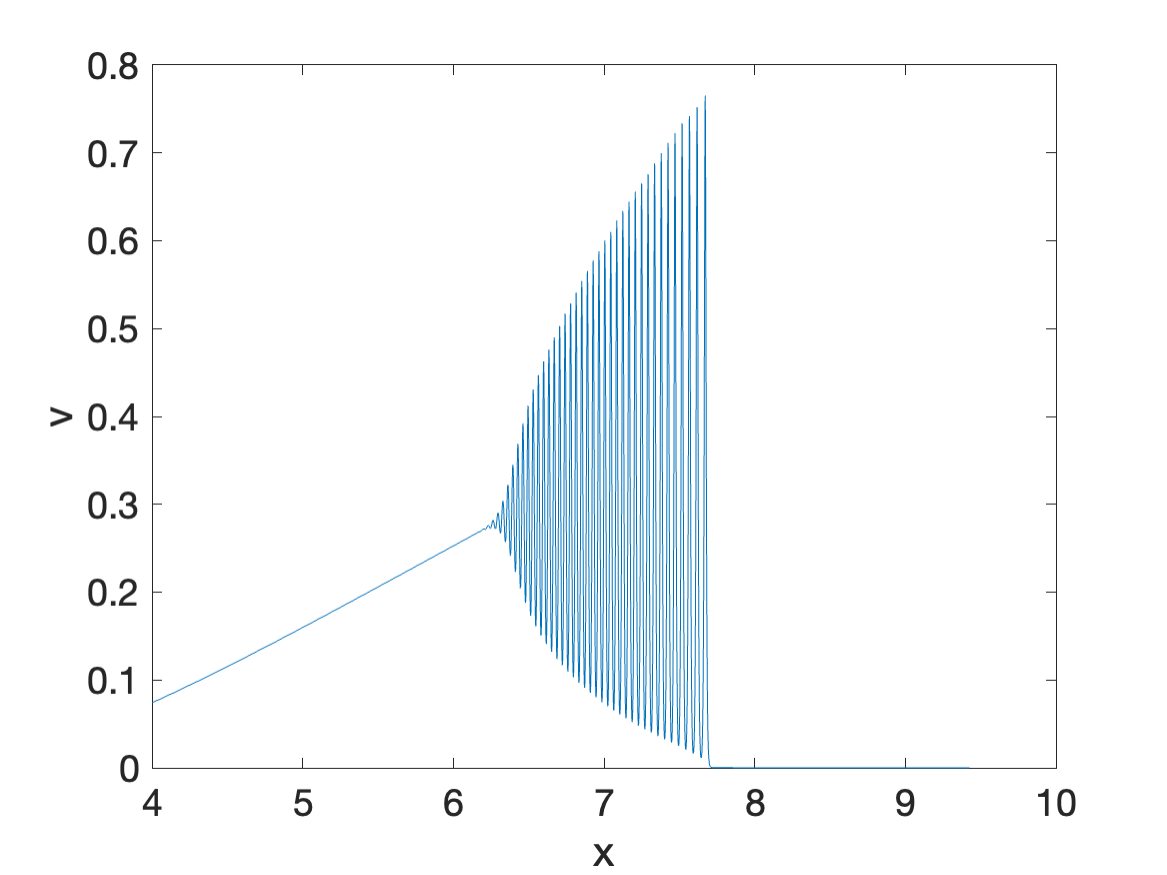}
 \caption{Close-up of the solutions to the Amick-Schonbeck system (\ref{ASe}) for 
 $\varepsilon=10^{-1.5}, 10^{-2},10^{-2.5}$ (from left to right) and  
 initial data $\eta(x,0)=\exp(-x^{2})$, $v(x,0)=0$  for $t=5$, upper row  
 $\eta$, lower row  $v$.}
 \label{Amickgaussecu}
\end{figure}

We show the solution for $t=5$ for several values of $\varepsilon$ in 
Fig.~\ref{Amickgaussecu} in a close-up of the oscillatory zone for 
positive $x$. It can be seen that the oscillations become more and 
more rapid for smaller $\varepsilon$, and that they are also more and 
more confined to a well defined zone sometimes referred to as Whitham 
zone. This behavior is reminiscent of
the well studied case of the KdV equation, see for instance 
\cite{KSbook}. 

\section{Conclusion}
In this paper we have numerically studied certain aspects of 
solutions to the Amick-Schonbeck system (\ref{AS}). The focus was on 
the asymptotic stability of the solitary waves for a wide range of 
the velocities $C$, and their appearence in the long time behavior of 
initial data from the Schwartz class. We have also addressed the 
appearence of DSWs that seem to be similar to DSWs in KdV 
solutions. In addition we have  established the appearence of a blow-up in 
solutions for initial data not satisfying the non-cavitation 
condition. 

It is an interesting question whether similar results can be obtained 
in two spatial dimensions. In \cite{GK} it was found that the 2D 
Serre-Green-Naghdi equations show a defocusing effect similar to the 
Kadomtsev-Petviashvili II equation, see \cite{KSbook}. No solitary 
waves localised in 2D were found. The question is whether there are 
such solutions for the 2D Amick-Schonbeck system (\ref{2D}). If they 
exist, there stability has to be studied as well as the long time 
behavior of solutions for localised initial data. The possibility of 
a blow-up in 2D has to be explored. This will be the topic of  future 
work. 
 
\begin{merci}

CK was partially supported by 
 the ANR-17-EURE-0002 EIPHI and by the 
European Union Horizon 2020 research and innovation program under the 
Marie Sklodowska-Curie RISE 2017 grant agreement no. 778010 
IPaDEGAN.
\end{merci}

\bibliographystyle{amsplain}

\end{document}